\documentstyle[12pt,amssymb]{article}
\newtheorem{Th}{Theorem}

\newtheorem{Lm}{Lemma}
\newtheorem{Lma}{Lemma}[section]

\newtheorem{Rm}{Remark}
\newcommand{\be}{\begin{equation}}
\newcommand{\ee}{\end{equation}}
\newcommand{\R}{\mathbb{R}}
\newcommand{\N}{\mathbb{N}}
\newcommand{\C}{\mathbb{C}}
\newcommand{\Z}{\mathbb{Z}}
\newcommand{\cqfd}
{%
\mbox{}%
\nolinebreak%
\hfill%
\rule{2mm}{2mm}%
\medbreak%
\par%
}
\newcommand\res{\mathop{\hbox{\vrule height 7pt width .5pt depth 0pt
\vrule height .5pt width 6pt depth 0pt}}\nolimits}

\newcommand{\reset}{\setcounter{equation}{0}\setcounter{Th}{0}\setcounter{Prop}{0}\setcounter{Co}{0}
\setcounter{Lm}{0}\setcounter{Rm}{0}}
\def\La{\Lambda}
\def\ti{\tilde}
\def\lf{\left}
\def\rg{\right}
\def\noi{\noindent}
\def\al{\alpha}
\def\la{\lambda}

\def\ep{\varepsilon}
\def\ds{\displaystyle}
\def\ov{\overline}
\def\Om{\Omega}
\def\om{\omega}
\def\p{\partial}
\def\vphi{\varphi}
\def\Txi{{{\ti{\Xi}}^{ij}}}
\begin{document}
\title{The Singular Set of 1-1 Integral Currents.}
\author{Tristan Rivi\`ere and Gang Tian}
\date{ }
\maketitle
{\bf Abstract }: {\it We prove that 2 dimensional Integral currents 
(i.e. integer multiplicity 2 dimensional rectifiable currents) which are almost complex cycles in an
almost complex manifold admitting locally a compatible symplectic form are
smooth submanifolds aside from isolated  points and therefore are
$J-$holomorphic curves. } 

\section{Introduction}

Let $(M^{2p},J)$ be an almost complex manifold. Let $k\in {\N}$, $k\le p$.
We shall adopt classical notations from Geometric Measure Theory. We say that a
 $2k-$current $C$ in $(M^{2p},J)$ is an
 almost complex integral cycle whenever it fulfills the following three conditions
\begin{itemize}
\item[i)] {\it rectifiability :} there exists an at most countable union of disjoint oriented
$C^1$ $2k-$submanifolds ${\mathcal C}=\cup_{i} N_i$ and an integer multiplicity
$\theta\in L^1_{loc}({\mathcal C})$ such that for any smooth compactly supported in $M$ $2k-$form
$\psi$ one has
\[
C(\psi)=\sum_{i}\int_{N_i}\theta\ \psi\quad.
\]
\item[ii)] {\it closedness :} $C$ is a cycle
\[
\p C=0\quad\mbox{ i. e. }\forall \al\in C_0^\infty(\wedge^{2k-1}M)\quad C(d\al)=0
\]
\item[iii)] {\it almost complex :} for ${\mathcal H}^{2k}$ almost every point $x$ in ${\mathcal C}$
the approximate tangent plane $T_x$ to the rectifiable set ${\mathcal C}$ is invariant under 
the almost complex structure $J$ i.e.
\[
J(T_x)=T_x
\]
\end{itemize}
In this work we address the question of the regularity of such a cycle : Does there exist a {\it smooth}
almost complex manifold $(\Sigma^{2k},j)$ without boundary and a {\it smooth} $j-J$-holomorphic map $u$
($\forall x\in \Sigma$ and $\forall X\in T_x\Sigma$ $du_xj\cdot X=J\cdot du_xX$) such that $u$ would
realise an {\it embedding} in $M^{2p}$ aside from a locally finite $2k-2$ measure closed subset of $M$
and such that $C=u_\ast [\Sigma^{2k}]$, i.e. $\forall \psi\in C^\infty_0(\wedge^{2k}M)$
\[
C(\psi)=\int_\Sigma u^\ast\psi
\]
In the very particular case where the almost complex structure $J$ is integrable, this regularity result
 is optimal ($C$ is the integral over multiples of algebraic subvarieties of $M$) and was established in 
\cite{HS}. 
There are numerous motivations for studiying the general case of arbitrary almost complex structures $J$.
First, as explained in \cite{RT}, the above regularity question for rectifiable almost complex
cycles is directly connected to the regularity question of $J-$holomorphic maps into complex projective
 spaces. It is conjectured, for instance, that the singular set of $W^{1,2}(M^{2p},N)$ $J-$holomorphic maps
between almost complex manifolds $M$ and $N$ should be of finite $(2p-4)$-Hausdorff measure. The resolution of that question leads, for instance, to the caracterisation of stable bundle almost complex structures
over almost K\"ahler manifolds via Hermitte-Einstein Structures and extends Donaldson, Uhlenbeck-Yau
characterisation in the integrable case (see \cite{Do}, \cite{UY}) to the non-integrable one.  
Another motivation for  studying the regularity of almost complex rectifiable cycles is the following.
In \cite{Li} and \cite{Ti} it is explained how the loss of compactness of
solutions to geometric PDEs having a given conformal invariant dimension $q$ (a dimension at which the
PDE is invariant under conforormal transformations - $q=2$ for harmonic maps, $q=4$ for Yang-Mills
Fields...etc) arises along $m-q$ rectifiable cycles (if $m$ denotes the dimension of the domain).
This cycles happen sometimes to be almost complex (see more details in \cite{Ri1}). 

By trying to produce in $({\R}^{2p}, J)$ an almost complex graph of real dimension $2k$ in a
 neighborhood of a point $x_0\in {\R}^{2p}$ as a perturbation
of a complex one ($J_{x_0}-$holomorphic), one realises easily that, for generic almost complex
 stuctures $J$,
 the problem is overdeterminate whenever $k>1$ and well posed for $k=1$. Therefore the case of 
$2-$dimensional integer rectifiable almost complex cycles is the generic one from the existence point 
of view. We shall restrict to that important case in the present paper. After complexification
of the tangent bundle to $M^{2p}$ a classical result asserts that a $2-$plane is invariant under $J$
if and only if it has a $1-1$ tangent $2-$vector. Therefore we shall also speak about $1-1$ integral
 cycles for the almost complex $2-$dimensional integral cycles.  
In the present work we consider the locally symplectic case : We say that $(M^{2p},J)$ has the {\it locally
symplectic property} if  at a neighborhood
of each point $x_0$ in $M^{2p}$ there exists a symplectic structure compatible with $J$ : there
exists a neighborhood $U$ of $x_0$ and a smooth 2-form $\om$ satisfying 
\begin{itemize}
\item[i)]
\[
d\om=0 \quad\quad\mbox{ in }U
\]
\item[ii)]
\[
\om^p(x_0)\ne 0
\]
\item[iii)]
\[
\om(\cdot,\cdot)=\om(J\cdot,J\cdot)
\]
\end{itemize}
Observe that iii) is equivalent to the fact that $\om(\cdot,J\cdot)$ is a symmetric form and together
with ii) it says that $\om(\cdot, J\cdot)$ is a scalar product.  
The local existence of such a symplectic form for arbitrary $J$ is still an open problem in dimension larger than 4. It was proved in \cite{RT} that arbitrary 4-dimensional almost complex manifolds fulfill
the {\it locally symplectic property}. This should not be the case in larger dimension anymore but to
our knowledge in 6 dimension or higher no $J$, which do not admit a locally symplectic structure,
has been exhibited.

Our main result is the following.

\begin{Th}
\label{th-I.1}
Let $(M^{2p},J)$ be an almost complex manifold satisfying the locally symplectic property above.
Let $C$ be an integral 2 dimensional almost complex cycle. Then, there exists $J$-holomorphic
curve $\Sigma$ in $M$, smooth aside from isolated points, and a smooth integer valued function $\theta$
on $\Sigma$ such that, for any 2 form $\psi\in C^\infty_0(M)$,
\[
C(\psi)=\int_\Sigma\theta\ \psi\quad .
\]
\end{Th}
In the ``locally symplectic case'' being an  almost-complex 2 cycle is equivalent for a 2-cycle of being
 calibrated by the local symplectic form $\om$ for the local metric $\om(\cdot, J\cdot)$. Therefore
the regularity question for almost complex cycles is embedded into the problem of calibrated current
and hence the theory of area minimizing
rectifiable 2-cycles. Therefore our result appears to be a consequence of the ``Big Regularity Paper''
of F.Almgren \cite{Alm} combined with the PhD thesis of his student S.Chang \cite{Ch}.
Our attempt here was to present an alternative proof independent of Almgren's monument \cite{Alm}
and adapted to the case we are interrested with. The motivation is to give a proof that could be modified
in order to solve the general case (non locally symplectic one) which cannot be ``embedded'' in the theory
of area minimizing cycles anymore. In \cite{PR} we explain how far the general case, without
necessarily the existence of a locally symplectic form, is a perturbation of the {\it locally symplectic case}
by terms of lower orders which have no significant influence on the strategy of the proof presented here
while extending theorem~\ref{th-I.1} to the general case. 
This attempt of writing a proof for the regularity of almost complex cycle in the locally symplectic case, 
independent of the regularity theory for area minimizing surfaces,
was also the main purpose of the work Gr$\Longrightarrow$SW of C.Taubes \cite{Ta} for $p=2$. 
In this work a proof of theorem~\ref{th-I.1} was given in the particular case where $p=2$.
Some argument happened to be incomplete in that proof and \cite{RT} fills the missing step in \cite{Ta}.
Theorem~\ref{th-I.1} has to be seen as the generalisation to higher dimension $(p>2)$ of these works.

One of the main difficulties arising in dimension $(p>2)$ is the non-necessary existence of
$J-$holomorphic foliations transverse to our almost complex  current $C$ in a neighborhood of a
point. This prevents then to describe the current as a $Q-$multivalued graph from $D^2$ into ${\C}^{p-1}$,
 $\{(a_i^k(z))_{k=1^{p-1}}\}_{i=1\cdots Q}$ in a neighborhood of a point of
density $N$ solving locally an equation of the form 
\be
\label{I.01}
\p_{\ov{z}}a_i^k= \sum_{l=1}^{p-1}A(z,a_i)^k_l\cdot \nabla a_i^l+\al^k(a_i,z)\quad,
\ee
where $A$ and $\al$ are small in $C^2$ norm, as we did for $p=2$ in \cite{RT}. What we can only ensure instead is to describe the current $C$,
in a neighborhood of a point of multiplicity $Q$, as a ``algebraic $Q-$valued graph'' from
$D^2$ into ${\C}^{p-1}$ : that is a familly of points in ${\C}^{p-1}$,
 $\{ a_1(z),\cdots ,a_P(z), b_1(z),\cdots,b_N(z)\}$ where only $P-N=Q$ is independent on $z$
(neither $P$ nor $N$ are a-priori independent on $z$), 
$a_i$ are the positive
intersection points and $b_j$ are the negative ones. This ``algebraic $Q-$valued graph'' solve
locally a much less  attractive equation than (\ref{I.01})
\be
\label{I.02} 
\p_{\ov{z}}a_i^k=\sum_{l=1}^{p-1}{A}^l_k(z,a_i,\nabla a_i)\cdot\nabla a_i^l +\sum_{l=1}^{p-1} B^l_k(z,a_i)
\cdot\nabla a_i^l + C^k(z,a_i)\quad .
\ee
where ${A}(z,a,p)$, $B(z,a)$ and $C(z,a)$ are also small in $C^2$ norm but the dependence in $p$ in $A(z,a,p)$
 is linear and therefore
as $\nabla a_i$ gets bigger, which can happen, the right-hand-side of (\ref{I.02}) can not be handled as
a perturbation of the left-hand one in steps such as the ``unique continuation argument'' which was
used in \cite{RT} for proving that singularities of mutiplicity $Q$ cannot have accumulation point
in the carrier ${\mathcal C}$ of $C$. 
 
The strategy of the proof goes as follows. A classical blow-up analysis tells us that  for an arbitrary point $x_0$ of the manifold
$M^{2p}$  the limiting density $\theta(x_0)=lim_{r\rightarrow 0}r^{-2}M(C\res B_r(x_0))$ - $M$ denotes the mass of a current and $\res$ is 
the restriction operator -  equal $\pi$ times an integer $Q$. Since the density function $r\rightarrow r^{-2}M(C\res B_r(x_0))$
 at every point is an increasing function, the complement of the set ${\mathcal C}_Q:=\{x\in M\ ;\ \theta(x)\le Q\}$ is closed
in $M$ and this permits to perform an inductive proof of theorem~\ref{th-I.1} restricting the current to ${\mathcal C}_Q$ 
and considering increasing integers $Q$. A point of multiplicity $Q$ is called a singular point of $C$ if it is in the closure
of points of non zero multiplicity strictly less than $Q$. The goal of the proof is then to show that singularities of multiplicity less than $Q$
are isolated. We assume this fact for $Q-1$ and the paper is devoted to the proof that this then holds for $Q$ itself.
From a now classical result of B. White (see \cite{Wh}), the dilated currents at a point
$x_0$ of density $Q\ne 0$ converge in flat norm to a sum of $Q$ flat $J_{x_0}-$holomorphic disks, moreover,
for any $\ep>0$ and $r$ sufficiently small $C\res B_r(x_0)$ is supported in the cones whose axis are the limiting disks
and angle $\ep$. For $Q>1$, if two of this limiting disks are different it is then easy to observe that $x_0$ cannot
be an accumulation point of singularities of multiplicity $Q$, this is the so called ``easy case''. If the limiting
disks are all identical, equal to $D_0$, then we are in the ``difficult case'' and much more work has to be done in order to reach the same
statement. Countrary to the special case of the 4-dimension (p=2) considered by the authors before in \cite{RT},
we could not find nice coordinates that would permit to write $C$ as a $Q-$valued graph over the limiting
disk $D_0$. Considering then some $J_{x_0}-$complex coordinates $(z,w_1,\cdots,w_{p-1})$ 
in a neighborhood $B^{2p}_{\rho_{x_0}}(x_0)$ such that $\cap_iw_i^{-1}\{0\}$ corresponds to $D_0$, by the mean of the
``lower-epiperimetric inequality'' proved by the first author in \cite{Ri2}, one can construct a Whitney-Besicovitch
covering, $\{B^2_{\rho_i}(z_i)\}_{i\in I}$, of the orthogonal projection on $D_0$ of the points in $B^{2p}_{\rho_{x_0}}(x_0)$
having a positive density strictly less than $Q$. This covering is such that for every $i\in I$ there exists 
$x_i=(z_i,w_i)\in B^{2p}_{\rho_{x_0}}(x_0)$ verifying that the restriction of $C$ to the tube $B^2_{\rho_i}(z_i)\times B^{2p-2}_{\rho_{x_0}}(0)$ is
in fact supported in the ball $B_{\rho_i}^{2p}(x_i)$ of radius $\rho_i$, the width of the tube. Moreover if one looks
inside $B_{\rho_i}^{2p}(x_i)$, $C$ is ``splitted'' : this last word means that $C$ restricted to $B_{\rho_i}^{2p}(x_i)$ is at a flat distance comparable
 to $\rho_i^3$ from the $Q$ multiple of any  graph 
over $B^2_{\rho_i}(z_i)$ - this comes from the fact that the density $\rho_i^{-2}M(B_{\rho_i}(x_i))$ is strictly less
than $\pi Q$ minus a constant $\al$ depending only on $p$, $Q$, $J$ and $\om$  .
We then construct an average curve for $C$. In the 4-dimensional case since $C$ was a $Q-$valued graph over $D_0$ we took simply
the average of the $Q$ points over any point in $D_0$. Here, in arbitrary dimension, the construction of the average curve is 
more delicate and uses the covering. We first approximate $C\res B_{\rho_i}^{2p}(x_i)$ by a $J_{x_i}-$holomorphic graph $C_i$
using a technique introduced in \cite{Ri3},
and choosing a $J_{x_i}-$holomorphic disk $D_i$ approximating $D_0$ we can express $C_i$ as a $Q-$valued graph over $D_i$ for
which we take the average $\ti{C}_i$ that happens to be Lipschitz with a uniformly bounded lipschitz constant. Therefore
the $J_{x_i}-$holomorphic curve $\ti{C}_i$ can be viewed as a graph $\ti{a}_i$ over $B^2_{\rho_i}(z_i)$. Patching the $\ti{a}_i$
together we get a graph  $\ti{a}$ that extends over the whole $B^{2}_{\rho_{x_0}}(0)$ as a $C^{1,\al}$ graph for any $\al<1$
which is almost $J-$holomorphic and which passes through all the $B_{\rho_i}^{2p}(x_i)$. The fact that the average curve 
is more regular than the $J-$holomorphic cycle $C$ from which it is produced is clear in the integrable case (since it it holomorphic)
- $(z,\pm\sqrt{z})$ is a $C^{0,\frac{1}{2}}$ 2-valued graph whereas its average $(z,0)$ is smooth. This was extended in the
non-integrable case in the particular case of the 4 dimension in \cite{ST}. The points of multiplicity $Q$ in $C$ are contained in
the average curve $\ti{a}$. We then show, by the mean of a unique continuation argument in the spirit of the one developped in \cite{Ta}
in 4 dimension, that the points where $C$ get to coincide with $\ti{a}$ are either isolated or coincide with the whole curve
$\ti{a}$. We have then showed that any point $x_0$ of multiplicity $Q$ is either surrounded by points of multiplicity $Q$ only, and in $B^{2p}_{\rho_{x_0}}$
$C$ coincides with $Q$ time a smooth graph over $D_0$ or $x_0$ is not an accumulation point of points of multiplicity
$Q$ and is surrounded in $B^{2p}_{\rho_{x_0}}$ by points of multiplicity strictly less than $Q$. It remains at the end to show that it 
cannot be an accumulation point of singularities of lower density. This is obtained
again using an approximation argument by holomorphic curves introduced in \cite{Ri3}.

The paper is organised as follows. In chapter II we establish preliminaries, introduce notations and give
the main statement, assertion ${\mathcal P}_Q$, we are going to prove by induction  in the rest of the paper.
 In chapter III, with the help of the ``upper-epiperimetric inequality'' of B. White, we establish 
the uniqueness of the tangent cone and a quantitative version of it, see lemma~\ref{lm-III.2}.
 In chapter IV we prove the relative lipschitz estimate together with a tilting control of the 
tangent cones of density $Q$ points in a neighborhood of a density $Q$ point - see lemma~\ref{lm-IV.2}.
 In chapter V we proceed to the covering argument lemma~\ref{lm-V.2} which is based on
the ``splitting before tillting'' lemma - see lemma~\ref{lm-IV.3} - proved in \cite{Ri2}.
In chapter VI we construct the approximated average curve and prove the $C^{1,\al}$ estimate
for this curve - lemma~\ref{lm-VII.2}. In chapter VII we perform the unique continuation argument that
shows that singularities of multiplicity $Q$ cannot be an accumulation point of singularities
of multiplicity $Q$. In chapter VIII we show that singularities of multiplicity
$Q$ cannot be an accumulation point of singularity of multiplicity less than $Q$ either. 

\section{Preliminaries}
\reset
{\bf Notations :}We shall adopt standard notations from the Geometric Measure Theory such as $M(A)$ for the Mass of a current $A$,
${\mathcal F}(A)$ for it's flat norm, $A\res E$ for it's restriction to a measurable subset $E$...etc, we refer
the reader to \cite{Fe}.

\medskip
 
{\bf Preliminaries : }Since our result is a local one we shall work in a neighborhood $U$ of a point $x_0$, use a symplectic form
$\om$ compatible with $J$. We denote by $g$ the metric generated by $J$ and $\om$ :
$g(\cdot,\cdot)=\om(\cdot,J\cdot)$. We also introduce normal coordinates $(x_1,x_2,\cdots,x_{2p-1},x_{2p})$
about $x_0$ in $U$ which can
 be chosed such that
\be
\label{II.0}
\mbox{ at }x_0\quad\quad J_{x_0}\cdot\frac{\p}{\p x_{2i+1}}=\frac{\p}{\p x_{2i+2}}\quad\quad\mbox{ for }
i=0\cdots p-1\quad .
\ee
Since $C$ is a calibrated current in $(U,\om,J)$, it is an area minimizing current and it's generalized
 mean curvature vanishes (see \cite{All} or \cite{Si}). One may isometrically embedd $(U,g)$ into an euclidian space ${\R}^{2p+k}$ and the generalised mean curvature of $C$ in ${\R}^{2p+k}$ coincides with the 
mean curvature of the embedding of $(U,g)$ and is therefore a bounded function. Combining this fact 
together with the monotonicity formula (17.3) of \cite{Si} we get that
\be
\label{II.1}
\frac{M(C\res B_r(x_0))}{r^2}=f(r)+O(r)\quad ,
\ee
where $f(r)$ is an increasing function, $M$ denotes the Mass of a current and $C\res B_r(x_0)$ is the 
restriction of $C$ to the geodesic ball of center $x_0$ and radius $r$. There exists in fact a constant $\al$
depending only on $g$ such that $e^{\al r}\frac{M(C\res B_r(x_0))}{r^2}$ is an increasing function
in $r$ (see \cite{Si}). The factor $e^{\al r}$ is a perturbation of an order which will have no influence
on the analysis below, therefore, by an abuse of notation we will often omitt to write it and consider straight
that $\frac{M(C\res B_r(x_0))}{r^2}$ is an increasing function.

By the mean of the coordinates $(x_1\cdots x_{2p})$ we shall identify $U$ with a subdomain in ${\R}^{2p}$
and use the same notation $C$ for the push forward of $C$ in ${\R}^{2p}$ by this chart. For small radii $r$
we introduce the dilation function $\la^{r,x_0}(x)=\frac{x-x_0}{r}$, and we introduce
the following dilation of $C$ about $x_0$ with rate $r$ as being the following current in ${\R}^{2p}$
\be
\label{II.2}
C_{r,x_0}:=\lf({\la^{r,x_0}}_\ast C\rg)\res B^{2p}_1(0)\quad . 
\ee
Observe that $r^{2}M_0(C_{r,x_0})=M_0(C\res B_r^{2p}(0))$ where $M_0$ denotes the mass in for
the flat metric $g_0$ in ${\R}^{2p}$. Since $g=g_0+O(r^2)$, we deduce from (\ref{II.1}) that 
$M_0(C_{r,x_0})$ is uniformly bounded as $r$ tends to zero. Again since $g$ and $g_0$ coincide
up to the second order, it does not hurt in the analysis below if one mixes the notations for
the two masses $M$ and $M_0$ and speak only about $M$. Since now $C$ is a cycle in $U$,
$\p C_{r,x_0}\res B_1^{2p}(0)=0$ and we can apply Federer-Fleming compactness theorem to deduce that,
from any sequence $r_i\rightarrow 0$ one can extract a subsequence $r_{i'}$ such that
$C_{r_{i'},x_0}$ converges in Flat norm to a limiting current $C_{0,x_0}$ called a tangent cone of $C$ at
$x_0$. One of the purpose of the next section will be to establish that $C_{0,x_0}$ is independent
of the subsequence and that the tangent cone is unique. The lower semi-continuity of the mass 
under weak convergence implies that
\be
\label{II.3}
\lim_{r\rightarrow 0}\frac{M(C\res B_r(x_0))}{r^2}=\lim_{r\rightarrow 0}M(C_{r,x_0})\ge M(C_{0,x_0})
\ee
In fact, from the fact that $C$ is calibrated by $\om$ we deduce now that the inequality (\ref{II.3}) 
is an equality. Indeed 
\[
M(C_{r,x_0})=r^{-2}C\res B_r^{2p}(0)(\om)=C_{r,x_0}\lf(r^2(\la^{r,x_0})^\ast\om\rg)
\]
It is clear that $\lim_{r\rightarrow 0}\|r^2(\la^{r,x_0})^\ast\om-\om_0\|_\infty$
where $\om_0=\sum_{i=1}^p dx_{2i-1}\wedge dx_{2i}$
therefore 
$C_{r,x_0}(r^2(\la^{r,x_0})^\ast\om-\om_0)\longrightarrow 0$ and we get that
\be
\label{II.3a}
\lim_{r\rightarrow 0}M(C_{r,x_0})=\lim_{i'\rightarrow +\infty}C_{r_{i'},x_0}(\om_0)=C_{0,x_0}(\om_0)
\ee
Since the comass of $\om_0$ is equal to 1, $C_{0,x_0}(\om_0)\le M(C_{0,x_0})$. Combining this last fact with
(\ref{II.3}) and (\ref{II.3a}) we have established that
\be
\label{II.3b}
\lim_{r\rightarrow 0}M(C_{r,x_0})=M(C_{0,x_0})=C_{0,x_0}(\om_0)
\ee
which means in particular that $C_{0,x_0}$ is calibrated by the K\"ahler form $\om_0$ in
 $(R^{2p},J_0)\simeq {\C}^p$ which is equivalent to the fact that $C_{0,x_0}$ is $J_0-$holomorphic.
Using the explicit form of the monotonicity formula (see \cite{Si} page 202), one observes that for any
$s\in{\R}_+^\ast$  
\[
C_{0,x_0}=\la^s_\ast C_{0,x_0}
\]
which means that ${\mathcal H}^2$ almost everywhere on the carrier ${\mathcal C}_{0,x_0}$ of $C_{0,x_0}$,
 $\frac{\p}{\p r}$ is in the approximate tangent plane to $C_{0,x_0}$, in other words, $C_{0,x_0}$
is a cone. Since it is $J_0-$holomorphic, ${\mathcal H}^2-$a.e. $x$ in ${\mathcal C}_{0,x_0}$,
 the approximate tangent cone is given by 
\[
T_xC_{0,x_0}=\mbox{Span}\lf\{\frac{\p}{\p r},J_0\frac{\p}{\p r}\rg\}
\]
Integral curves of $J_0\frac{\p}{\p r}$ are great-circles, fibers of the Hopf fibration
\[
(z_1=x_1+ix_2,\cdots,z_p)\longrightarrow [z_1,\cdots,z_p]
\]
therefore we deduce that $C_{0,x_0}$ is the sum of the integrals over radial extensions
of such great circles $\Gamma_1\cdots \Gamma_Q$ in $S^{2p-1}$ which is the integral over a sum of
$Q$ flat holomorphic disks. We adopt the following notation for the radial extensions
in $B^{2p}_1(0)$ of curents supported in $\p B^{2p}_1(0)$
\[
C_{0,x_0}=\oplus_{i=1}^Q0\sharp \Gamma_i\quad .
\]
Then we deduce that
\be
\label{II.3c}
\lim_{r\rightarrow 0} M(C_{r,x_0})=\pi Q\in \pi{\Z}\quad .
\ee
For any $x\in U$ one denotes $Q_x$ the integer such that
\[
\lim_{r\rightarrow 0}\frac{M(C\res B_r(x))}{\pi r^2}=Q_x
\]
Using the monotonicity formula, it is straightforward to deduce that for any $Q\in {\N}$
\[
{\mathcal C}_Q=\lf\{ x\in U\quad\mbox{ s. t. }\quad 0<Q_x\le Q\rg\}
\]
is an open subset of ${\mathcal C}_\ast=\lf\{ x\in U\quad\mbox{ s. t. }\quad 0<Q_x\rg\}$.
For $Q>1$, let us also denote 
\[
Sing^Q=\lf\{x\in {\mathcal C}_\ast\quad\mbox{ s. t. }\quad Q_x=Q\quad\mbox{ and }x
\mbox{ is an acc. point of }{\mathcal C}_{Q-1}\rg\}
\]
 Observe that, from Allard's theorem, it is clear that $C\res (U\setminus \cup_{Q}Sing^Q)$ is the integral
along a smooth surface with a smooth integer multiplicity. Although we won't make use of Allard's theorem this
justifies a-priori our notation. The whole purpose of our paper is to show
that $\cup_Q Sing^Q$ is made of isolated points. As we said, we won't make use of Allard's paper below
since the relative Lipschitz estimate we establish in lemma~\ref{lm-IV.2} gives Allard's result in our case
which is more specific.
Because of this nice stratification of $C$ ( ${\mathcal C}_Q$ is open in ${\mathcal C}_\ast$) we can argue
 by induction on $Q$ . Let ${\mathcal P}_Q$ be the following assertion
\be
\label{zz}
{\mathcal P}_Q\quad :\quad \cup_{q\le Q}Sing^q\mbox{ is made of isolated points}\quad .
\ee
From the begining of chapter IV until chapter VII we will assume either $Q=2$ or that ${\mathcal P}_{Q-1}$
holds and the goal will be to  establish ${\mathcal P}_Q$.

\section{The uniqueness of the tangent cone.}
\reset

The uniqueness of the tangent cone means that the limiting cone $C_{0,x_0}$, obtained in the 
previous section while
dilating at a point following a subsequence of radii $r_{i'}$, is independent of the subsequence
and is unique. Since our calibrated two dimensional rectifiable cycle is area minimizing,
this fact is a consequence of B. White upper-epiperimetric inequality in \cite{Wh}
(see also \cite{Ri2} for the justification of the prefix ``upper''). 
We need, however, a more quantitative version of this uniqueness of the tangent cone and
express how far we are from the unique tangent cone in terms of the closedness of the density of area
$M(C\res B_r(x_0))/\pi r^2$ to the limiting density $Q$. Precisely the goal of this section is to prove
the following lemma
\begin{Lm}
\label{lm-III.1}
{\bf (Uniqueness of the tangent cone.)}
For any $\ep>0$ and $Q\in{\N}$ there exists $\delta>0$ and $\rho_\ep\le 1$ such that, for any compatible
pair $(J,\om)$ almost complex structure-symplectic form over $B^{2p}_1(0)$ satisfying $J(0)=J_0(0)$,
$\om(0)=\om_0(0)$
\be
\label{III.1za}
\|J-J_0\|_{C^2(B_1)}+\|\om-\om_0\|_{C^2(B_1)}\le\delta\quad ,
\ee
for any $J-$holomorphic integral 2-cycle $C$ in $B_1(0)$ such that $Q_0=Q$, if
\[
M(C\res B_1^{2p}(0))\le \pi Q+\delta
\]
then, there exists $Q$ $J_0$ holomorphic flat discs $D_1\cdots D_Q$ passing through 0, intersection of holomorphic lines
of ${\C}^p$ with $B_1^{2p}(0)$, such that, for any $\rho\le\rho_\ep$
\be
\label{III.1aa}
{\mathcal F}(C_{\rho,0}-\oplus_{i=1}^QD_i)\le\ep
\ee
and for any $\psi\in C^\infty_0(B^1\setminus\{x\in B^1\ ;\ \mbox{dist}(x,\cup_i D_i)\le \ep |x|\})$,
\be
\label{III.1b}
C_{\rho,0}(\psi)=0
\ee
\end{Lm}
Before proving lemma~\ref{lm-III.1}, we first establish the following intermediate result
\begin{Lm}
\label{lm-III.2}
For any $\ep>0$ and $Q\in{\N}$ there exists $\delta>0$ such that, for any compatible
pair $(J,\om)$ almost complex structure-symplectic form over $B^{2p}_1(0)$ satisfying $J(0)=J_0(0)$,
$\om(0)=\om_0(0)$
\[
\|J-J_0\|_{C^2(B_1)}+\|\om-\om_0\|_{C^2(B_1)}\le\delta\quad ,
\]
for any $J-$holomorphic integer rectifiable 2-cycle $C$ such that $Q_x=Q$, if
\[
M(C\res B_1^{2p}(0))\le \pi Q+\delta
\]
then, there exist $Q$ $J_0$ holomorphic flat discs $D_1\cdots D_Q$ passing through 0, intersection of holomorphic lines
of ${\C}^p$ with $B_1^{2p}(0)$, such that, 
\[
{\mathcal F}(C\res B_1(0)-\oplus_{i=1}^QD_i)\le\ep
\]  
\end{Lm}
\begin{Rm}
Lemma~\ref{lm-III.2} give much less information than lemma~\ref{lm-III.1}. Since a-priori in lemma~\ref{lm-III.2}
the disks $D_i$ may vary a lot as one dilate $C$ about $0$, whereas lemma~\ref{lm-III.1} prevents such a tilting
as one dilates the current further.
\end{Rm}
{\bf Proof of lemma~\ref{lm-III.2} :}
We prove lemma~\ref{lm-III.2} by contradiction. Assume there exists $\ep_0>0$ , $\delta_n\rightarrow 0$,
compatible $J_n$ and $\om_n$ and $C_n$ such that
\begin{itemize}
\item[i)]
\[
\|J_n-J_0\|_{C^2}+\|\om_n-\om_0\|_{C^2}\le\delta_n
\]
\item[ii)]
\[
\lim_{r\rightarrow 0}\pi^{-1}r^{-2}M(C_n\res B_r(0))= Q
\]
\item[iii)] 
\[
M(C_n\res B_1)\le\pi Q+\delta_n
\]
\item[iv)]
\be
\label{III.1a}
\inf\lf\{{\mathcal F}(C_n\res B_1-\oplus_{i=1}^Q D_i)\ \mbox{ s.t. }D_i\mbox{ flat holom discs}, 0\in D_i\rg\}\ge \ep_0
\ee 
\end{itemize}
Since $\p C_n\res B_1=0$ and since the mass of  $C_n$ is uniformly bounded, one may assume, modulo extraction
of a subsequence if necessarily, that $C_n$ converges to a limiting rectifiable cycle $C_\infty$.
Exactly like in section III we have the fact that for any $0<r\le 1$
\be
\label{III.2}
\lim_{n\rightarrow +\infty}M(C_n\res B_r)=M(C_\infty\res B_r)= C_\infty\res B_r(\om_0)
\ee
We deduce then that  $C_\infty$ is calibrated by $\om_0$ and is therefore a $J_0-$holomorphic cycle.
Using ii) we deduce also that
\[
\lim_{r\rightarrow 0}\pi^{-1}r^{-2}M(C_\infty\res B_r(0))= Q
\]
and finally, from iii) and the lower semicontinuity of the mass, we have that $M(C_\infty\res B_1)=\pi Q$
Thus, since $\pi^{-1}r^{-2}M(C_\infty\res B_r(0))$ is an increasing function, we have established that on $[0,1]$
\be
\label{III.3}
\pi^{-1}r^{-2}M(C_\infty\res B_r(0))\equiv Q
\ee
Let, for almost every $r$,  $S^r_\infty=<C_\infty,dist(\cdot,0),r>$ be the slice current obtained by slicing
 $C_\infty$ with $\p B_r(0)$ (see \cite{Fe} 4.2.1).
By Fubini, we have that for a.e. $0<r<1$
\[
M(<C_\infty, dist(\cdot,0),r>)\le 2\pi Q r
\]
Let $0\sharp S^r_\infty$ be the radial extension of $S^r_\infty$ in $B_r(0)$.
\[
M(0\sharp S^r_\infty)=\frac{r}{2}M(S^r_\infty)=\pi Q r^2=M(C_\infty\res B_r(0))
\]
Since $\p(C_\infty\res B_r(0)-0\sharp S^r_\infty)=0$, and since $C_\infty\res B_r(0)$ is area minimizing
we have that $0\sharp S^r_\infty$ is also area minimizing. Let $\al$ such that $d\al=\om$.
\[
M(0\sharp S^r_\infty)=M(C_\infty\res B_r(0))=C_\infty\res B_r(0)(\om_0)=S^r_\infty(\al)=(0\sharp S^r_\infty)(\om_0)
\]
Therefore $0\sharp S^r_\infty$ is an holomorphic cone which is a cycle. So we deduce like in section II
that $0\sharp S^r_\infty$ is a sum of flat holomorphic disk for any $r$. Thus $C_\infty$ is also a sum
\[
C_\infty\res B_1(0)=\sum_{i=1}^Q D_i
\]
where each $D_i$ is the intersection of a complex straight line in ${\C}^p$ with $B^{2p}_1$.
From Federer-Fleming compactness theorem we have the fact that the weak convergence of $C_n$ to $C_\infty$
holds in flat norm
\[
{\mathcal F}(C_n\res B_1-\sum_{i=1}^Q D_i)\longrightarrow 0
\]
which contradicts iv) and lemma~\ref{lm-III.2} is proved. \cqfd

\medskip

{\bf Proof of lemma~\ref{lm-III.1} :}

We first prove assertion (\ref{III.1aa}). We first recall Brian White's upper-epiperimetric inequality
adapted to our present context : Brian White's upper-epiperimetric inequality was proved
for area minimizing surfaces in ${\R}^{2p}$. Here, in the present situation, we are dealing with
area minimizing currents which are $J-$holomorphic for a metric $g=\om(\cdot,J\cdot)$ which get's
as close as we want to the standard one because of assumption (\ref{III.1za}). Therefore very minor
changes have to be provided to adapt B.White theorem to the present context. An adaptation of
the epiperimetric inequality for ambiant non flat metric is also given in \cite{Ch} Appendix A.
So we have the following result.

 Given an integer $Q$, there exists a positive number $\ep_Q>0$,
such that, for any compatible pair $\om$, $J$ in $B^{2p}_2(0)$ satisfying
$\|\om-\om_0\|_{C^2(B_2)}+\|J-J_0\|_{C^2(B_2)}\le \ep_Q$ and for any
 $C$ $J-$holomorphic 2-rectifiable integral current in $B_2^{2p}(0)$, satisfying
 $\p C\res B^{2p}_2=0$,
 assuming there exist $Q$ flat holomorphic disks $D_1\dots D_Q$ in 
$(B^{2p}_1(0),J)\simeq {\C}^p\cap B^{2p}_1(0)$ passing through the origin such that
\be
\label{I.0}
{\mathcal F}\lf(C_{2,0}\res B_1(0)-\sum_{i=1}^QD_i\rg)\le\ep_Q
\ee
(where we used a common notation for the oriented 2-disks $D_i$ and the corresponding 2-currents)
 then
\be
\label{I.1}
M(C\res B_1^{2p})-\pi Q\le (1-\ep_Q)\ \lf(\frac{1}{2}M\lf(\p(C\res B_1^{2p})\rg)-\pi Q\rg)
\ee
\begin{Rm}
\label{rm-III.1}
Observe that in the statement of the epiperimetric property in Definition 2 of \cite{Wh} the epiperimetric
 constant $\ep_Q$ may a-priori also depend of the cone $\sum_{i=1}^Q D_i$. It is however elementary to 
observe that this
space of cones made of the intersection of $Q$ holomorphic straight lines passing through the origin
with $B^{2p}_1(0)$ is compact for the flat distance. Now using a simple finite covering argument
for this space of cones by balls (for the flat distance) permits to obtain a constant $\ep_Q>0$
for which the epiperimetric property holds independently of the cone $\sum_{i=1}^Q D_i$.  
\end{Rm}

Once again we shall ignore the factor $e^{\al r}$ in front of $r^{-2}M(C\res B_r)$ which induces lower
order perturbations and argue as $r^{-2}M(C\res B_r)$ itself would be an increasing function
(observe also that $\al$ is as small as we want as $J$ and $\om$ are chosen as close
as we want from $J_0$ and $\om_0$).

Let then $\ep>0$ such that $\ep<\ep_Q$ and let $\delta>0$ given by lemma~\ref{lm-III.2} for that $\ep$.
 Assuming then $M(C\res B_1)\le\pi Q+\delta$
implies from the monotonicity formula that for any $r<1$ $r^{-2}M(C\res B_r(0))=M(C_{r,0})\le(\pi Q+\delta)$.
 Applying then Lemma~\ref{lm-III.2} to $C_{2r,0}$ for $r<1/2$ we deduce the existence of $Q$ flat disks
 $D_1\cdots D_Q$ such that
\be
\label{III.5a} 
{\mathcal F}(C_{2r,0}-\sum_{i=1}^Q D_i)\le \ep
\ee
We can then apply the epiperimetric inequality to $C_{r,0}$ and we get, after rescaling, that 
\be
\label{III.5b}
M(C\res B_r(0))-\pi Q r^2\le (1-\ep_Q) \lf(\frac{r}{2} M(\p(C\res B_r(0)))-\pi Q r^2\rg)
\ee 
Denote $f(r)=M(C\res B_r(0))-\pi Q r^2$. $f'(r)\ge M(\p(C\res B_r(0)))-2\pi Q r$. Therefore (\ref{III.5b})
implies
\[
\frac{1-\ep_Q}{2}r\ f'(r)\ge f(r)\quad .
\]
Integrating this differential inequality between  $s$ and $\sigma$ ($1/2>s>\sigma$),
 $f(s)\ge\lf(\frac{s}{\sigma}\rg)^\frac{2}{1-\ep_Q}
 f(\sigma)$. Let $\nu=\frac{2}{1-\ep}-2>0$, we then have
\be
\label{III.6}
\frac{f(s)}{s^2}\ge\lf(\frac{s}{\sigma}\rg)^\nu\ \frac{f(\sigma)}{\sigma^2}\quad .
\ee
Let $F(x)=\frac{x}{|x|}$. We have
\[
M(F_\ast\lf(C\res B_s(0)\setminus B_\sigma(0)\rg))=\int_{B_S(0)\setminus B_\sigma(0)}
\frac{1}{|x|^3}\lf|\tau\wedge\frac{x}{|x|}\rg|\ |\theta| \ d{\mathcal H}^{2}\res {\mathcal C}
\]
where $\tau$ denotes the unit $2-$vector associated to the oriented approximate
tangent plane to $C$ and  defined ${\mathcal H}^2-$a.e. along the carrier ${\mathcal C}$ of 
the rectifiable current, $\theta$ is the $L^1({\mathcal C})$ integer-valued multiplicity of
$C$ (i.e using classical GMT notations :$C=<{\mathcal C},\theta,\tau>$) and 
$d{\mathcal H}^2\res{\mathcal C}$ is the restriction to ${\mathcal C}$ of the $2-$dimensional
Hausdorff measure. Using Cauchy-Schwarz and 5.4.3 (2) of \cite{Fe} (the explicit formulation of the monotonicity
formula) and (\ref{III.6}), we have
\be
\label{III.6b}
\begin{array}{l}
\ds M(F_\ast\lf(C\res B_s(0)\setminus B_\sigma(0)\rg))\\[5mm]
\ds\le\lf[\frac{M(C\res B_s(0))}{s^2}
-\frac{M(C\res B_\sigma(0))}{\sigma^2}\rg]^\frac{1}{2}\ \lf[\frac{M(C\res B_s(0))}{\sigma^2}\rg]^\frac{1}{2}\\[5mm]
\ds\quad\le\lf[\frac{M(C\res B_s(0))}{s^2}
- \pi Q\rg]^\frac{1}{2}\ \lf[\frac{M(C\res B_s(0))}{\sigma^2}\rg]^\frac{1}{2}\\[5mm]
\ds\quad\le\lf[\frac{f(s)}{s^2}\rg]^\frac{1}{2}\
 \lf[\frac{s^2}{\sigma^2}\frac{M(C\res B_s(0))}{s^2}\rg]^\frac{1}{2}\le K\ s^\frac{\nu}{2}\ \frac{s}{\sigma} 
\end{array}
\ee
Let $r<\rho<1/2$, applying (\ref{III.6b}) for $s=2^{-k}\rho$ and $\sigma=2^{-k-1}\rho$
for $k\le \log_2\frac{\rho}{r}$ and summing over $k$ we get
\be
\label{III.6c}
 M(F_\ast\lf(C\res B_\rho(0)\setminus B_r(0)\rg))\le C\rho^\frac{\nu}{2}
\ee
Observe that $\p(F_\ast\lf(C\res B_\rho(0)\setminus B_r(0)\rg))=\p C_{\rho,0}-\p C_{r,0}$. Therefore
we deduce
\be
\label{III.6ba}
{\mathcal F}\lf((C_{\rho,0}-C_{r,0})\res B_1(0)\setminus B_{\frac{1}{2}}(0)\rg)\le C\rho^\frac{\nu}{2}
\ee
Since
\[
{\mathcal F}\lf((C_{\rho,0}-C_{r,0})\res B_\frac{1}{2}(0)\setminus B_\frac{1}{4}(0)\rg)
\le\lf(\frac{1}{3}\rg)^\frac{1}{3}{\mathcal F}\lf((C_{\frac{\rho}{2},0}-C_{\frac{r}{2},0})\res B_1(0)\setminus B_\frac{1}{2}(0)\rg)
\]
applying then (\ref{III.6ba}) for $\rho$, $r$ replaced by $2^{-k}\rho$, $2^{-k} r$ and summing over
$k=1,\cdots\infty$ we finally obtain
\be
\label{III.6bb}
{\mathcal F}\lf((C_{\rho,0}-C_{r,0})\res B_1(0)\rg)\le C\rho^\frac{\nu}{2}
\ee
which is the desired inequality (\ref{III.1aa}).

 It remains to show (\ref{III.1b}) in order to finish
the proof of lemma~\ref{lm-III.1}. We argue by contradiction. Assume there exists $\ep_0>0$ , 
$\rho_n\rightarrow 0$ and $\psi_n\in C^\infty_0(\wedge^2 B_1)$ such that 
\[
\mbox{supp}\, \psi_n\subset E_0=\{x\in B^1\ ;\ \mbox{dist}(x,\cup_i D_i)\le \ep_0 |x|\}\quad,
\]
where $C_{0,0}=\oplus_{i=1}^Q D_i$, and 
\[
C_{\rho_n,0}(\psi_n)\ne 0
\]
This later fact implies in particular that there exists $x^n\in E_0$ such that 
$\lim_{r\rightarrow 0}M(C_{r,x_n})\ne 0$. Using the monotonicity formula we deduce then that
\[
M(C_{\rho_n|x_n|,0}\res B_{\ep_0/2}(\frac{x_n}{|x_n|}))\ge \frac{\pi}{4}\ep_0^2
\]
We may then extract a subsequence such that $\frac{x_n}{|x_n|}\rightarrow x_\infty$. Thus we have 
\[
M(C_{\rho_n|x_n|,0}\res B_{3\ep_0/4}(x_\infty))\ge \frac{\pi}{4}\ep_0^2
\]
We have
\[
M(C_{\rho_n|x_n|,0}\res B_{3\ep_0/4}(x_\infty))=C_{\rho_n |x_n|,0}\res B_{3\ep_0/4}(x_\infty)
\lf({\frac{x}{\rho_n |x_n|}}^\ast\om_n\rg)
\]
Since $\|\om_n-\om_0\|_{C^2}\rightarrow 0$ and since $\om_n(0)=\om_0(0)$ we clearly have that
\[
\|{\frac{x}{\rho_n |x_n|}}^\ast\om_n-\om_0\|_\infty\longrightarrow 0
\]
Therefore
\[
\begin{array}{l}
\ds\lf|C_{\rho_n |x_n|,0}\res B_{3\ep_0/4}(x_\infty)
\lf({\frac{x}{\rho_n |x_n|}}^\ast\om_n-\om_0\rg)\rg|\\[5mm]
\le\quad\le M(C_{\rho_n|x_n|,0}\res B_{3\ep_0/4}(x_\infty))
\ \|{\frac{x}{\rho_n |x_n|}}^\ast\om_n-\om_0\|_\infty\longrightarrow 0\quad ,
\end{array}
\]
Thus
\[
C_{0,0}\res B_{3\ep_0/4}(x_\infty)(\om_0)=\lim_{n\rightarrow +\infty}C_{\rho_n|x_n|,0}\res B_{3\ep_0/4}(x_\infty)(\om_0)\ge \frac{\pi}{4}\ep_0^2
\]
which contradicts the fact that $B_{3\ep_0/4}(x_\infty)\subset E_0$. Therefore (\ref{III.1b}) holds and
lemma~\ref{lm-III.1} is proved.\cqfd

\section{Consequences of lemma~\ref{lm-III.1} : No accumulation of points in $Sing^Q$ in the easy case -
The relative Lipschitz estimate in the difficult case.}
\reset

In this section we expose two important consequences of lemma~\ref{lm-III.1}. Before to expose them we 
first observe that in the way to show ${\mathcal P}_{Q-1}\Longrightarrow {\mathcal P}_Q$, two cases
have to be considered separately. The first case (the easy one) is the case where the tangent cone
at the point $x_0$ of multiplicity $Q$ (i.e. $\pi^{-1}r^{-2}M(B_r(x_0))\rightarrow Q$) is not
made of $Q$ times the same disks. The second case is the case where the tangent case is made of $Q$
times the same disk. In the first case we will deduce almost straight from lemma~\ref{lm-III.1} that such an $x_0$
cannot be an accumulation point of points of multiplicity $Q$ also, see lemma~\ref{lm-IV.1a} below.
 In the second case, much more analysis will be needed to reach the same statement and this is the purpose
of chapters IV until IX. We can nevertheless deduce in this chapter an important consequence of our
quantitative version of the
uniqueness of the tangent cone (lemma~\ref{lm-III.1}) for the difficult case : this is the so called
``relative lipshitz estimate'' (see lemma~\ref{lm-IV.2} below). This property says that, given a point
$x_0$ of multiplicity $Q$  whose tangent cone is $Q$ times a flat disk and given an $\ep>0$, there exists a radius 
$r_{\ep,x_0}>0$ such that given any two points of ${\mathcal C}_\ast\cap B_{r_{\ep,x_0}}(x_0)$, one of the
two being also of multiplicity $Q$, the slope they realise relative to the tangent cone of $x_0$ is
less than $\ep$. The condition that one of the two points has multiplicity $Q$ (this could be $x_0$ itself
for instance) is a crucial assumption. It is indeed straightforward to find counterexamples to 
any Lipschitz estimates of multivalued graphs of holomorphic curves : take for instance $w^2=z$
in ${\C}^2\simeq\{(z,w)\ z,w\in {\C}\}$ viewed as a 2-valued graph over the line $\{w=0\}$, all points
have multiplicity 1, $(0,0)$ included of course, but the best possible estimate is a H\"older one 
$C^{0,\frac{1}{2}}$. We cannot exclude that such a configuration exists as we dilate at a point $x_0$
 of mupltiplicity $Q>1$.

We first then prove the following consequence of lemma~\ref{lm-III.1}

\begin{Lm}
\label{lm-IV.1a}
{\bf(no accumulation - the easy case)}
Let $Q\in {\N}$, $Q\ge 2$. Let $x_0$ be a point in ${\mathcal C}_{Q}\setminus {\mathcal C}_{Q-1}$
(i.e.$\pi^{-1}r^{-2}M(B_r(x_0))\rightarrow Q$ as $r\rightarrow 0$). Assume that the tangent cone
at $x_0$, $C_{0,x_0}$, contains at least two different flat $J_{x_0}-$holomorphic disks (i.e. $C_{0,x_0}\ne
Q\, D$ for any flat $J_{x_0}-$holomorphic disk $D$) then, there exists $r>0$ sucht that  
\[
B_r(x_0)\cap ({\mathcal C}_{Q}\setminus {\mathcal C}_{Q-1})=\{x_0\}\quad .
\]
\end{Lm}
{\bf Proof of lemma~\ref{lm-IV.1a}.}

Let $x_0$ being as in the statement of the lemma : $x_0\in {\mathcal C}_{Q}\setminus {\mathcal C}_{Q-1}$
and $C_{0,x_0}=\oplus_{i=1}^K Q_i\ D_i$ (where $D_i\ne D_j$ for $i\ne j$ and $K>1$). Let $\ep>0$ be a positive number  smaller than
 1/4Q times the
maximal angle $\al$ between the various diks $D_1,\cdots ,D_K$ in the tangent cones in such a way that there exists
$i\ne j$ such that
\[
E_\ep(D_i)\cap E_\ep(D_j)=\{x_0\}\quad ,
\]
where we are using the following notation 
\[
E_\ep(D_i)\{x\in {\R}^{2p}\ ;\ \mbox{dist}(x, D_i)\le \ep |x-x_0|\})\quad .
\]
By taking $\ep$ as small as above, we even have ensured that $\cup E_\ep(D_i)\setminus x_0$ has at least two 
connected components
whose intersections with $\p B_1(x_0)$ are at a distance larger than $\al/2$.
We prove now lemma~\ref{lm-IV.1a} by contradiction : we assume there exists 
$x_n\in{\mathcal C}_{Q}\setminus {\mathcal C}_{Q-1}$ such that $x_n\rightarrow x_0$
and $x_n\ne x_0$.
Let $\delta>0$ given by lemma~\ref{lm-III.1} for $\ep$ chosen just above. Let $\rho>0$
such that $\rho^{-2}M(C\res B_\rho(x_0))\le \pi Q+\delta/2$. For any
 $x\in B_{\rho\frac{\delta}{4\pi Q}}(x_0)$ we have
\be
\label{IV.aa}
\begin{array}{l}
M(C\res B_{\rho(1-\frac{\delta}{4\pi Q})}(x))\le M(C\res B_\rho(x_0))\le \rho^2\ 
\lf(\pi Q+\frac{\delta}{2}\rg)\\[5mm]
\quad\le (\rho(1-\frac{\delta}{4\pi Q}))^2\  (1-\frac{\delta}{4\pi Q})^{-2} \lf(\pi Q+\frac{\delta}{2}\rg)\\[5mm]
\quad\le (\rho(1-\frac{\delta}{4\pi Q}))^2\ \lf(\pi Q+\delta\rg)
\end{array}
\ee
Choose then $x_n\in  B_{\rho\frac{\delta}{8\pi Q}}(x_0)$. Applying (\ref{III.1b}) for $x_0$ we know
that $x_n$ is contained in one of the $E_\ep(D_i)$, say $E_\ep(D_1)$. Denote $E_1$ the connected component
of $\cup E_\ep(D_i)\setminus \{x_0\}$ that contains $E_\ep(D_1)$. $\ep$ has been chosen small enough such
 that $\cup E_\ep(D_i)$
has at least two connected components. Therefore we can chose $D_j$ such that $E_\ep(D_j)$ is disjoint
from the component containing $E_\ep(D_1)\setminus \{x_0\}$. Let $\al$ be the angular distance, relative
to $x_0$, from $E_\ep(D_j)$
and the component containing $E_\ep(D_1)\setminus \{x_0\}$. $al$ is clearly bounded from below by a positive
number as one choses 
$\ep$ smaller and smaller. Applying lemma~\ref{lm-III.1} this time to $x_n$,
we know that in $B_{4|x_n|}(x_n)\setminus B_{|x_n|}(x_n)$ the support of $C$ is at the $\ep |x_n|$ distance from 
a union of flat disks passing through $x_n$ (the tangent cone at $x_n$). This imposes that the angular distance
 between the tangent cone at $x_n$ and $D_1$ is less than $C_Q\ep$, where $C_Q$ depends on
$Q$ only . Therefore
\be
\label{IV.aaa}
\mbox{supp}(C\res B_{2|x_n|}(x_n)\setminus B_{|x_n|}(x_n))\subset \ti{E}_1=\{x\ ;\ dist(x, D_1)\le C_Q\ep|x_n|\}  
\ee
Observe that $dist\{E_\ep(D_j)\cap \lf (B_{|x_n|}(x_0)\setminus B_{|x_n|/2}(x_0)\rg); \ti{E}_1\}\ge \al/4$.
This later fact combined with (\ref{IV.aaa}) contradicts (\ref{III.1aa}). Lemma~\ref{lm-IV.1a} is then proved.\cqfd

From now on until the begining of chapter X we will be dealing with the difficult case only : 
the cas where the point $x_0$ of multiplicity $Q$ has a tangent cone which is made of $Q$
time the same disk. As we use to do since chapter II, we will work in a neighborhood of 
$x_0$ where a compatible simplectic form $\om$ for $J$ exists, and we shall use normal coordinates
for $g(\cdot,cdot)=\om(\cdot,J\cdot)$ about $x_0$, compatible with $J_{x_0}$ at $x_0$, satisfying (\ref{II.0}),
and we can also assume that the tangent cone at $x_0$ is   
\be
\label{IV.aab}
C_{0,x_0}\res B_1(0)=Q [D_0]\quad 
\ee
where $D_0$ is the flat oriented disk whose tangent $2-$vector is $\frac{\p}{\p x_1}\wedge \frac{\p}{\p x_2}$.
From now on we will also use the following notations for complex coordinates about $x_0$
\be
\label{IV.aac}
z=x_1+ix_2\quad\mbox{ and }\quad w_i=x_{2k+1}+i x_{2k+2}\mbox{ for }k=1\cdots p-1\quad .
\ee
We will also denote $w=(w_1,\cdots, w_{p-1})$. 
A second consequence to lemma~\ref{lm-III.1} is the following result :
\begin{Lm}
\label{lm-IV.2}
{\bf (the relative Lipschitz estimate)} :
Let $x_0$ be a point of multiplicity $Q$ (i.e. $x_0\in {\mathcal C}_Q\setminus {\mathcal C}_{Q-1}$),
assume the tangent cone $C_{0,x_0}\res B_1(0)$ at $x_0$ is $Q$ times a flat disk (i.e. of the form (\ref{IV.aab})).
Let $\ep>0$, then there exists $r_{\ep,x_0}$ such that for any $r\le r_{\ep,x_0}$
\be
\label{IV.1a}
\forall x\in B_r(x_0)\cap ({\mathcal C}_Q\setminus{\mathcal C}_{Q-1})\quad\quad
{\mathcal F}((C_{0,x}-C_{0,x_0})\res B_1(0))\le\ep
\ee
and $\forall x=(z,w)\in B_r(x_0)\cap ({\mathcal C}_Q\setminus {\mathcal C}_{Q-1})$ and $x'=(z',w')\in
B_r(x_0)\cap {\mathcal C}_\ast$ we have
\be
\label{IV.1b}
|w-w'|\le \ep |z-z'|
\ee
\end{Lm}
{\bf Proof of lemma~\ref{lm-IV.2}.} 
Let $\ep>0$ and $\delta>0$ given by lemma~\ref{lm-III.1}. Choose $r_1$ such that
\[
M(C\res B_{r_1}(x_0))\le r^2\ \lf(\pi Q+\frac{\delta}{2}\rg)
\]
This implies in particular that for any $r<r_1$
\be
\label{IV.01}
{\mathcal F}((C_{r,x_0}-C_{0,x_0})\res B_1(0))\le\ep
\ee

As in the proof of lemma~\ref{lm-IV.1a}, see (\ref{IV.aa}) we have that for any
 $x\in B_{r_1\frac{\delta}{4\pi Q}}(x_0)$ and $r<r_1(1-\frac{\delta}{4\pi Q})$
\[
M(C\res B_r(x))\le r^2(\pi Q+\delta)
\]
Let then $x\in B_{r_1\frac{\delta}{4\pi Q}}(x_0)\cap ({\mathcal C}_Q\setminus {\mathcal C}_{Q-1})$,
applying lemma~\ref{lm-III.1},
we have then for $r<r_1(1-\frac{\delta}{4\pi Q})=r_2$
\be
\label{IV.1}
{\mathcal F}((C_{r,x}-C_{0,x})\res B_1(0))\le\ep
\ee
Choose now $x\in B_{\ep r_1}(x_0)\cap ({\mathcal C}_Q\setminus {\mathcal C}_{Q-1})$. $C_{r_2,x}$
is fast eine $\ep-$translation from $C_{r_1,x_0}$, therefore, since $M(C_{r,x_0})\le 2\pi Q$, we
have
\be
\label{IV.2}
{\mathcal F}((C_{r_2,x_0}-C_{r_2,x})\res B_1(0))\le 2\pi Q
\ee
Take $r_{\ep,x_0}=\min \{\ep r_1,\frac{\delta}{4\pi Q} r_1\}$. Combining (\ref{IV.01}), (\ref{IV.1})
 and (\ref{IV.2}), we deduce that 
\be
\label{IV.3}
\forall x\in B_{r_{\ep,x_0}}(x_0)
\cap ({\mathcal C}_Q\setminus{\mathcal C}_{Q-1})\quad\quad
{\mathcal F}((C_{0,x}-C_{0,x_0})\res B_1(0))\le(2+2\pi Q)\ep
\ee
It remains to check (\ref{IV.1b}) which is in fact an almost direct consequence of (\ref{III.1b})
and (\ref{IV.1a}). Lemma~\ref{lm-IV.2} is then proved.\cqfd

\section{The covering argument.}
\reset

Let $x_0$ a point of multiplicity $Q>1$ whose tangent cone $C_{0,x_0}$ is $Q$ times the integral
over the flat disk $D_0$ given by $w_i=0$ for $i\cdots Q-1$ (we use the system
of coordinates introduced in the begining of chapter IV) . The purpose of this section is to construct
a Whitney-Besicovitch covering $B^2_{r_i}(z_i)$ of $\Pi({\mathcal C}_{Q-1})\cap B^2_{\rho}(x_0)$ 
where $\Pi$ is the projection 
on $D_0$ which gives the first complex coordinate of each point ($\Pi(z,w_1\cdots w_{p-1})=z$), for
some radius $\rho$, small enough depending on $x_0$. This covering will be chosed in such a way
the following striking facts hold : first $C\res\Pi^{-1}(B^2_{r_i}(z_i))$ is in fact supported in a ball
of radius $2r_i$, $B^{2p}_{2r_i}(x_i)$,
moreover $C\res B^{2p}_{2r_i}(x_i)$ is ``splitted''. This last word means that the flat distance between $C\res\Pi^{-1}(B^2_{r_i}(z_i))$
and the $Q$ multiple of any single valued graph over $D_0$  is larger than $K\,r_i^3$ where $K$ only depends 
on $p$, $J$ and $\om$. This will come from the fact that $r_i$ may be chosed i such a way that
$r_i^{-2}M(B^{2p}_{r_i}(x_i)\le \pi- K'$ where again $K'>0$ only depends 
on $p$, $Q$, $J$ and $\om$. The existence of such a covering is a consequence of the ``splitting before tilting'' lemma
proved in \cite{Ri2}.

Let $\al$ given by lemma~\ref{lm-IV.3} and let $\ep>0$ to be chosen small enough, compare to $\al$ later.
Let $r_{\ep,x_0}$ be the radius given by lemma~\ref{lm-IV.2}. We may chose also $r_{\ep,x_0}$ small 
enough in such a way that
\be
\label{VI.a1}
\forall r\le r_{\ep,x_0}\quad\quad M(C_{r,x_0}\res B_1(0))\le \pi Q+\ep^2\quad .
\ee
Using the proof of lemma~\ref{lm-III.1} (from (\ref{III.6b}) until the end of the proof), we deduce that
\be
\label{VI.a2}
\forall r\le r_{\ep,x_0}\quad\quad{\mathcal F}((C_{r,x_0}-C_{0,x_0})\res B_1(0))\le K\ep\quad .
\ee
(In fact $\delta =O(\ep^2)$ works in the statement of lemma~\ref{lm-III.1}).
 In one hand, as in the proof of lemma~\ref{lm-IV.1a}, see (\ref{IV.aa}) we have that for any
 $x\in B_{\ep^2 r_{\ep,x_0}}(x_0)$ and $r\le r_{\ep,x_0}(1-\ep^2)$ 
\be
\label{VI.1a}
M(C\res B_r(x))\le r^2\ (\pi Q+\ep^2)\quad ,
\ee
and, on the other hand, arguing like in the proof
of lemma~\ref{lm-IV.2}, between (\ref{IV.1}) and (\ref{IV.3}), we have, using also (\ref{V.a11}) 
\be
\label{VI.1b}
\forall x\in B_{\ep^2}(x_0)\quad\quad{\mathcal F}(C_{\frac{r_{\ep,x_0}}{2},x_0}\res B_1(0) -\pi Q\ [D_0])
\le K\ep\quad. 
\ee
Having chosen $K\ep<\al$ we are in a position to apply the ``Splitting before tilting'' lemma of \cite{Ri2} 
which is a key step in our proof of the regularity of 1-1 rectifiable cycles.

\begin{Lm}
\label{lm-IV.3}
{\bf( splitting before tilting)} \cite{Ri2}
There exists  $\al>0$ such that for any $x_0\in {\mathcal C}_{Q-1}$ and for any radius
$0<\rho<\al$ satisfying
\be
\label{V.a1}
{\mathcal F}(C_{2\rho,x_0}\res B_1(0)-Q\ [D_0])\le\al
\ee
where $D_0$ is a flat $J_{x_0}-$holomorphic disk passing through $x_0$. Then, for any $r<\rho$ and any 
 $J_{x_0}-$holomorphic flat disk, $D_1$, passing through $x_0$ and satisfying
\be
\label{V.a2}
{\mathcal F}([D_0]-[D_1])\ge\frac{1}{4}
\ee
we have
\be
\label{V.a3z}
{\mathcal F}(C_{r,x_0}\res B_1(0)- Q[D_1])\ge \al\quad .
\ee
Moreover, there exists $r_0<\rho$ and $K_0$ a constant depending only on $\|\om\|_{C^1}$ and $\ep_Q^\pm$,
the epiperimetric constants, such that
\be
\label{V.b3}
M(C_{r_0,x_0}\res B_1(0))=\pi Q-K_0\al\quad ,
\ee
\be
\label{V.b4}
{\mathcal F}(C_{{r}_0,x_0}\res B_1(0)- Q [D_0])\le K\ \sqrt{\al}\quad
\ee
for some constant $K$ depending also only on $\|\om\|_{C^1}$ and $\ep_Q^\pm$. Finally,
 for any $J_{x_0}-$holomorphic disk $D$ passing through $0$
\be
\label{V.aa3}
\forall r\le r_0\quad\quad{\mathcal F}(C_{r,x_0}\res B_1(0)- Q\ [D])\ge \al\quad .
\ee 
\end{Lm}
For
any $x\in {\mathcal C}_{p-1}\cap B_{\ep^2 r_{\ep,x_0}}(x_0)$ we denote by $r_x$ the radius
$r_0$ given by the lemma. We then have
\be
\label{VI.1c}
M(C_{r_x,x}\res B_1(0))=\pi Q-K_0\al\quad ,
\ee
\be
\label{VI.1d}
{\mathcal F}(C_{{r}_x,x}\res B_1(0)-\pi Q [D_0])\le K\ \sqrt{\al}\quad
\ee
for some constant $K$ depending also only on $\|\om\|_{C^1}$ and $\ep_Q^\pm$.
Moreover, for $\al$ chosen smal enough and $\ep$ 
small enough compare to $\al$, the following lemma holds
\begin{Lm}
\label{lm-VI.1}
Under the above notations we have that for any $x\in {\mathcal C}_{p-1}\cap B_{\ep r_{\ep,x_0}}(x_0)$
\be
\label{V.a3}
 \mbox{supp}\lf(C\res \Pi^{-1}(B^2_{r_x}(\Pi(x))\cap B_{r_{\ep,x_0}}(x_0)\rg)
\subset B^2_{r_x}(\Pi(x))\times B^{2p-2}_{r_x}(0)\quad,
\ee
and that
\be
\label{V.a4}
\mbox{supp}\lf(C\res \Pi^{-1}(B^2_{r_x}(\Pi(x))\cap B_{r_{\ep,x_0}}(x_0)\rg)\subset{\mathcal C}_{p-1}
\ee
\end{Lm}
{\bf Proof of lemma~\ref{lm-VI.1}.}
We claim that for any $r$ between $\frac{r_{\ep,x_0}}{2}$ and $r_x$ one has
\be
\label{V.a5}
supp(C_{r,x}\res B_1(0))\subset E(\al^\frac{1}{16})
\ee
where we use the notation
\be
\label{V.ab5}
E(\la)=\lf\{ y=(z,w)\in B_1(0)\quad\mbox{ s. t. }\quad|w|\le\la \rg\}
\ee
We show (\ref{V.a5}) argueing by contradiction. First of all from the proof
of lemma~\ref{lm-IV.3} that we apply to $x$ we have the fact that for any $r\in[r_x,\frac{r_{\ep,x_0}}{2}]$
\be
\label{V.a6}
{\mathcal F}(C_{r,x}\res B_1(0)-\pi Q [D_0])\le K\ \sqrt{\al}\quad.
\ee
Let $\om_\al=\chi_\al\ \om=\chi\lf(\frac{|w|}{\al^\frac{1}{4}}\rg)\om$ where $\chi$ is a smooth cut-off function
on ${\R}_+$ satisfying $\chi\equiv 1$ on $[0,1/2]$ and $\chi\equiv 0$ in $[1,+\infty)$.
Let $S$ and $R$ be a 3 and a $2-$current satisfying $(C_{r,x}\res B_1(0)-\pi Q [D_0])=\p S+R$
with $M(S)+M(R)\le 2 K\sqrt{\al}$. We have
\[
\begin{array}{l}
\lf|(C_{r,x}\res B_1(0)-\pi Q [D_0])(\om_\al)\rg|=\lf|S(\om\wedge d\chi_\al)+R(\om_\al)\rg|\\[5mm]
\ds\quad\le
\|\nabla\chi_\al\|_\infty\ \|\om_\infty\|\ M(S)+\|\om_\al\|_\infty\ M(R)\le K\al^\frac{1}{4}\quad .
\end{array}
\]
Thus we get in particular 
\be
\label{V.ab6}
\begin{array}{l}
\ds\pi Q-K\al^\frac{1}{4}\le |C_{r,x}\res B_1(0)(\om_\al)|\le
 M(C_{r,x}\res B_1(0)\cap E(\al^\frac{1}{4})) \|\om_\al\|_\infty\\[5mm]
\ds\quad\le  
M(C_{r,x}\res B_1(0)\cap E(\al^\frac{1}{4}))
\end{array}
\ee
Assuming now there exists 
$y\in ({\mathcal C}_{r,x})_\ast\cap B_1(0)\cap ({\R}^{2p}\setminus E(\al^\frac{1}{16}))$. 
From the monotonicity formula we deduce that 
\be
\label{V.a7}
M(C\res B_{\frac{\al^\frac{1}{16}}{2}r}(y))\ge \frac{\pi}{4}\al^\frac{1}{8}\ r^2\quad .
\ee
Combining (\ref{V.ab6}) and (\ref{V.a7}), we obtain that
\be
\label{V.a8}
M(C\res B_r(x))\ge r^2\lf(\pi Q-K\al^\frac{1}{4}+\frac{\pi}{4}\al^\frac{1}{8}\rg)
\ee
For $\al$ small enough (\ref{V.a8}) contradicts (\ref{VI.1c}) and (\ref{V.a5}) holds true for any
$r\in[r_x,\frac{r_{\ep,x_0}}{2}]$.
From this later fact one deduces (\ref{V.a3}).

It remains to prove (\ref{V.a4}). Again we argue by contradiction. Assume there exists 
$y\in ({\mathcal C}_p\setminus {\mathcal C}_{p-1})\cap \Pi^{-1}(B^2_{r_x}(\Pi(x))\cap B_{r_{\ep,x_0}}(x_0)$.
Because of (\ref{VI.1a}) and since $y\in B_{r_{x_0,\ep}}(x_0)$ we can apply lemma~\ref{lm-IV.2}
  to $y$ in order to deduce 
 that ${\mathcal C}_\star\cap B_{r_x}(x)$
is included in a cone of center $y$, axis parallel to $D_0$ and angle $\ep$. This cone of course contains
$x$ and then we can deduce that
\be
\label{V.a9}
supp(C_{r_x,x}\res B_1(0))\subset E(4\ep)\quad .
\ee
(The notation $E(\la)$ is introduced in (\ref{V.ab5})). We have $\p C_{r,x}\res B_1(0)=0$ moreover, because
of (\ref{V.a6}), for $\al$ small enough we deduce that the intersection number of $C\res B^2_{r_x}(\Pi(x))\times
 B^{2p-2}_{r_x}(0)$ with any vertical current $\Pi^{-1}(z)$ for $z\in B^2_{r_x}$ is $Q$.
 Combining this fact with (\ref{V.a9}), using Fubini, one deduces that
\be
\label{V.a10}
M(C_{r_x,x}\res B_1(0))\ge \pi Q-O(\ep^2)\quad .
\ee
For $\ep$ small enough compare to $\al$ we get a contradiction while comparing (\ref{V.a10}) and (\ref{VI.1c})
and (\ref{V.a4}) is proved. This concludes the proof of lemma~\ref{lm-VI.1}.\cqfd

In the following second lemma of this chapter, we show that the covering 
$(B^2_{r_x}(\Pi(x)))_{x\in{\mathcal C}_{p-1}\cap B^2_{\ep^2 r_{\ep,x_0}}(x_0)}$ of 
$\Pi({\mathcal C}_{p-1}\cap B^{2p}_{\ep^2 r_{\ep,x_0}}(x_0)$ has the ``Whitney'' property : two balls
intersecting each-other have comparable size. From now on we adopt the following  notation granting the fact 
that $\al$ and $\ep$ are fixed small enough for the constraint mentionend above to be fulfilled :
\be
\label{V.a11}
\rho_{x_0}:=\ep^2 r_{\ep,x_0}\quad .
\ee 
Precisely we have.
\begin{Lm}
\label{lm-V.2}
{\bf (Whitney property of the covering.)}
There exists  $\gamma>0$ depending only on $Q$ such that, given 
$x_0\in {\mathcal C}_p\setminus {\mathcal C}_{p-1}$ whose tangent cone is $Q[D_0]$ and let 
$(B^2_{r_x}(\Pi(x)))$ for $x\in{mathcal C}_{p-1}\cap B^2_{\rho_{x_0}}(x_0)$ the covering 
of $\Pi({mathcal C}_{p-1}\cap B^{2p}_{\rho_{x_0}}(x_0)$  described above, assuming for
some $x,y\in{mathcal C}_{p-1}\cap B^2_{\rho_{x_0}}(x_0)$ 
\[
B^2_{r_x}(x)\cap B^2_{r_y}(y)\ne\emptyset\quad ,
\]
then
\be
\label{V.a12}
r_{x}\ge \al^\gamma r_{y}\quad .
\ee
\end{Lm}
{\bf Proof of lemma~\ref{lm-V.2}.}
This lemma is again a consequence of the upper and lower-epiperimetric inequalities.
Assume for instance that $r_x\le r_y$. From (\ref{V.aa3}) we have
\be
\label{V.ab12}
{\mathcal F}(C_{r_y,y}\res B_1(0)- Q\ [D_0])\ge \al\quad .
\ee 
Which implies that forall $r\le r_y$
\be
\label{V.ab13}
{\mathcal F}(C\res B^2_{r}(z_y)\times B^{2p}_{\rho_{x_0}}(0)- Q[B^2_{r}(z_y)\times\{0\}])\ge\al r^3
\ee
where $y=(z_y,w_y)$. Since $|z_x-z_y|\le 2\max\{r_x,r_y\}=2r_y$ and $B_{r_y}(z_y)\subset B_{3 r_y}(z_x)$,
 (\ref{V.ab13}) implies that 
\be
\label{V.ab14}
{\mathcal F}(C\res B^2_{4r_y}(z_x)\times B^{2p}_{\rho_{x_0}}(0)- Q[B^2_{4r_y}(z_x)\times\{w_y\}])\ge 
\frac{\al}{3}r_y^3
\ee
This passage from (\ref{V.ab13}) to (\ref{V.ab14}) is obtained by applying some Fubini type argument.
Indeed, let $A=C\res B^2_{4r_y}(z_x)\times B^{2p}_{\rho_{x_0}}(0)- Q[B^2_{4r_y}(z_x)\times\{w_y\}]$
and let $S$ and $R$ such that $A=\p S+R$ and $M(S)+M(R)^\frac{3}{2}\le 2 {\mathcal F}(A)$. For almost every $r$
in $[r_y/2,r_y]$ we have 
\[
\p(S\res B^2_{r}(z_y)\times B^{2p}_{\rho_{x_0}}(0))=\p S\res B^2_{r}(z_y)\times B^{2p}_{\rho_{x_0}}(0))  
+\lf< S,dist(\cdot,\{z=z_y\},r\rg>
\]
where $\lf< S,dist(\cdot,\{z=z_y\}),r\rg>$ is the slice current between $S$ and the boundary of the cylinder
$B^2_{r}(z_y)\times B^{2p}_{\rho_{x_0}}(0)$ and $dist(\cdot,\{z=z_y\})$ denotes the distance function to
the axis to this cylinder (see \cite{Fe} 4.2.1 pages 395...). Thus
\be
\label{V.ac15}
\begin{array}{l}
\ds A\res B^2_{r}(z_y)\times B^{2p}_{\rho_{x_0}}(0)\\[5mm]
\ds=\p(S\res B^2_{r}(z_y)\times B^{2p}_{\rho_{x_0}}(0))-
\lf< S,dist(\cdot,\{z=z_y\},r\rg> + R\res B^2_{r}(z_y)\times B^{2p}_{\rho_{x_0}}(0)
\end{array}
\ee
We have, see \cite{Fe} 4.2.1 page 395,
\be
\label{V.ac16}
\int_{\frac{r_y}{2}}^{r_y}M\lf(\lf< S,dist(\cdot,\{z=z_y\},r\rg>\rg)\le 
M\lf(S\res (B^2_{r_y}\setminus B^2_{\frac{r_y}{2}})\times B^{2p}_{\rho_{x_0}}(0)\rg)
\ee
Using Fubini theorem we may then find $r=r_1\in [r_y/2,r_y]$ such that 
\be
\label{V.ac17}
M\lf(\lf< S,dist(\cdot,\{z=z_y\},r_1\rg>\rg)\le \frac{2}{r_y}M\lf(S\res (B^2_{r_y}\setminus B^2_{\frac{r_y}{2}})\times B^{2p}_{\rho_{x_0}}(0)\rg)\le\frac{2}{r_y}M(S) 
\ee
Combining (\ref{V.ac15}), (\ref{V.ac16}) and (\ref{V.ac17}) we deduce that 
\be
\label{V.ac18}
\begin{array}{l}
\ds {\mathcal F}(C\res B^2_{r_1}(z_y)\times B^{2p}_{\rho_{x_0}}(0)- Q[B^2_{r_1}(z_y)\times\{0\}])
={\mathcal F}(A\res B^2_{r}(z_y)\times B^{2p}_{\rho_{x_0}}(0))\\[5mm]
\ds\quad \le M(S)+(M(R)+\frac{2}{r_y} M(S))^\frac{3}{2}\quad .
\end{array}
\ee
Thus, combining (\ref{V.ab13}) for $r=r_1$ and (\ref{V.ac18}), we have 
$M(S)+(M(R)+\frac{2}{r_y} M(S))^\frac{3}{2}\ge\al r_y^3$ and since 
${\mathcal F}(C\res B^2_{4r_y}(z_x)\times B^{2p}_{\rho_{x_0}}(0)- Q[B^2_{4r_y}(z_x)\times\{w_y\}])\ge
\frac{1}{2}\lf[M(S)+M(R)^\frac{3}{2}\rg]$, we obtain (\ref{V.ab14}).
Therefore we deduce that
\be
\label{V.ab}
{\mathcal F}(C_{4r_y,x}\res B_1(0)-Q\ [D_0])\ge \frac{1}{3\times 4^3}\al\quad .
\ee
 Let $\frac{\rho_{x_0}}{\ep^2}>s_x>r_x$
such that
\be
\label{V.a13}
M(C_{s_x,x}\res B_1(0))=\pi Q\quad .
\ee
Because of (\ref{VI.1a}), argueing like in the prof of lemma~\ref{lm-IV.2}, between {\ref{IV.1}) and (\ref{IV.3}), we have
\be
\label{V.ab13c}
\forall s_x\le r\le \frac{\rho_{x_0}}{\ep^2}\quad\quad{\mathcal F}(C_{r,x}\res B_1(0)-Q\ [D_0])\le K\ep\quad .
\ee
Assuming, as we did above that $\al>>\ep$, comparing (\ref{V.ab}) and (\ref{V.ab13}) we deduce that $s_x>4 r_y$.
Let $\la=4\frac{r_x}{r_y}$. From the proof of lemma~\ref{lm-IV.3}, in fact from (\ref{V.b4}) precisely,
for any $r\in [r_x,s_x]$ we have
\be
\label{V.a14}
{\mathcal F}(C_{{2r},x}\res B_1(0)-\pi Q [D_0])\le K\ \sqrt{\al}\le \ep_Q^-\quad,
\ee
which means in particular that we are in the position to apply the lower-epiperimetric inequality and the differential
inequality (III.27) in \cite{Ri2} deduced from it. Integrating then this inequality between $r_x$ and $4r_y$  we have
\be
\label{V.12}
\begin{array}{l}
\ds\pi Q-M(C_{4r_y,x}\res B_1(0))\le\lf(\frac{1}{2^{2\ep_Q^-}}\rg)^{\log_2\la}\ [\pi Q- M(C_{r_x,x}\res B_1(0))]\\[5mm]
\ds\quad\quad =\lf(\frac{1}{2^{2\ep_Q^-}}\rg)^{\log_2\la}\al\quad .
\end{array}
\ee
Using now (III.28) from \cite{Ri2}, we deduce from (\ref{V.12})
\be
\label{V.13}
\begin{array}{l}
{\mathcal F}((C_{s_x,x}-C_{4r_y,x})\res B_1(0))\le\sqrt{\pi Q-M(C_{4r_y,x}\res B_1(0))}\\[5mm]
\quad\quad =\lf(\frac{1}{2^{\ep_Q^-}}\rg)^{\log_2\la}\al^\frac{1}{2}\quad. 
\end{array}
\ee
Combining (\ref{V.ab13c}) and (\ref{V.13}) one gets that
\be
\label{V.15}
{\mathcal F}(C_{4r_y,x}\res B_1(0)-Q[D_0])\le \lf(\frac{1}{2^{\ep_Q^-}}\rg)^{\log_2\la}\al^\frac{1}{2}+K\ep\quad .
\ee
Comparing (\ref{V.ab}) and (\ref{V.15}) we obtain
\be
\label{V.16}
\frac{1}{3\times 4^3}\al\le\lf(\frac{1}{2^{\ep_Q^-}}\rg)^{\log_2\la}\al^\frac{1}{2}\quad .
\ee
Since $K\ep \le \frac{1}{6\times 4^3}\al$ we have $\frac{1}{6\times 4^3}\al\le\lf(\frac{1}{2^{\ep_Q^-}}\rg)^{\log_2\la}\al^\frac{1}{2}$.
Taking the log of this last inequality we obtain
\[
\frac{\log\la}{\log 2}\ep_Q^-\log\frac{1}{2}+\frac{1}{2}\log\al\ge\log\al -\log(6\times 4^3)
\]
Thus
\[
\frac{1}{2}\log\frac{1}{\al}+\log(6\times 4^3)\ge\ep_Q^-\log\la
\] 
By taking $\al$ small enough, we may always assume that $\log\frac{1}{\al}\ge 4 \times \log(6\times 4^3)$ and we finally get that
\[
\frac{1}{\ep_Q^-}\log\frac{1}{\al}\ge \log \la
\]
Which leads to the desired inequality (\ref{V.a12}) and lemma~\ref{lm-V.2} is proved.\cqfd

\medskip

\noi {\bf Constructing a partition of unity adapted to the covering.}

\medskip

From the covering $(B^2_{r_x}(x))$ for $x\in {\mathcal C}_{p-1}\cap B_{\rho_{x_0}}(x_0)$ of $\Pi({\mathcal C}_{p-1}\cap B_{\rho_{x_0}}(x_0))$
we extract a Besicovitch covering $(B^2_{r_{x_i}}(x_i))$ for $i\in I$ ($I$ is a countable set) of $\Pi({\mathcal C}_{p-1}\cap B_{\rho_{x_0}}(x_0))$
that is a covering such that
\be
\label{V.zz1}
\forall z\in B^2_{\rho_{x_0}}(x_0)\quad\quad Card\lf\{i\in I\quad\mbox{ s. t. }\quad z\in B^2_{r_{x_i}}(x_i)\rg\}\le n\quad ,
\ee
where $N$ is some  universal number (see \cite{Fe}). To simplify the notation we will simply write $r_i$ for $r_{x_i}$.
Remark that since balls intersecting each-other have comparable size
(see lemma~\ref{lm-V.2}), each ball $B_{r_i}^2(z_i)$ intersects a uniformly bounded number of other balls : there
exists $N_{Q,\al}$ such that
\be
\label{V.16a}
\forall i\in I\quad\quad \mbox{Card}\lf\{j\in I\quad\mbox{ s. t. }\quad B^2_{r_j}(z_j)\cap B^2_{r_i}(z_i)\ne\emptyset\rg\}\le N_{Q,\al}\quad .
\ee 
 We now construct a partition of unity adapted to a slightly modified covering.
Considering the covering $(B^2_{r_{z_i}}(z_i))$ for $i\in I$ ($I$ is a countable set) of 
$\Pi({\mathcal C}_{p-1}\cap B^{2p}_{\rho_{x_0}}(x_0))$, we can apply lemma~\ref{lm-V.3} and obtain $\delta$ depending
on $\al$ and $Q$ such that (\ref{V.ax3}) holds true for some $P\in {\N}$. Let $i\in I$ we can deduce
from (\ref{V.ax3}) and (\ref{V.a12}) that the radii $r_j$
 of balls $B^2_{r_j}(z_j)$ intersecting $B^2_{r_i(1+\delta)}(z_i)$ satisfy 
$\al^{\gamma P} r_i\le r_j\le \al^{-\gamma P} r_i$. From this later fact we deduce that there exists a
number $M\in {\N}$ depending only on $\al$ and $Q$ such that
\be
\label{V.ax9}
\mbox{Card}\lf\{j\in I\quad
\mbox{ s. t. }B_{r_{i}}(z_{i})\cap B_{(1+\delta) r_{j}}(z_{j})\ne\emptyset\rg\}\le M\quad .
\ee
Indeed, assuming $B_{r_i}(z_i)\cap B_{r_j}(z_j)=\emptyset$,
 if $B_{r_{i}}(z_{i})\cap B_{(1+\delta) r_{j}}(z_{j})$ we just have seen that 
$\al^{\gamma P} r_j\le r_i\le \al^{-\gamma P} r_j$ :  the two radii have comparable size which is of course 
also comparable with the distance $|z_i-z_j|$. From (\ref{V.zz1}) it is then clear that the number of such
ball $B_{r_j}(z_j)$ is bounded by a constant depending only of the variables $\al$ and $Q$. 
It is now not difficult to deduce that $(B^2_{r_i(1+\frac{\delta}{2})})_{i\in I}$ realizes a locally finite
covering of $\Pi({\mathcal C}_{p-1}\cap B^{2p}_{\rho_{x_0}}(x_0))$ satisfying
\be
\label{V.ax10}
\forall i\in I\quad\quad\ \mbox{Card}\lf\{j\in I\quad
\mbox{ s. t. }B_{(1+\frac{\delta}{2})r_{i}}(z_{i})\cap B_{(1+\frac{\delta}{2}) r_{j}}(z_{j})\ne\emptyset\rg\}\le
 M+P\quad .
\ee
Indeed assuming for instance that $r_i\ge r_j$, then 
$B_{(1+\frac{\delta}{2})r_{i}}(z_{i})\cap B_{(1+\frac{\delta}{2}) r_{j}}(z_{j})\ne\emptyset$ implies clearly 
that $B_{(1+\delta)r_i}(z_i)\cap B_{r_j}(z_j)\ne \emptyset$ and the numer of such a $j$ is controled by
$P$ (see \ref{V.ax3}), whereas if $r_i\le r_j$, 
$B_{(1+\frac{\delta}{2})r_{i}}(z_{i})\cap B_{(1+\frac{\delta}{2}) r_{j}}(z_{j})\ne\emptyset$ implies clearly 
that $B_{r_i}(z_i)\cap B_{r_j(1+\delta)}(z_j)\ne \emptyset$ and the numer of such a $j$ is controled by
$M$ (see \ref{V.ax9}). Thus (\ref{V.ax10}) holds true.
 And we shall use from now on the notation
\be
\label{V.az10}
\forall i\in I\quad\quad \rho_i:=r_i\lf(1+\frac{\delta}{2}\rg)\quad .
\ee

For any $i\in I$ we define $\chi_i$ to be a smooth non-negative function satisfying
\begin{itemize}
\item[i)]
 \be
\label{V.ax11}
\chi_i\equiv 1\quad\mbox{ in }B^2_{r_i}(z_i)\quad .
\ee
\item[ii)]
\be
\label{V.ax12}
\chi_i\equiv 0\quad\mbox{ in }{\R}^2\setminus B^2_{(1+\frac{\delta}{2}r_i)}(z_i)\quad
\ee
\item[iii)]
\be
\label{V.ax15}
\forall k\in {\N}\quad\quad \|\nabla^k \chi_i\|_\infty\le \frac{K_k}{r_i^k}\quad,
\ee
where $K_k$ depends only on $k$ and $Q$.
\end{itemize}
We define now
\be
\label{V.ax14}
\vphi_i:=\frac{\chi_i}{\sum_{i\in I}\chi_i}\quad.
\ee
It is clear that $(\vphi_i)$ defines a partition of unity adapted to $B^2_{(1+\frac{\delta}{2}r_i)}(z_i)$
and satisfying the following estimates 
\be
\label{V.ax13}
\forall k\in {\N}\quad\quad \|\nabla^k \vphi_i\|_\infty\le \frac{K_k}{r_i^k}\quad,
\ee
where $K_k$ depends only on $k$ and $Q$.

\section{The approximated average curve.}
\reset

This chapter is another step towards the proof that ${\mathcal P}_{Q-1}\Longrightarrow {\mathcal P}_Q$ which
goes until chapter VIII.
We then assume that ${\mathcal P}_{Q-1}$ holds (or that $Q=1$).
Again in this part we consider the difficult case which is the case where we are blowing-up the current at
a point $x_0$ of multiplicity $Q>1$ whose tangent cone $C_{0,x_0}$ is $Q$ times the integral
over the flat disk $D_0$ given by $w_i=0$ for $i\cdots Q-1$ (we use the system
of coordinates introduced in the begining of chapter II) and where $x_0$ belongs to the closure of 
${\mathcal C}_{Q-1}$. 
The purpose of this chapter is to approximate first our current over each ball of the covering introduced
in the previous section $C\res \Pi^{-1}(B^2_{r_i}(z_i))$ by a $Q-$valued graph $\{a_i^k\}_{k=1\cdots Q}$
over $B^2_{r_i}(z_i)$ which is almost $J-$holomorphic ($J_{x_i}-$holomorphic in fact where 
$x_i\in {\mathcal C}_\ast$ and $\Pi(x_i)=z_i$)
and glueing the average curves $\ti{a}_i=\frac{1}{Q}\sum_{k=1}^Q a_i^k$ of each of these $J_{x_i}-$holomorphic 
$Q-$valued graphs together we shall produce
a single-valued graph $\ti{a}$ over $B^2_{\rho_{x_0}}(x_0)$ which approximates $C$ and for which we will study 
regularity properties that will be used in the following chapter VII devoted to the unique continuation argument .
Finally in the second subsection of this chapter we construct new coordinates adapted to the average curve.

\subsection{Constructing the average curve.}

Let $\rho_{x_0}$ given by (\ref{V.a11}) and let $(B^2_{\rho_i}(z_i))_{i\in I}$ be the Besicovitch-Whitney
 covering of $\Pi({\mathcal C}_{Q-1}\cap B_{\rho_{x_0}}(x_0))$ obtained at the end of the previous chapter.
As we have seen above, for any $i\in I$,
 ${\mathcal C}_\ast\cap\Pi^{-1}(B^2_{\rho_i}(z_i))\subset {\mathcal C}_{Q-1}\cap
B^{2}_{\rho_i}(z_i)\times B^{2p-2}_{2\rho_i}(w_i)$ where $x_i=(z_i,w_i)$ is in ${\mathcal C}_\ast$
(see lemma~\ref{lm-VI.1}). For convenience we shall adopt the following notation
\be
\label{VII.1}
N_{r}^i:=B^{2}_{r}(z_i)\times B^{2p-2}_{2\rho_i}(w_i)\quad .
\ee
Assuming ${\mathcal P}_{Q-1}$, $C\res N_{2\rho_i}$ is a $J-$holomorphic curve : there exists a smooth Riemannian 
surface and a smooth $J-$holomorphic map
\be
\label{VII.2}
\begin{array}{rcl}
\ds \Psi_i\ :\ \Sigma_{2,i} &\ds \longrightarrow &\ds N^i_{2\rho_i}\\[5mm]
\ds \xi &\ds \longrightarrow &\ds \Psi_i(\xi)
\end{array}
\ee
such that $\Psi_\ast[\Sigma_{2,i}]=C\res N_{2\rho_i}$. Let $H^0_{\pm}(\Sigma_{2,i})$ be the sets
respectively of holomorphic and antiholomorphic functions on $\Sigma_{2,i}$
. We introduce now $\eta_i$ the map from $\Sigma_{2,i}$ into ${\R}^{2p-2}$
chosen such that the perturbation $\Psi_i+\eta_i$ is $J_{x_i}-$holomorphic, precisely $\eta_i$ is given
by
\be
\label{VII.3}
\lf\{
\begin{array}{l}
\ds\frac{\p}{\p\xi_1}(\Psi_i+\eta_i)+J_{x_i}\frac{\p}{\p\xi_2}(\Psi_i+\eta_i)=0\quad\mbox{ in }\Sigma_{2,i}\\[5mm]
\ds \forall h\in H(\Sigma_{2,i})\quad\quad\int_{\p \Sigma_{2,i}}h\ d\eta_i=0
\end{array}
\rg.
\ee
where $(\xi_1,\xi_2)$ are local coordinates on $\Sigma_{2,i}$ compatible with the complex structure.
The existence of $\eta_i$ is justified few lines below.
 To this aim we shall make use of the following notations.
 Since $J$ is smooth, Local inversion theorem gives the 
existence of a smooth map $\La\ :\  B_{\rho_{x_0}}^{2p}(x_0)\times {\R}^{2p}\longrightarrow {\R}^{2p}$
- for $\rho_{x_0}$ chosen small enough such 
that
\begin{itemize}
\item[i)]
\be
\label{VII.3a}
\La_x:=\La(x,.)\quad\quad\mbox{ is a linear isomorphism of }{\R}^{2p}\quad,
\ee
\item[ii)] 
\be
\label{VII.3b}
\La_{x_0}=id\quad,
\ee
\item[iii)]
\be
\label{VII.3c}
\forall x\in B_{\rho_{x_0}}^{2p}(x_0)\quad\quad 
J_{x_0}=
\lf(
\begin{array}{ccccc}
0 & -1 & 0\cdots & 0 &0\\[5mm]
1 & 0  & 0\cdots & 0 &0\\[5mm]
. & . & . & . & .\\[5mm]
. & . & . & . & .\\[5mm]
. & . & . & . & .\\[5mm]  
0 & 0 & .\cdots & 0 & -1\\[5mm]
0 & 0 & .\cdots & 1 & 0
\end{array}
\rg)=\La_x\  J_{x}\ \La_x^{-1}
\ee
\end{itemize}
We shall denote by $(z^i=x^i_1+ix^i_2,w^i_1=x^i_3+ix^i_4,\cdots,w^i_p=x^i_{2p-1}+ix^i_{2p})$ the following complex coordinates in $N_{2\rho_i}$
\be
\label{VII.3d}
\lf(
\begin{array}{c}
x^i_1\\[5mm]
x^i_2\\[5mm]
.\\[5mm]
.\\[5mm]
.\\[5mm]
x^i_{2p-1}\\[5mm]
x^i_{2p}
\end{array}
\rg)= \La_{x_i}\cdot\lf(
\begin{array}{c}
x_1\\[5mm]
x_2\\[5mm]
.\\[5mm]
.\\[5mm]
.\\[5mm]
x_{2p-1}\\[5mm]
x_{2p}
\end{array}
\rg)-x_i 
\ee
We will also denote by $\Pi_i$ the map that assign to any point $x$ in $N^i_{\rho_i}$ the complex coordinate
$z^i$ and by $D^i$ we denote the $J_{x_i}-$holomorphic 2-disk 
\be
\label{VII.3e}
D^i:=\lf\{x\ ;\ \forall k=1\cdots p-1\quad
w^i_k=0\rg\}=
\La_{x_i}^{-1} D_0
\ee  
Using these complex coordinates in $N_{2\rho_i}^i$, (\ref{VII.3}) means
\be
\label{VII.4}
\lf\{
\begin{array}{l}
\ds \ov{\p}\eta_i=-\ov{\p}\Psi_i\quad\mbox{ in }\Sigma_{2,i}\\[5mm]
\ds \forall h\in H(\Sigma_{2,i})\quad\quad\int_{\p \Sigma_{2,i}}h\ d\eta_i=0\quad.
\end{array}
\rg.
\ee 
The existence and uniqueness of $\eta_i$ is given by proposition A.3 of \cite{Ri3}.
Since $\Psi_i$ is $J-$holomorphic we have $\p_{\xi_1}\Psi_i+J(\Psi_i(\xi))\p_{\xi_2}\Psi_i=0$, thus
$|\p_{\xi_1}\Psi_i+J(x_i)\p_{\xi_2}\Psi_i|\le |J(\Psi_i(\xi))-J(x_i)|\ |\nabla\psi|$. Combining this fact with
the second part of proposition A.3 (i.e. estimate (A.13) of \cite{Ri3}) we obtain
\be
\label{VII.5}
\int_{\Sigma_{2,i}}|\nabla\eta_i|^2\le K r_i^2\int_{\Sigma_{2,i}}|\nabla\Psi_i|^2\le K r_i^4\quad .
\ee
For $\la\le 2$, we denote by $\Sigma_{\la,i}$  the surface $\Sigma_{\la,i}=\Sigma_{2,i}\cap 
\Psi^{-1}_i(N^i_{\la\rho_i})$ so that
\be
\label{VII.6}
{\Psi_i}_\ast[\Sigma_i]=C\res B^2_{\rho_i}(z_i)\times B^{n-2}_{\rho_{x_0}}(0)\quad .
\ee
We then prove in \cite{Ri3} the following lemma
\begin{Lm}
\label{lm-VII.1}
Under the above notations one has
\be
\label{VII.7}
\|\eta_i\|_{L^\infty(\Sigma_{\frac{3}{2},i})}\le K r_i^2
\ee
where $K$ is a constant depending only on $\|\nabla J\|_\infty$ and the choice of 
$\al$ made in the previous chapter.
\end{Lm}

Consider now the $J_{x_i}-$holomorphic curve $C^i$ given by the image by $\Psi_i+\eta_i$ of $\Sigma_{\frac{3}{2},i}$.
Since $\p {\Psi_i}_\ast[\Sigma_{\frac{3}{2},i}]$ is supported in $\Pi^{-1}(\p B^2_{\frac{3}{2}\rho_i})$. Now from
Lemma~\ref{lm-VII.1} we know that $|\eta_i|_\infty\le Cr_i^2$ therefore $\p {\Psi_i+\eta_i}_\ast[\Sigma_{\frac{3}{2},i}]$
is supported in an $r_i^2$ neighborhood of $\Pi^{-1}(\p B^2_{\frac{3}{2}\rho_i}(0))$ and thus (for $r_i$ small enough
: that holds if $\ep$ has been chosed small enough in section VI) we have that 
${\Psi_i}_\ast[\Sigma_{\frac{3}{2},i}])$ is a cycle in $ \Pi_i^{-1}(B^2_{\frac{5}{4}\rho_i}(0))$ and the part of the image of $\Sigma_{\frac{3}{2},i}$
included in $\Pi^{-1}_i(B^2_{\frac{5}{4}\rho_i}(0))$ by $\Psi_i+\eta_i$ is a $J_{x_i}-$holomorphic cycle and therefore it is a $Q-$valued
graph over $D^i$ for the complex coordinates given by $(z^i,w^i)$. We denote by $\{a^i_k\}_{k=1\cdots Q}$ this $Q-$valued graph
(i.e. $a_k^i(z^i_0)$ are the $w^i$ coordinates, in the chart $(z^i,w^i)$, of the $Q$ intersection points between the $J_{x_i}-$holomorphic
curve $\Psi_i+\eta_i(\Sigma_{\frac{3}{2},i})$ and the $J_{x_i}-$holomorphic submanifold given by $z^i=z^i_0$).
We now define $\ti{C}^i$ to be the $J_{x_i}$ holomorphic curve in $\Pi_i^{-1}(B^2_{\rho_i})$ given by 
\be
\label{VII.27}
\ti{C}^i:=\lf\{x=\La_{x_i}^{-1}\lf( (z^i,\ti{a}^i_i(z_i)=\frac{1}{Q}\sum_{k=1}^Q a_k^i(z^i))+x_i\rg)\quad\forall z^i\in B^2_{\frac{5}{4}\rho_i}(0)\rg\}\quad .
\ee
Observe that 
\be
\label{VII.28}
\frac{\p}{\p z^i}\ti{a}^i_i=0\quad\quad\mbox{ in }{\mathcal D}'(B_{\frac{5}{4}\rho_i}^2(0))\quad .
\ee
Moreover, the conformal invariance of the Dirichlet energy gives
\be
\label{VII.29}
\begin{array}{l}
\ds\int_{B_{\frac{5}{4}\rho_i}^2(0)}\sum_{k=1}^Q|\nabla a^i_k|^2(z^i)\ dz^i\wedge d\ov{z}^i\\[5mm]
\ds\quad\quad=\int_{(\Psi_i+\eta_i)^{-1}(C_i\cap\Pi^{-1}_i(B^2_{\frac{5}{4}\rho_i}(0)))}
|\nabla(\Psi_i+\eta_i)|^2(\xi)\ d\xi\le K r_i^2\quad .
\end{array}
\ee
We then deduce that
\be
\label{VII.30}
 \int_{ B_{\frac{5}{4}\rho_i}^2(0)}|\nabla\ti{a}^i_i|^2(z^i)\ dz^i\le K r_i^2\quad.
\ee
Combining (\ref{VII.29}) and (\ref{VII.30}) and using standard elliptic estimates we get that for any $l\in {\N}$
\be
\label{VII.31}
\|\nabla^l\ti{a}^i_i\|_{L^\infty(B^2_{\frac{6}{5}\rho_i}(0))}\le K_l r_i^{-l+1}\quad .
\ee
The subscript $i$ in the notation $\ti{a}^i_i$ is here to recall that we express $\ti{C}^i$ as a graph in the 
$(z^i,w^i)$ coordinates. The same $J_{x_i}-$holomorphic curve $\ti{C}^i$ can also, due to (\ref{VII.31}), be expressed as a graph
in a neighborhing system of coordinate $(z^j,w^j)$ where $B_{\rho_i}(z_i)\cap B_{\rho_j}(z_j)\ne\emptyset$ (indeed
the passage from $(z^i,w^i)$ to $(z^j,w^j)$ is given by a transformation matrix in ${\R}^{2p}$ close to the identity at a distance
of the order $r_i$). In such system of coordinates $(z^j,w^j)$, we shall denote $\ti{a}_j^i(z^j)$ the graph corresponding to $\ti{C}^i$.

Since $\ti{C}^i$ is a graph over $w^i=0$ given by $(z^i,\ti{a}^i_i(z^i))$ whose gradient is bounded (see (\ref{VII.31})),
and since the passage
from the $(z,w)$ coordinates to $(z^i,w^i)$ coordinates is given by a transformation $\La_{x_i}$ whose distance to the
identity is bounded by $|x_i|$ that tends to zero, $\ti{C}^i$ is then also realised by a graph over $B^2_{\frac{7}{6}\rho_i}(\Pi(x_i))$ that we shall now 
denote $(z,\ti{a}_i(z))$ :
\be
\label{VII.59a}
\ti{C}^i\res \Pi^{-1}(B^2_{\frac{7}{6}\rho_i}(\Pi(x_i)))=(z,\ti{a}_i(z))_\ast[B^2_{\frac{7}{6}\rho_i}(\Pi(x_i))]\quad.
\ee
Consider now $i$ and $j$ such that $B_{\rho_i}(z_i)\cap B_{\rho_j}(z_j)\ne\emptyset$. We shall compare $\ti{a}_i$ and $\ti{a}_j$ in 
$\Pi^{-1}(B_{\rho_i}(z_i)\cap B_{\rho_j}(z_j))$. Precisely we have the following lemma
\begin{Lm}
\label{lm-VII1a}
Under the above notations one has
\be
\label{VII.78}
\forall l\in {\N}\quad\quad\|\nabla^l(\ti{a}_i-\ti{a}_j)\|_{L^\infty(B^2_{\rho_i}(z_i)\cap B^2_{\rho_j}(z_j))}\le K_l
\rho_i^{2-l}\quad .
\ee
\end{Lm}
{\bf Proof of lemma~\ref{lm-VII1a}.}

 First of all we compare $C_i$ and $C_j$ in $\Pi^{-1}(B_{\rho_i}(z_i)\cap B_{\rho_j}(z_j))$.
We can allways assume that $\Sigma_{2,i}$ and $\Sigma_{2,j}$ are part of a same Riemannian surface $\Sigma$ with a joint parametrization
$\Psi=\Psi_i$ on $\Sigma_{2,i}$ and $\Psi=\Psi_j$ on $\Sigma_{2,j}$ and such that
 ${\Psi}_\ast[\Sigma]=C\res N^i_{2\rho_i}\cup N^j_{2\rho_j}$. Let $\Sigma^{ij}:=\Psi^{-1}(supp(C)\cap N^i_{2\rho_i}\cap N^j_{2\rho_j}$. We consider the following mapping
\be
\label{VII.32}
\begin{array}{rcl}
\ds\Xi^{ij}&\ds \Sigma^{ij}\times[0,1]\longrightarrow &\ds N^i_{3\rho_i}\cap N^j_{3\rho_j}\\[5mm]
 &\ds (\xi,t)\longrightarrow &\ds \Psi(\xi)+t\eta_j(\xi)+(1-t)\eta_i(\xi)\quad .
\end{array}
\ee
Clearly for any $\la<\frac{3}{2}$
\be
\label{VII.33}
\p {\Xi^{ij}}_\ast[\Sigma^{ij}]\times[0,1]\res N^i_{\la\rho_i}\cap N^j_{\la\rho_j}=
C^j-C^i\res N^i{\la\rho_i}\cap N^j_{\la\rho_j}\quad .
\ee
We have
\be
\label{VII.34}
M({\Xi^{ij}}_\ast[\Sigma^{ij}]\times[0,1])=\int_0^1\int_{\Sigma^{ij}}J_3\Xi^{ij}\quad,
\ee
where $(J_3\Xi^{ij})^2$ is the sum of the squares of the determinants of the $3\times 3$ submatrices
of $D\Xi^{ij}$. Clearly
\be
\label{VII.35}
|J_3\Xi^{ij}|(\xi,t)\le K\ \lf[\|\eta_i\|_\infty+\|\eta_j\|_\infty\rg]\ \lf[|\nabla \Psi|^2(\xi)+
|\nabla\eta_i|^2(\xi)+|\nabla\eta_j|^2(\xi)\rg]\quad .
\ee
Combining Lemma~\ref{lm-V.2},  Lemma~\ref{lm-VII.1}, (\ref{VII.34}) and (\ref{VII.35}), we get that
for any $\la<\frac{3}{2}$
\be
\label{VII.36}
M({\Xi^{ij}}_\ast[\Sigma^{ij}]\times[0,1]\res N^i{\la\rho_i}\cap N^j_{\la\rho_j})\le K\ r_i^4 \quad .
\ee
Therefore, combining (\ref{VII.33}) and (\ref{VII.36}), using a standard slicing and Fubini type argument,
 we may find $\la\in (\frac{5}{4},\frac{3}{2}$ such that
\be
\label{VII.37}
{\mathcal F}((C^i-C^j)\res N^i_{\la\rho_i}\cap N^j_{\la\rho_j})\le K r_i^4\quad .
\ee
We shall now compare $\ti{C}^i$ and $\ti{C}^j$. Denote $(x_1^t,x_2^t\cdots x_{2p}^t)$ the coordinates given by
$(x_1^t,x_2^t\cdots x_{2p}^t)^T=\La_{x_t}\cdot [(x_1,x_2\cdots x_{2p})^T-x_i]$ where we keep denoting
$(x_1,x_2\cdots x_2p)$ our original normal coordinates vanishing at the center $x_0$ introduced in (\ref{II.0})
and $\La_{x_t}$ is the transformation matrix  introduced in (\ref{VII.3a}). Observe that with these notations
$(x_1^0,x_2^0\cdots x_{2p}^0)=(x_1^i,x_2^i\cdots x_{2p}^i)$ that $(x_1^1,x_2^1\cdots x_{2p}^1)
=(x_1^j,x_2^j\cdots x_{2p}^j)+\La_{x^j}\cdot(x_i-x_j)$ and that $(x_1^t,x_2^t\cdots x_{2p}^t)$ has been chosed 
in order to vanish at a fixed point $x_i$. Observe also that 
\be
\label{VII.38}
\lf|\frac{d}{dt}\La_{x^t}\rg|\le K r_i\quad,\quad \|\frac{d}{dt}x^t\|_{L^\infty(B^{2p}_{4r_i}(x_i))}+\|\frac{d}{dt}y^t\|_{L^\infty(B^{2p}_{4r_i}(x_i))}\le K r^2_i\quad .
\ee
We also adopt the notations $z^t:=x^t_1+ix^t_2$ and for $k=1\cdots p-1$
$w^t:=x^t_{2k+1}+ix^t_{2k+2}$. Observe then that $z^t=$constant or $w^t_k=$constant are $J_{x^t}-$holomorphic
$2p-2$ submanifolds, or simply complex variety in $({\R}^{2p},J_{x^t})$. In order to
compare $\ti{C}^i$ and $\ti{C}^j$ we shall perturb $\Xi^{ij}$ in the following
way : denote first $\Psi^t$, $\eta_i^t$ and $\eta^t_j$ the maps $\Psi$,$\eta_i$ and $\eta_j$ expressed in the
 coordinates $(z^t,w^t)$, and consider the map $s^t\ :\ (\Sigma')^{ij}\rightarrow {\C}^p$ solving
\be
\label{VII.39}
\begin{array}{l}
\ds\p_{\ov{\xi}}s^t=\p_{\ov{\xi}}(\Psi^t+t\eta^t_j+(1-t)\eta^t_i)\quad\mbox{ in }(\Sigma')^{ij}\\[5mm]
\ds\forall h\in H(\Sigma^{ij})\quad\quad\int_{\p\Sigma^{ij}}s^t\,dh=0\quad.
\end{array}
\ee
where $\Sigma^{ij}:=\Psi^{-1}(supp(C))\cap N^i_{\frac{3}{2}\rho_i}\cap N^j_{\frac{3}{2}\rho_j}$. The existence
and uniqueness of $s^t$ is given by proposition A.3 of \cite{Ri3}
We shall now replace the map $\Xi^{ij}$ on $(\Sigma')^{ij}$ by the map
\be
\label{VII.40}
\begin{array}{rl}
\ds(\Xi')^{ij}\ :\  (\Sigma')^{ij}\times[0,1]\longrightarrow &\ds N^i_{3\rho_i}\cap N^j_{3\rho_j}\\[5mm]
 (\xi,t)\longrightarrow &\ds \La_{x^t}^{-1}\cdot[\Psi^t(\xi)+t\eta^t_j(\xi)+(1-t)\eta^t_i(\xi)-s^t]+x_i\quad .
\end{array}
\ee
Observe that for each $t\in[0,1]$ the map $\Xi'^{ij}(\cdot,t)$ is a $J_{x^t}-$holomorphic curve. Observe also that,
for $t=0$ $s^t=0$ and that for $t=1$, $\p_{\ov{\xi}}(\Psi^1+\eta_j^1)=0$, since $\Psi+\eta_j$ is $J_{x^j}-$holomorphic
and $(z^1,w^1)$ are $J_{x^j}$ coordinates, thus we have also $s^1=0$. One can easily verify, like for 
proving (\ref{VII.5}) that $forall t\in [0,1]$
\[
\int_{(\Sigma')^{ij}}|\p_{\ov{\xi}}(\Psi^t+t\eta^t_j+(1-t)\eta^t_i)|^2\le K r^4_i\quad ,
\]
and using lemma~\ref{lm-VII.1} we have 
$\|s^t\|_\infty((\Sigma'')^{ij})\le K r_i^2$ where $(\Sigma'')^{ij}:=
\Psi^{-1}(supp(C))\cap N^i_{\frac{5}{4}\rho_i}\cap N^j_{\frac{5}{4}\rho_j}$. Therefore for any $\la<\frac{6}{5}$
we have
\be
\label{VII.41}
\p {(\Xi')^{ij}}_\ast[(\Sigma')^{ij}]\times[0,1]\res N^i{\la\rho_i}\cap N^j_{\la\rho_j}=
C^j-C^i\res N^i{\la\rho_i}\cap N^j_{\la\rho_j}\quad .
\ee
We consider now the following interpolation between $\ti{C}^i$ and $\ti{C}^j$ : let $\ti{\Xi}^{ij}$ be the 
following map
\be
\label{VII.42}
\begin{array}{rcl}
\ds\ti{\Xi}^{ij}&\ds \Pi(\Psi((\Sigma'')^{ij}))\times[0,1]\longrightarrow &\ds N^i_{3\rho_i}\cap N^j_{3\rho_j}\\[5mm]
 &\ds (z,t)\longrightarrow &\ds \La_{x^t}^{-1}\cdot[(z,\ti{a}^t(z)]+x_i\quad ,
\end{array}
\ee
where the $p-1$ complex components are given by the slices of $C^t:={(\Xi')^{ij}}_\ast[(\Sigma')^{ij}]\times\{t\}$
by $z^t=\xi$ evaluated on the functions $w_k^t$. Precisely using the notations of \cite{Fe} 4.3
\be
\label{VII.43}
\ti{a}^t_k(z):=\lf<C^t,z^t,z\rg>(w^t_k)\quad .
\ee
It is clear that for any $\la<\frac{6}{5}$
\be
\label{VII.44}
\Txi_\ast[\Pi_i(\Psi((\Sigma')^{ij}))]\times[0,1]\res N^i{\la\rho_i}\cap N^j_{\la\rho_j}=
\ti{C}^j-\ti{C}^i\res N^i{\la\rho_i}\cap N^j_{\la\rho_j}\quad .
\ee
In order to get a bound for ${\mathcal F}(\ti{C}^j-\ti{C}^i\res N^i{\la\rho_i}\cap N^j_{\la\rho_j})$, it remains
 to evaluate the mass of 
$\Txi_\ast[\Pi_i(\Psi((\Sigma')^{ij}))]\times[0,1]\res N^i{\la\rho_i}\cap N^j_{\la\rho_j}$ for
 $\la=\frac{6}{5}$
for instance. We have
\be
\label{VII.45}
|J_3\ti{\Xi}^{ij}|(z,t)\le \lf|\frac{\p}{\p t}\La_{x^t}^{-1}\cdot[(z,\ti{a}^t(z)]\rg|(z,t)\ 
\lf[1+|\nabla_z\ti{a}^t(z)|^2\rg]\quad .
\ee
Because of the same arguments developped to prove (\ref{VII.31}), since $\ti{a}^t(\xi)$ is holomorphic, we have
\be
\label{VII.46}
\|\nabla_{z}\ti{a}^t(z)\|_{L^\infty((\Sigma'')^{ij})}\le K\quad .
\ee
Thus
\be
\label{VII.47}
\begin{array}{l}
\ds M(\Txi_\ast[\Pi\Psi((\Sigma')^{ij})]\times[0,1]\res N^i{\la\rho_i}\cap N^j_{\la\rho_j})\le
 \int_0^1\int_{\Pi(\Psi((\Sigma'')^{ij}))} |J_3\ti{\Xi}^{ij}|\\[5mm]
\ds \quad\le K\ \int_0^1\int_{\Pi(\Psi((\Sigma'')^{ij}))}
\lf|\frac{\p}{\p t}\La_{x^t}^{-1}\cdot[(z,\ti{a}^t(z)]\rg|\ dz\wedge d\ov{z}\wedge dt\\[5mm]
\ds\quad
\le K\ \int_0^1\int_{\Pi(\Psi((\Sigma'')^{ij}))} \lf[r_i |(z,\ti{a}^t(z))|+\lf|\frac{\p\ti{a}^t}{\p t}\rg|\rg]\quad .
\end{array}
\ee
 In one hand, since $\Pi_i(\Psi((\Sigma'')^{ij}))$ $|(z,\ti{a}^t(z))|\le K\ r_i$, we have 
\be
\label{VII.48}
\int_0^1\int_{\Pi(\Psi((\Sigma'')^{ij}))} r_i |(z,\ti{a}^t(z))|\le K r_i^4\quad .
\ee
In the other hand 
\be
\label{VII.49}
\int_0^1\int_{\Pi(\Psi((\Sigma'')^{ij}))}\lf|\frac{\p\ti{a}^t}{\p t}\rg|=
\lim_{N\rightarrow +\infty}\frac{1}{N}\sum_{l=1}^{N-1}
\int_{\Pi(\Psi((\Sigma'')^{ij}))}N\ |\ti{a}^{\frac{l}{N}}(z)-\ti{a}^{\frac{l+1}{N}}(z)|\ dz\wedge d\ov{z}\quad.
\ee
We have 
\be
\label{VII.50}
\begin{array}{l}
\ds|\ti{a}^{\frac{l}{N}}_k(z)-\ti{a}^{\frac{l+1}{N}}_k(z)|
=|<C^{\frac{l}{N}}, z^{\frac{l}{N}},z>(w^{\frac{l}{N}}_k)-
<C^{\frac{l+1}{N}}, z^{\frac{l+1}{N}},z>(w^{\frac{l+1}{N}}_k)|\\[5mm]
\ds\le|<C^{\frac{l}{N}}-C^{\frac{l+1}{N}}, z^{\frac{l}{N}},z>(w^{\frac{l}{N}}_k)|\\[5mm]
\ds+|<C^{\frac{l+1}{N}}, z^{\frac{l}{N}},z>(w^{\frac{l}{N}}_k)-
<C^{\frac{l+1}{N}}, z^{\frac{l+1}{N}},z>(w^{\frac{l}{N}}_k)|\\[5mm]
\ds+|<C^{\frac{l+1}{N}}, z^{\frac{l+1}{N}},z>(w^{\frac{l}{N}}_k-w^{\frac{l+1}{N}}_k)|
\end{array}
\ee
We have to control the sum over $l$ of the integral over $\Pi(\Psi((\Sigma'')^{ij}))$ of the 3 absolute values
in the right-hand-side of (\ref{VII.50}) one by one. For the first term of the r.h.s. of(\ref{VII.50}) we have, 
using \cite{Fe} 4.3.1, since $\|w^{\frac{l}{N}}_k\|_\infty+\|dw^{\frac{l}{N}}_k\|_\infty\le 1$,
\be
\label{VII.51}
\begin{array}{l}
\ds\int_{\Pi(\Psi((\Sigma'')^{ij}))}|<C^{\frac{l}{N}}-C^{\frac{l+1}{N}}, z^{\frac{l}{N}},z>(w^{\frac{l}{N}}_k)|\
 dz\wedge d\ov{z}\\[5mm]
\ds\le Lip(z^{\frac{k}{N}})\ {\mathcal F}_{N^i_{\la\rho_i}\cap N^j_{\la\rho_j}}(C^{\frac{k}{N}}-C^{\frac{k+1}{N}})
\ds\le K\ {\mathcal F}_{N^i_{\la\rho_i}\cap N^j_{\la\rho_j}}(C^{\frac{k}{N}}-C^{\frac{k+1}{N}})
\end{array}
\ee
Similarly as the way we have establishe estimate (\ref{VII.37}) we can show that
\be
\label{VII.52}
{\mathcal F}_{N^i{\la\rho_i}\cap N^j_{\la\rho_j}}(C^{\frac{k}{N}}-C^{\frac{k+1}{N}})\le K\ \frac{1}{N}\ r_i^4\quad .
\ee
Thus
\be
\label{VII.53}
\lim_{N\rightarrow +\infty}\sum_{l=1}^{N-1}
\int_{\Pi(\Psi((\Sigma'')^{ij}))}|<C^{\frac{l}{N}}-C^{\frac{l+1}{N}}, z^{\frac{l}{N}},z>(w^{\frac{l}{N}}_k)|\le K\ r_i^4
\quad .
\ee
For the second term of the r.h.s. of(\ref{VII.50}) we use 4.3.9 (3) of \cite{Fe} and we get
\be
\label{VII.54}
\begin{array}{l}
\ds\int_{\Pi(\Psi((\Sigma'')^{ij}))}\ |<C^{\frac{l+1}{N}}, z^{\frac{l}{N}},z>(w^{\frac{l}{N}}_k)-
<C^{\frac{l+1}{N}}, z^{\frac{l+1}{N}},z>(w^{\frac{l}{N}}_k)|\\[5mm]
\ds K\int_{\frac{l}{N}}^{\frac{l+1}{N}} \ dt\int_{(z^t)^{-1}(\Pi(\Psi((\Sigma'')^{ij})))}
|z^{\frac{l}{N}}-z^{\frac{l+1}{N}}|\ d\|C^{\frac{l+1}{N}}\|\\[5mm]
\ds K\ \frac{r_i^2}{N} M(C^{\frac{l+1}{N}}\res N^i{\la\rho_i}\cap N^j_{\la\rho_j})\\[5mm]
\ds\le K\frac{r_i^4}{N}
\end{array}
\ee
where we have used (\ref{VII.38}). Therefore we obtain
\be
\label{VII.55}
\ds\lim_{l\rightarrow +\infty}\sum_{l}^{N-1}\int_{\Pi(\Psi((\Sigma'')^{ij}))}\ |<C^{\frac{l+1}{N}}, z^{\frac{l}{N}},z>(w^{\frac{l}{N}}_k)-
<C^{\frac{l+1}{N}}, z^{\frac{l+1}{N}},z>(w^{\frac{l}{N}}_k)|\le K r_i^4\quad .
\ee
Finally for the second term of the r.h.s. of(\ref{VII.50}), we use again (\ref{VII.38}) and 4.3.2 (2) of \cite{Fe}
 to obtain that
\be
\label{VII.56}
\begin{array}{l}
\ds\int_{\Pi(\Psi((\Sigma'')^{ij}))}\ |<C^{\frac{l+1}{N}}, z^{\frac{l+1}{N}},z>(w^{\frac{l}{N}}_k-w^{\frac{l+1}{N}}_k)|
\\[5mm]
\ds \le M(C^{\frac{l+1}{N}}\res N^i{\la\rho_i}\cap N^j_{\la\rho_j}) \ \frac{r_i^2}{N}\quad \le K\ \frac{r_i^4}{N}.
\end{array}
\ee
 Combining (\ref{VII.47}), (\ref{VII.48}), (\ref{VII.49}), (\ref{VII.50}), (\ref{VII.53}), (\ref{VII.55})
and (\ref{VII.56}), we obtain that
\be
\label{VII.57}
\ds M(\Txi_\ast[\Pi\Psi((\Sigma')^{ij})]\times[0,1]\res N^i_{\la\rho_i}\cap N^j_{\la\rho_j})\le r_i^4\quad .
\ee
Combining this last inequality with (\ref{VII.44}) and a Fubini type argument we obtain that 
there exists $\la\in[\frac{7}{6},\frac{6}{5}]$ such that
\be
\label{VII.58}
{\mathcal F}((\ti{C}_i-\ti{C}_j)\res N^i{\la\rho_i}\cap N^j_{\la\rho_j})\le r_i^4\quad.
\ee
From this fact we then deduce, since $\ti{C}_i$ and $\ti{C}_j$ are single valued graphs with uniformly bounded 
gradients
\be
\label{VII.59}
\int_{\Pi_i(N^i{\frac{7}{6}\rho_i}\cap N^j_{\frac{7}{6}\rho_j})}|\ti{a}^i_j(z^i)-\ti{a}^j_j|\le K r_i^4\quad.
\ee

 Using the notations introduced in (\ref{VII.3a}) and (\ref{VII.27}), we have that for any $z$ there exists $\xi$ such that
\be
\label{VII.60}
(z-z_i,\ti{a}_i(z)-w_i)=\La_{x_i}^{-1}(\xi,\ti{a}_i^i(\xi))\quad,
\ee
where $|\La_{x_i}-id|\le K\ |x_i|\le K\ \rho_{x_0}$. Let $z':= \rho_i^{-1}(z-z_i)$ and $\hat{a}_i(z'):=\rho_i^{-1}(\ti{a}_i(z)-w_i)$. Let also
$\xi':=\rho_i^{-1}\xi$ and $\hat{a}_i^i(\xi'):=\ti{a}_i^i(\xi)$. Since $\ti{a}_i^i$ is holomorphic (see (\ref{VII.28}), $\hat{a}_i^i$ is also clearly
holomorphic and since $\|\hat{a}_i^i\|_{L^\infty(B_\frac{3}{2}(0))}\le K$, we have that for any $l\in {\N}$
\be
\label{VII.60a}
\|\nabla^l \hat{a}_i^i\|_{L^\infty(B_\frac{5}{4}(0))}\le K_l
\ee
Using the above notations we have 
\be
\label{VII.60b}
(z',\hat{a}_i(z'))=\La_{x_i}^{-1}(\xi',\hat{a}_i^i(\xi'))\quad .
\ee
From the local inversion theorem, since $|\La_{x_i}-id|\le K\ |x_i|\le K\ \rho_{x_0}$ can be taken as small as we want by taking $\rho_{x_0}$
small enough, we have that for all $l\in {\N}$ there exists $K_l$ such that
\be
\label{VII.60c}
\|\nabla_{z'}^l\xi'\|_\infty\le K_l\quad .
\ee
Therefore, combining (\ref{VII.60a}), (\ref{VII.60b} and (\ref{VII.60c}, we obtain
\be
\label{VII.60d}
\|\nabla_{z'}^l\hat{a}_i(z')\|_\infty\le K_l\quad .
\ee
From that estimate we then deduce
\be
\label{VII.60e}
\|\nabla^l_z\ti{a}_i\|_{L^\infty(B^2_{\frac{7}{6}\rho_i}(z_i))}\le K_l\ r_i^{l-1}\quad .
\ee
Since $\ti{C}^i$ is $J_{x_i}-$holomorphic, we have the existence of $\la_1^i,\mu_1^i,\la_2^i,\mu_2^i$ such that
\be
\label{VII.61}
\lf\{
\begin{array}{l}
\ds J_{x_i}\cdot\lf(
\begin{array}{c}
1\\[3mm]
0\\[3mm]
\ds \frac{\p\ti{a}_i}{\p x}
\end{array}
\rg)=\la_1^i
\lf(
\begin{array}{c}
1\\[3mm]
0\\[3mm]
\ds\frac{\p\ti{a}_i}{\p x}
\end{array}
\rg)+\mu_1^i
\lf(
\begin{array}{c}
0\\[3mm]
1\\[3mm]
\ds\frac{\p\ti{a}_i}{\p y}
\end{array}
\rg)\\[5mm]
\ds J_{x_i}\cdot\lf(
\begin{array}{c}
0\\[3mm]
1\\[3mm]
\ds \frac{\p\ti{a}_i}{\p y}
\end{array}
\rg)=\la_2^i
\lf(
\begin{array}{c}
1\\[3mm]
0\\[3mm]
\ds\frac{\p\ti{a}_i}{\p x}
\end{array}
\rg)+\mu_2^i
\lf(
\begin{array}{c}
0\\[3mm]
1\\[3mm]
\ds\frac{\p\ti{a}_i}{\p y}
\end{array}
\rg)
\end{array}
\rg.
\ee
Writing $J_{x_i}=J_0+\delta(x_i)$, we first observe that 
\be
\label{VII.62}
|\delta(x_i)|\le \|J\|_{C^1}|x_i|\le K\rho_{x_0}\quad.
\ee
Using this notation we deduce from (\ref{VII.61})
\be
\label{VII.63}
\lf\{
\begin{array}{l}
\ds\la_1^i=\delta_{1,1}(x_i)+\sum_{l=3}^{2p}\delta_{1,l}(x_i)\frac{\p \ti{a}_i^l}{\p x}\\[5mm]
\ds \mu_1^i= 1+\delta_{2,1}(x_i)+\sum_{l=3}^{2p}\delta_{2,l}(x_i)\frac{\p \ti{a}_i^l}{\p x}\\[5mm]
\ds \la_2^i=-1+\delta_{1,2}(x_i)+\sum_{l=3}^{2p}\delta_{1,l}(x_i)\frac{\p \ti{a}_i^l}{\p y}\\[5mm]
\ds\mu_2^i=\delta_{2,2}(x_i)+\sum_{l=3}^{2p}\delta_{2,l}(x_i)\frac{\p \ti{a}_i^l}{\p y}
\end{array}
\rg.
\ee
Therefore the equation solved by $\ti{a}_i$ is for any $k=1\cdots p-1$
\be
\label{VII.64}
\lf\{
\begin{array}{rl}
\ds\frac{\p\ti{a}_i^{2k+1}}{\p x}-\frac{\p\ti{a}_i^{2k+2}}{\p y} &\ds=\lf[\delta_{1,1}(x_i)+\sum_{l=3}^{2p}\delta_{1,l}(x_i)\frac{\p \ti{a}_i^l}{\p x}\rg]\
\frac{\p\ti{a}_i^{2k+2}}{\p x}\\[3mm]
 &\ds +\lf[\delta_{2,1}(x_i)
 +\sum_{l=3}^{2p}\delta_{2,l}(x_i)\frac{\p \ti{a}_i^l}{\p x}\rg]\ \frac{\p\ti{a}_i^{2k+2}}{\p y}\\[3mm]
 &\ds-\delta_{2k+2,1}(x_i)-\sum_{l=3}^{2p}\delta_{2k+2,l}\frac{\p \ti{a}_i^{l}}{\p x}\\[5mm]
\ds\frac{\p\ti{a}_i^{2k+1}}{\p y}+\frac{\p\ti{a}_i^{2k+2}}{\p x}&\ds=\lf[\delta_{1,2}(x_i)+\sum_{l=3}^{2p}\delta_{1,l}(x_i)\frac{\p \ti{a}_i^l}{\p y}\rg]
\frac{\p\ti{a}_i^{2k+2}}{\p x}\\[3mm]
 &\ds+\lf[\delta_{2,2}(x_i)+\sum_{l=3}^{2p}\delta_{2,l}(x_i)\frac{\p \ti{a}_i^l}{\p y}\rg]\ \frac{\p \ti{a}_i^{2k+2}}{\p y}\\[3mm]
 &\ds-\delta_{2k+2,1}(x_i)-\sum_{l=3}^{2p}\delta_{2k+2,l}(x_i)\frac{\p \ti{a}_i^{l}}{\p y}
\end{array}
\rg.
\ee
Then we deduce that there exists a linear map
\[
A(x_i,\cdot)\ :\ {\R}^2\otimes{\R}^{2p-2}\longrightarrow{\C}^{p-1}\otimes_{\R}\lf({\R}^2\otimes{\R}^{2p-2}\rg)^\ast\quad ,
\]
and an element
\[
B(x_i,\cdot)\in{\C}^{p-1}\otimes_{\R}\lf({\R}^2\otimes{\R}^{2p-2}\rg)^\ast\quad,
\]
such that $\ti{a}_i$ solves
\be
\label{VII.65}
\frac{\p\ti{a}_i}{\p\ov{z}}=A(x_i,\nabla\ti{a}_i)\cdot\nabla\ti{a}_i+B(x_i,\nabla\ti{a}_i)+D(x_i,\ti{a}_i)\quad.
\ee
Observe also that the dependence of $A$ and $B$ in $B_{\rho_{x_0}}(x_0)$ is smooth and that $A(x_0,\cdot)=0$, $B(x_0,\cdot)=0$, $D(x_0)=0$ and because
of (\ref{VII.62}) we have an estimate
of the sort
\be
\label{VII.66}
\forall p\in{\R}^2\otimes{\R}^{2p}\quad\quad |A(x_i,p)|+|B(x_i,p)|\le K\ |x_i|(1+|p|)\quad.
\ee
Consider now $i$ and $j$ such that $B^2_{\rho_i}(z_i)\cap B^2_{\rho_j}(z_j)\ne\emptyset$. On $B^2_{\frac{7}{6}\rho_i}(z_i)\cap B^2_{\frac{7}{6}\rho_j}(z_j)$
$\ti{a}_i-\ti{a}_j$ solves the following equation
\be
\label{VII.67}
\begin{array}{rl}
\ds \p_{\ov{z}}(\ti{a}_i-\ti{a}_j)&\ds = A(x_i,\nabla\ti{a}_i)\cdot\nabla\ti{a}_i+B(x_i,\nabla\ti{a}_i)\\[5mm]
 &\ds-A(x_j,\nabla\ti{a}_j)\cdot\nabla\ti{a}_j-B(x_j,\nabla\ti{a}_j)\\[5mm]
 &\ds=C(x_i,\nabla\ti{a}_i,\nabla\ti{a}_j)\cdot \nabla(\ti{a}_i-\ti{a}_j)\\[5mm]
 &\ds + E(x_i,\nabla\ti{a}_j)-E(x_j,\nabla\ti{a}_j)\quad,
\end{array}
\ee
where 
\be
\label{VII.68}
\begin{array}{rl}
\ds C(x_i,\nabla\ti{a}_i,\nabla\ti{a}_j)\cdot \nabla(\ti{a}_i-\ti{a}_j)&\ds :=A(x_i,\nabla \ti{a}_i)\cdot\nabla\ti{a}_i-A(x_i,\nabla \ti{a}_j)\cdot
\nabla\ti{a}_j\\[5mm]
 &\ds +B(x_i,\nabla\ti{a}_i)-B(x_i,\nabla\ti{a}_j)\quad,
\end{array}
\ee
(where we have used the linear dependance in p of $A(x_i,p)$ and $B(x_i,p)$), and where
\be
\label{VII.69}
E(x,p):= A(x,p)\cdot p+B(x,p)+D(x)\quad.
\ee
Observe, in one hand, that $C(x,p,q)$ has a linear dependence in $p$ and $q$ in ${\R}^2\otimes{\R}^{2p-2}$, that
\be
\label{VII.70}
|C(x,p,q)|\le K\ |x|\ (1+|p|+|q|)\quad .
\ee
and that following estimates hold for $D(x,p)$, forall $l\in {\N}$
\be
\label{VII.71}
|\nabla_x^lE(x,p)|\le K_l\ (1+|p|^2)\quad .
\ee
On $B^2_{\frac{7}{6}}(\rho_i^{-1}z_i)\cap B^2_{\frac{7}{6}}(\rho_i^{-1}z_j)$ the function
$f(z'):=a_i(\rho_i z')-a_j(\rho_i z')$ solves
\be
\label{VII.72}
\p_{\ov{z'}}f-C(z')\cdot\nabla f= g(z')\quad,
\ee
where 
\[
C(z'):=C(x_i,\nabla\ti{a}_i,\nabla\ti{a}_j)(\rho_i z')\quad, 
\]
and
\[
g(z'):= \rho_i\lf[D(x_i,\nabla \ti{a}_j(\rho_iz'))-D(x_j,\nabla\ti{a}_j(\rho_iz'))\rg]\quad .
\]
Using (\ref{VII.60e}),  observe that for any $l\in{\N}$
\be
\label{VII.73}
\begin{array}{l}
\ds\|\nabla^l_{z'}(C(x_i,(\nabla_z \ti{a}_i)(\rho_i z'),(\nabla_z \ti{a}_j)(\rho_i z')))\|_\infty\\[5mm]
\ds\le 
K\ |x_i|\ \rho_i^l\lf[\|\nabla_z^{l+1}\ti{a}_i\|_\infty+\|\nabla_z^{l+1}\ti{a}_j\|_\infty\rg]
\le \ K_l\ \rho_{x_0}\quad.
\end{array}
\ee
Therefore, for $\rho_{x_0}$ small enough,
 $L:=\p_{\ov{z'}}-C(z')\cdot\nabla_{z'}$ is an elliptic coercive first order operator with smooth coefficients
whose derivatives are uniformly bounded.
Observe also that, using again (\ref{VII.60e}),
\be
\label{VII.74}
\begin{array}{l}
\ds\|\nabla^l_{z'}\rho_i\lf[D(x_i,\nabla \ti{a}_j(\rho_iz'))-D(x_j,\nabla\ti{a}_j(\rho_iz'))\rg]\|_\infty
\\[5mm]
\ds\le K\rho_i^2\ \lf[\rho_i^l\sum_{s=0}^{[l/2]}\|\nabla_z^{s+1}\ti{a}_j\|_\infty
\ \|\nabla_z^{l-s+1}\ti{a}_j\|_\infty+\rho_i^l\ \|\nabla_{z'}^{l+1}\ti{a}_j\|_\infty\rg]
\end{array}
\ee
Then we have
\be
\label{VII.75}
\|\nabla^l_{z'}g\|_\infty\le K_l\ \rho_i^2\quad.
\ee
From (\ref{VII.58}) we deduce that
\be
\label{VII.76}
\int_{B^2_{\frac{7}{6}}(\rho_i^{-1}z_i)\cap B^2_{\frac{7}{6}}(\rho_i^{-1}z_j)}|f|\le K\ \rho_i^2\quad.
\ee
Thus combining (\ref{VII.72})...(\ref{VII.76}) and using standard elliptic estimates we obtain that 
for any $l\in {\N}$
\be
\label{VII.77}
\|\nabla^l f\|_{L^\infty(B^2_{1}(\rho_i^{-1}z_i)\cap B^2_{1}(\rho_i^{-1}z_j))}\le K_l\ \rho_i\quad,
\ee
which yields, going back to the original scale the estimate (\ref{VII.78}) and lemma~\ref{lm-VII1a} is proved.\cqfd

{\bf Definition of the approximated average curve.} On $B^2_{\rho_{x_0}}(0)=\Pi(B_{\rho_{x_0}}(x_0))$ we define
the approximated average curve as follows. Let $\vphi$ be the partition of unity of $\Pi({\mathcal C}_{Q-1}\cap B_{\rho_{x_0}}(x_0))$
defined in (\ref{V.ax14}). We denote 
\be
\label{VII.79}
\lf\{
\begin{array}{l}
\ds\ti{a}(z_0):=\sum_{i\in I}\vphi(z_0)\ti{a}_i(z_0)\quad\quad\forall z_0\in\Pi({\mathcal C}_{Q-1}\cap B_{\rho_{x_0}}(x_0))\quad .\\[5mm]
\ds\ti{a}(z_0):=<C,z,z_0>(w)\quad\quad\forall z_0\in\Pi(({\mathcal C}_Q\setminus{\mathcal C}_{Q-1})\cap B_{\rho_{x_0}}(x_0))
\end{array}
\rg.
\ee
Observe that, because of lemma~\ref{lm-IV.2}, for any $z_0\in\Pi(({\mathcal C}_Q\setminus{\mathcal C}_{Q-1})\cap B_{\rho_{x_0}}(x_0))$
the slice $<C,z,z_0>$ consists of exactly one point and $\ti{a}(z_0)$ is simply the $w$ coordinates of that point.
The following estimates for $\ti{a}$ holds :
\begin{Lm}
\label{lm-VII.2}
Under the above notations, for any $q<+\infty$ there exists a constant $K_q$ independent of $\rho_{x_0}$ such that
\be
\label{VII.96b}
\int_{B^2_{\frac{\rho_{x_0}}{2}}(0)}|\nabla^2\ti{a}|^q\le K_q\ \rho_{x_0}^2\quad.
\ee 
\end{Lm}
{\bf Proof of lemma~\ref{lm-VII.2}.}
 
We claim first that $\ti{a}$ is a Lipschitz map over $B^2_{\rho_{x_0}}(0)$. In $\Pi({\mathcal C}_{Q-1}\cap B_{\rho_{x_0}}(x_0))$,
by the assumptions of the inductive procedure, we know that $\ti{a}$ is smooth and using both (\ref{VII.60e}) and (\ref{VII.78}),
because also of (\ref{V.ax13}), we have
\be
\label{VII.80}
\|\nabla\ti{a}\|_{L^\infty(\Pi({\mathcal C}_{Q-1}\cap B_{\rho_{x_0}}(x_0)))}\le K\quad .
\ee
Consider now two arbitrary points $x$ and $y$ of $B^2_{\rho_{x_0}}(0)$, either the segment $[x,y]$ in $B^2_{\rho_{x_0}}(0)$
is included in $\Pi({\mathcal C}_{Q-1}\cap B_{\rho_{x_0}}(x_0))$ and then we can integrate (\ref{VII.80}) all along that
segment to get
\be
\label{VII.81}
|\ti{a}(x)-\ti{a}(y)|\le K|x-y|\quad,
\ee
or there exists $z\in [x,y]\cap\Pi(({\mathcal C}_Q\setminus{\mathcal C}_{Q-1})\cap B_{\rho_{x_0}}(x_0))$ and using this time
(\ref{IV.1b}) we also get (\ref{VII.81}), which proves the desired claim.
Using the equation (\ref{VII.65}) solved by the $\ti{a}_i$s we obtain that, in 
$\Pi({\mathcal C}_{Q-1}\cap B_{\rho_{x_0}}(x_0))$, $\ti{a}$ is a solution of 
\be
\label{VII.82}
\p_{\ov{z}}\ti{a}-A((z,\ti{a}(z)),\nabla\ti{a})\cdot\nabla\ti{a}-B((z,\ti{a}(z)),\nabla\ti{a})-D(z,\ti{a}(z))=\zeta(z)\quad ,
\ee  
where
\be
\label{VII.83}
\begin{array}{rl}
\ds\zeta(z)&\ds:=\sum_{i\in I}\vphi_i\ \lf[A(x_i,\nabla\ti{a}_i)\cdot\nabla\ti{a}_i
-A((z,\ti{a}(z)),\nabla\ti{a})\cdot\nabla\ti{a}\rg]\\[5mm]
 &\ds+\sum_{i\in I}\ \lf[B(x_i,\nabla\ti{a}_i)-B((z,\ti{a}(z)),\nabla\ti{a})\rg]\\[5mm]
 &\ds+\sum_{i\in I}\p_{\ov{z}}\vphi_i\ \lf[\ti{a}_i-\ti{a}\rg]\quad .
\end{array}
\ee
Observe that $\sum_{i\in I}\p_{\ov{z}}\vphi_i\ti{a}$ was added at the end because this quantity vanishes.
Since for any $x_i$ and $r_i$, (\ref{V.aa3}) holds, since also for any 
$z\in\Pi(({\mathcal C}_Q\setminus{\mathcal C}_{Q-1})\cap B_{\rho_{x_0}}(x_0))$ we have for an $\ep$
as small as we want (recall $\al$ was fixed independently of $\ep$), the relative lipschitz estimate
(\ref{IV.1b}) holds and granting the fact that 
$\Pi(({\mathcal C}_Q\setminus{\mathcal C}_{Q-1})\cap B_{\rho_{x_0}}(x_0))$ is a compact subset of 
$B_{\rho_{x_0}}(x_0)$, it is clear that
\be
\label{VII.84}
\begin{array}{l}
\ds\forall\eta>0\quad\exists\delta>0\quad\quad\mbox{ s. t. }\forall i\in I\\[5mm]
\quad\ds dist(B_{\rho_i}(z_i),\Pi(({\mathcal C}_Q\setminus{\mathcal C}_{Q-1})))\le\delta\quad
\Longrightarrow\quad |\rho_i|\le\eta\quad .
\end{array}
\ee
Observe that 
\be
\label{VII.85}
\begin{array}{rl}
\ds|\zeta(z)| &\ds\le K\sum_{i,j\in I}\vphi_i(z)\vphi_j(z)\ |\ti{a}_j-w_i|\\[5mm]
 &\ds\sum_{i,j\in I}\ \vphi_i(z)\vphi_j(z)\ |\nabla(\ti{a}_i-\ti{a}_j)|\\[5mm]
 &\ds\sum_{i,j\in I}\ |\nabla\vphi_i(z)|\vphi_j(z)\ |\ti{a}_i(z)-\ti{a}_j(z)|\quad,
\end{array}
\ee
using (\ref{VII.78}), we have that
\be
\label{VII.86}
|\zeta(z)|\le\ K \max\lf\{r_i\ i\in I;\ \mbox{ s. t. }z\in B^2_{\rho_i}(z_i)\rg\}\le K\ 
dist(z,\Pi(({\mathcal C}_Q\setminus{\mathcal C}_{Q-1})))\quad.
\ee
Combining (\ref{VII.84}) and (\ref{VII.86}), we have that $\zeta(z)$ converge uniformly to 0 as $z$
tends to $\Pi(({\mathcal C}_Q\setminus{\mathcal C}_{Q-1}))$. We then extend $\zeta$ by 0 in 
$\Pi(({\mathcal C}_Q\setminus{\mathcal C}_{Q-1}))$. $\zeta$ is now a continuous function in 
$B^2_{\rho_{x_0}}(x_0)$ and we claim that
\be
\label{VII.87}
\p_{\ov{z}}\ti{a}-A((z,\ti{a}(z)),\nabla\ti{a})\cdot\nabla\ti{a}-B((z,\ti{a}(z)),\nabla\ti{a})=\zeta(z)
\quad\quad\mbox{ in }{\mathcal D}'(B^2_{\rho_{x_0}}(x_0))\quad.
\ee
Since $\ti{a}$ is a lipschitz function  in $B^2_{\rho_{x_0}}(x_0)$, $\p_{\ov{z}}\ti{a}-A((z,\ti{a}(z)),\nabla\ti{a})\cdot\nabla\ti{a}-B((z,\ti{a}(z)),\nabla\ti{a})$ is a bounded function in $B^2_{\rho_{x_0}}(x_0)$ and 
therefore, in order to prove (\ref{VII.87}), it suffices to prove that for ${\mathcal H}^2$ almost every $z$ in
$\Pi(({\mathcal C}_Q\setminus{\mathcal C}_{Q-1}))\cap B^2_{\rho_{x_0}}(x_0)$ 
\be
\label{VII.88}
\p_{\ov{z}}\ti{a}-A((z,\ti{a}(z)),\nabla\ti{a})\cdot\nabla\ti{a}-B((z,\ti{a}(z)),\nabla\ti{a})-D(z,\ti{a}(z))=0\quad .
\ee
This later equality, because of the computations leaded between (\ref{VII.61})...(\ref{VII.65}), is just equivalent
to the fact that the tangent plane of the graph $(z,\ti{a}(z))$ at that point is $J-$holomorphic, which is the case
at every point of $\Pi(({\mathcal C}_Q\setminus{\mathcal C}_{Q-1}))\cap B^2_{\rho_{x_0}}(x_0)$ due to 
lemma~\ref{lm-IV.2}. Thus (\ref{VII.87}) is proved.
Using again (\ref{VII.78}) we observe that
\be
\label{VII.89}
\begin{array}{rl}
\ds|\nabla\zeta(z)| &\ds\le K\sum_{i,j\in I}\vphi_i(z)\vphi_j(z)\ \lf[|\nabla\ti{a}_j|+|\nabla\ti{a}_j|^3\rg]\\[5mm]
 &\ds\sum_{i,j\in I}\ \vphi_i(z)\vphi_j(z)\ |\nabla\ti{a}_i|(z)\ |\nabla^2(\ti{a}_i-\ti{a}_j)|(z)\\[5mm]
 &\ds\sum_{i,j\in I}\ \lf[|\nabla^2\vphi_i(z)|\vphi_j(z)+|\nabla\vphi_i|(z)|\nabla\vphi_j(z)|\rg]\ |\ti{a}_i(z)-\ti{a}_j(z)|+|\ti{a}_j-w_i|\quad,\\[5mm]
 &\le K
\end{array}
\ee
We claim now that $\zeta$ is Lipschitz in $B^2_{\rho_{x_0}}(x_0)$. Indeed, argueing like for $\ti{a}$, given $x$ and $y$ in $B^2_{\rho_{x_0}}(x_0)$,
if the segment $[x,y]$ has no intersection with $\Pi(({\mathcal C}_Q\setminus{\mathcal C}_{Q-1}))$, then, integrating (\ref{VII.85}) on that segment
gives
\be
\label{VII.90}
|\zeta(x)-\zeta(y)|\le K\ |x-y|\quad.
\ee
Otherwise, if  there exist $z\in [x,y]\cap\Pi(({\mathcal C}_Q\setminus{\mathcal C}_{Q-1}))$, then, (\ref{VII.86}) gives
\[
|\zeta(x)|+|\zeta(y)|\le K\ |x-z|+|y-z|\quad,
\]
which gives (\ref{VII.90}) and the claim is proved. We claim now that $\ti{a}\in W^{2,q}(B^2_{\rho_{x_0}/2}(0))$
for any $q<+\infty$. Let $e$ be a unit vector  in 
$B^2_{\rho_{x_0}/2}(0)$. For small $h$ we denote $\ti{a}_h(z):=\ti{a}(z+h)$ and $\zeta_h(z):=\zeta(z+h)$. We have, using the linear dependencies
of 
\be
\label{VII.91}
\begin{array}{rl}
\ds\p_{\ov{z}}(\ti{a}-\ti{a}_h)&\ds -A((z,\ti{a}),\nabla\ti{a})\cdot\nabla(\ti{a}-\ti{a}_h)-A((z,\ti{a}),\nabla(\ti{a}-\ti{a}_h))\cdot\nabla\ti{a}_h\\[5mm]
 &\ds-B((z,\ti{a}),\nabla(\ti{a}-\ti{a}_h))\\[5mm]
 &\ds=\zeta-\zeta_h +A((z,\ti{a}),\nabla\ti{a}_h)\cdot\nabla\ti{a}_h-A((z+h,\ti{a}_h),\nabla\ti{a}_h)\cdot\nabla\ti{a}_h\\[5mm]
 &\ds B((z,\ti{a}),\nabla\ti{a}_h)-B((z+h,\ti{a}_h),\nabla\ti{a}_h)+D(z,\ti{a})-D(z+h,\ti{a}_h)\quad.
\end{array}
\ee
Denote $\Delta_h$ the right-hand side of (\ref{VII.91}) and observe that there exists a constant $K$ such that
\be
\label{VII.92}
|\Delta_h|\le K\ h\quad.
\ee
Let $\chi_{\rho_{x_0}}$ be a cut-off function  such that $\chi_{\rho_{x_0}}\equiv 1$ in $B^2_{\frac{\rho_{x_0}}{2}}(0)$ and $\chi_{\rho_{x_0}}\equiv 0$
in ${\R}^2\setminus B^2_{{\rho_{x_0}}}(0)$. Let $f_h:=\chi_{\rho_{x_0}}\ (\ti{a}-\ti{a}_h)$, and let $L_h$ be the operator such that
\be
\label{VII.93}
L_hf:=\p_{\ov{z}}f-A((z,\ti{a}),\nabla\ti{a})\cdot\nabla f-A((z,\ti{a}),\nabla f)\cdot\nabla\ti{a}_h-B((z,\ti{a}),\nabla f)\quad,
\ee
we have
\be
\label{VII.94}
\begin{array}{rl}
\ds L_hf_h &\ds=\chi_{\rho_{x_0}}\Delta_h+(\ti{a}-\ti{a}_h)\p_{\ov{z}}\chi_{\rho_{x_0}}
-A((z,\ti{a}),\nabla\ti{a})\cdot(\ti{a}-\ti{a}_h)\nabla\chi_{\rho_{x_0}}\\[5mm]
 &\ds-A((z,\ti{a}),(\ti{a}-\ti{a}_h)\nabla\chi_{\rho_{x_0}})\cdot\nabla\ti{a}_h-B((z,\ti{a}),(\ti{a}-\ti{a}_h)\nabla\chi_{\rho_{x_0}})\quad .
\end{array}
\ee
Since $\ti{a}$ is Lipschitz, using (\ref{VII.92}) and (\ref{VII.66}), we have that
\be
\label{VII.95}
|L_hf_h|\le K\ h\quad .
\ee
Observe that $|L_hf_h|\ge |\p_{\ov{z}}f_h|-K\rho_{x_0}|\nabla f_h|$. Since $f_h=0$ on $\p B^2_{\rho_{x_0}}(0)$,
for any $p<+\infty$ we have that
\be
\label{VII.96}
\int_{B^2_{\rho_{x_0}}(0)}|\nabla f_h|^q\le K_q\int_{B^2_{\rho_{x_0}}(0)}|\p_{\ov{z}}f_h|^q\le 
K_q\int_{B^2_{\rho_{x_0}}(0)}|L_h f_h|^q
+K\rho_{x_0}^q|\nabla f_h|^q
\quad .
\ee
Dividing by $h^q$ and making $h$ tend to zero, we get that for $\rho_{x_0}$ small enough we get
(\ref{VII.96b}) and lemma~\ref{lm-VII.2} is proved.\cqfd

\subsection{Constructing adapted coordinates to $C$ in a neighborhood of $x_0\in{\mathcal C}_Q\setminus
{\mathcal C}_{Q-1}$.}

We consider a point $x_0$ in the support of our $J-$holomorphic current $C$ and  we assume,
as above,  that the multiplicity at $x_0$ is $Q$ and that the tangent cone
$C_{0,x_0}$ is $Q$ times a $J_{x_0}-$holomorphic disk $D$. We start with the coordinates $(z,w_1,\cdots,w_{p-1})$
chosed in (\ref{II.0}) such that $C_{0,x_0}=Q[D]$ is $Q$ times the ``horizontal'' disk given by $w_i=0$
for $i=1,\cdots,p-1$ and we work in the ball $B^{2p}_{\rho_{x_0}}(x_0)$ whose radius $\rho_{x_0}$ is given by
(\ref{V.a11}). The purpose of this subsection is to construct new coordinates $(\xi,\la_1,\cdots,\la_{p-1})$
 in $B^{2p}_{\rho_{x_0}}(x_0)$ such that the set $\la_i=0$ for $i=1,\cdots,p-1$ coincides with the graph
of the average map $\ti{a}$ constructed in the previous subsection. 

On the graph $\ti{A}(z):=(z,\ti{a}(z))$ we consider the complex structure $j$ given by the metric induced by 
$g:=\om(J\cdot,\cdot)$. Let $X$ be a vector tangent to $\ti{A}$ at $(z,\ti{a}(z))$. We compare
$jX$ and $JX$ in ${\R}^{2p}$. Let $n(z,\ti{a}(z))$ be the $2p-2$-unit vector normal to 
$T_{(z,\ti{a}(z))}\ti{A}$, making the identification between $2p-1$-vectors and vector given by the
ambiant metric, we have
\[
jX:=n\wedge X\quad .
\]
We have seen that
\[
|Jn-n|(z,\ti{a}(z))\le r_{(z,\ti{a}(z))}\quad .
\]
Therefore
\[
|jJX-JjX|\le|n\wedge JX-J(n\wedge X)|\le |n\wedge JX-Jn\wedge JX|\le r_{(z,\ti{a}(z))}|X|\quad.
\]
Thus we get that $|(J-j)(J+j)X|\le r_{(z,\ti{a}(z))}|X|$ where we extend $j$ to the normal bundle
to $\ti{A}$ again by the mean of the induced metric. Since $|(J+j)X|$ and $|X|$ are comparable independant of $X$, we have
\be
\label{VIII.0}
\forall X\in T_{(z,\ti{a}(z))}\ti{A}\quad\quad|(J-j)X|\le K  r_{(z,\ti{a}(z))}|X|\quad.
\ee
We choose now coordinates $\xi=(\xi_1,\xi_2)$ on $\ti{A}$ compatible with $j$ (i.e. 
$j\frac{\p}{\p\xi_1}=\frac{\p}{\p\xi_2}$).
Let $(z',\hat{a}(z')):=(\rho_{x_0}^{-1},\rho_{x_0}^{-1}\ti{a}(\rho_{x_0}z'))$ and in $B^{2p}_1(0)$ we consider the
 metric $\hat{g}(z',w'):=\rho_{x_0}^{-2}(\rho_{x_0}z',\rho_{x_0}w')^\ast g$ where $g$ in $B^{2p}_{\rho_{x_0}}(0)$
is the original metric $g(\cdot,\cdot)=\om(J\cdot,\cdot)$. After this scaling we have, using (\ref{VII.96b})
\be
\label{VIII.01}
\int_{B^2_1}|\nabla^2\hat{a}|^p\ dz'\le K_q\rho_{x_0}^q\quad ,
\ee
and
\be
\label{VIII.02}
\hat{g}^{ij}=\delta^{ij}+h^{ij}\quad\mbox{ where }h^{ij}(0,0)=0\quad\mbox{ and }\|\nabla h^{ij}\|_\infty\le K\rho_{x_0}
\quad .
\ee
We look for isothermal coordinates $(\xi'_1,\xi'_2)$ in $\hat{A}=\{(z',\hat{a}(z'))\quad z'\in B^2_1(0)\}$ of the
 form $\xi'=z'+\delta(z')$ where $\delta$ will be small in $W^{2,p}$. On $B^2_1(0)$ we consider the metric
$\hat{k}=(z',\hat{a}(z'))^\ast \hat{g}=(1+k_{11})(dx_1')^2+2k_{12}dx_1'dx_2'+(1+k_{22})(dx_2')^2$. 
From the estimates above we have for any $q>0$, (since $\nabla\hat{a}(0,0)=0$ and 
$\|\nabla^2\hat{a}\|_q\le \rho_{x_0}$ we have $\|\nabla\hat{a}\|_\infty\le\rho_{x_0}$)
\be
\label{VIII.02b}
\int_{B_1^2}|\nabla k|^q\le K_q \rho_{x_0}^{2q}\quad .
\ee 
Following \cite{DNF} page 110-111
it suffices to find $\delta_1$ solving 
\be
\label{VIII.03}
-\frac{\p}{\p x_1'}\lf[
\frac{(1+k_{11})\frac{\p\delta_1}{\p x'_1}-k_{12}\frac{\p\delta_1}{\p x'_2}}{\sqrt{(1+k_{11})(1+k_{22})-k_{12}^2}} 
\rg]-\frac{\p}{\p x_2'}\lf[
\frac{(1+k_{22})\frac{\p\delta_1}{\p x'_2}-k_{12}\frac{\p\delta_1}{\p x'_1}}{\sqrt{(1+k_{11})(1+k_{22})-k_{12}^2}} 
\rg]=0\quad.
\ee
Taking $\delta_1=0$ on $\p B_1^2(0)$ we get a well posed elliptic problem and we obtain the existence of
$\delta_1$ satisfying
\[
\sum_{i=1}^2a_{ij}\frac{\p^2\delta_1}{\p x'_i\p x'_j}=F\cdot\nabla\delta_1\quad,
\]
where $a_{ij}$ are H\"older continuous $\|a_{ij}-\delta_{ij}\|_{C^{0,\al}(B_1^2)}\le K_\al\rho^2$ and $F\in L^q$ with $\int|F|^q\le K_p \rho_{x_0}^{2q}$.
Standard elliptic estimates give then
\be
\label{VIII.04}
\|\delta_1\|_{W^{2,q}(B_1^2)}\le K\rho_{x_0}^2\quad.
\ee
Therefore, going back to the original scale, we have found coordinates 
$\xi_i=x_i+\rho_{x_0}\delta_i(\rho^{-1}_{x_0}z)=x_i+\al_i(z)$ such that 
\be
\label{VIII.05}
\|\nabla\al\|_\infty\le K \rho_{x_0}^2\quad\mbox{ and }\quad j\frac{\p}{\p\xi_1}=\frac{\p}{\p\xi_2}\quad.
\ee
We translate these coordinates in such a way that $\al(0,0)=(0,0)$.

Inside $Gl({\R}^{2p})$, the space of invertible $2p\times 2p$ matrices with real coefficients we denote
$U(p)$ the subspace of matrices $M$ which commute with $J_0$. $U(p)$ is a compact submanifold of $Gl({\R}^{2p})$
and for some metric in $Gl({\R}^{2p})$ we denote $\pi_{U(p)}$ the orthogonal projection from a neighborhood of
$U(p)$ onto $U(p)$. We consider $M(z)$ the matrix which is given by 
\be
\label{VIII.1}
M(z):=\La^{-1}_{(z,\ti{a}(z))}\pi_{U(p)}(\La_{(z,\ti{a}(z))})\quad ,
\ee
where we recall that $\La_x$ is given by (\ref{VII.3c}). We have clearly, since $\|\nabla\ti{a}\|\le K$ and 
$\int_{B^2_{\rho_{x_0}}(0)}|\nabla^2\ti{a}|^q\le K_q\rho_{x_0}^2$ for any $p<+\infty$,
\be
\label{VIII.2}
\|\nabla M(z)\|_{L^\infty(B^2_{\rho_{x_0}}(0))}\le K\quad\mbox{ and }\quad\|\nabla^2M(z)\|_{L^{q}(B^2_{\rho_{x_0}}(0))}\le K_q\rho_{x_0}^{\frac{2}{q}}\quad .
\ee
We keep denoting $e_1$, $e_2,\cdots, e_{2p}$ the canonical
basis of ${\R}^{2p}$. Let
\be
\label{VIII.3}
\varepsilon_k(z):=M(z)\cdot e_i\quad.
\ee
We have forall $i=1\cdots p$ that
\be
\label{VIII.4}
\begin{array}{l}
\ds J((z,\ti{a}(z))\cdot \varepsilon_{2i-1}(z)
=J((z,\ti{a}(z))\cdot\La^{-1}_{(z,\ti{a}(z))}\pi_{U(p)}(\La_{(z,\ti{a}(z))})
\cdot e_{2i-1}\\[5mm]
\ds=\La_{(z,\ti{a}(z))}^{-1}J_0\pi_{U(p)}(\La_{(z,\ti{a}(z))})\cdot e_{2i-1} \\[5mm]
\ds=\La_{(z,\ti{a}(z))}^{-1}\pi_{U(p)}(\La_{(z,\ti{a}(z))})J_0\cdot e_{2i-1}=
\La_{(z,\ti{a}(z))}^{-1}\pi_{U(p)}(\La_{(z,\ti{a}(z))})\cdot e_{2i}=\varepsilon_{2i}(z)\quad .
\end{array}
\ee
In $B^{2p}_{\rho_{x_0}}(x_0)$ we consider the new coordinates $(\xi,\la)$ given by 
\be
\label{VIII.5}
\Psi\ : (\xi,\la)\longrightarrow \Psi(\xi,\la):=
(z(\xi),\ti{a}(z(\xi)))+\sum_{l=1}^{2p}\la_l\varepsilon_{l+2}(z(\xi))\quad .
\ee
Let $\ti{J}$ be the expression of the almost complex structure in these coordinates
(i.e. $\ti{J}_{(\xi,\la)}\cdot X:=d\Psi^{-1} J_{\Psi(\xi,\la)}\cdot d\Psi\cdot X$), we shall now estimate 
$|\ti{J}-J_0|$
for points satisfying $|\la|\le r_{\Psi(\xi,\la)}$ - recall that $r_x$ was defined in the beginning of 
chapter VI (see (\ref{VI.1c}) and (\ref{VI.1d})) - 
which corresponds in $B^{2p}_{\rho_{x_0}}(x_0)$ to a neighborhood
of $\ti{A}(z)=(z,\ti{a}(z))$ containing the support of $C$.
We have first for $i=1,2$, using (\ref{VIII.0}), 
\be
\label{VIII.6}
\begin{array}{l}
d\Psi\ti{J}_{(\xi,0)}e_i=J_{\Psi(\xi,0)}\cdot\frac{\p}{\p\xi_i}=j_{\Psi(\xi,0)}\cdot\frac{\p}{\p\xi_i}
+(J_{\Psi(\xi,0)}-j_{\Psi(\xi,0)})\cdot\frac{\p}{\p\xi_i}\\[5mm]
\ds\quad\quad=(-1)^{i+1}\frac{\p}{\p\xi_{i+1}} + O(r_\Psi(\xi,0))\quad,
\end{array}
\ee
where we are using the convention $\frac{\p}{\p\xi_{i+1}}=\frac{\p}{\p\xi_{i-1}}$. We have then for $1<l\le {p}$
\be
\label{VIII.7}
d\Psi\ti{J}_{(\xi,0)}e_{2l}=J_{\Psi(\xi,0)}\frac{\p}{\p\la_{2l}}=J_{\Psi(\xi,0)}\varepsilon_{2l}=-\varepsilon_{2l-1}
=\frac{\p}{\p \xi_{2l+1}}\quad .
\ee
For $|\la|\le r_{\Psi(\xi,\la)}$ we have for $i=1,2$
\be
\label{VIII.8}
\begin{array}{l}
\ds d\Psi\ti{J}_{(\xi,\la)}e_i=J_{\Psi(\xi,\la)}\cdot d\Psi_{(\xi,\la)}e_i=J_{\Psi(\xi,0)}
\cdot d\Psi_{(\xi,0)}e_i\\[5mm]
\ds +J_{\Psi(\xi,0)}\cdot\lf[d\Psi_{(\xi,\la)}e_i-d\Psi_{(\xi,0)}e_i\rg]
+\lf[J_{\Psi(\xi,\la)}-J_{\Psi(\xi,0)}\rg]\cdot d\Psi_{(\xi,\la)}e_i\quad .
\end{array}
\ee
Using the fact that $|J_{\Psi(\xi,\la)}-J_{\Psi(\xi,0)}|\le \|J\|_{C^1}\ |\Psi(\xi,\la)-\Psi(\xi,0)|\le K
r_{\Psi(\xi,\la)}$ and that $d\Psi_{(\xi,\la)}e_i-d\Psi_{(\xi,0)}e_i=\frac{\p\Psi}{\p\xi_i}(\xi,\la)-
\frac{\p\Psi}{\p\xi_i}(\xi,0)=\sum_{l=1}\la_l\p_{\xi_i}\varepsilon_{l+2}$, we have then
\be
\label{VIII.9}
|d\Psi\ti{J}_{(\xi,\la)}e_1-J_{\Psi(\xi,0)}
\cdot d\Psi_{(\xi,0)}e_1|\le O(r_{\Psi(\xi,\la)})
\ee
Using (\ref{VIII.0}) again, we have $|[J_{\Psi(\xi,0)}-j(\Psi(\xi,0))]
\cdot d\Psi_{(\xi,0)}e_i|\le O(r_{\Psi(\xi,0)})$, thus, since $\Psi$ is lipschitz
$r_{\Psi(\xi,\la)}$ and $r_{\Psi(\xi,0)}$ from lemma~\ref{lm-V.2} are comparable and
(\ref{VIII.9}) implies
\be
\label{VIII.10}
 |d\Psi\ti{J}_{(\xi,\la)}e_1- d\Psi_{(\xi,\la)}\cdot e_2|\le
|d\Psi\ti{J}_{(\xi,\la)}e_1- d\Psi_{(\xi,0)}\cdot e_2|+|d\Psi_{(\xi,0)}\cdot e_2-d\Psi_{(\xi,\la)}\cdot e_2|
\ee
and using again the fact that $|d\Psi_{(\xi,0)}\cdot e_2-d\Psi_{(\xi,\la)}\cdot e_2|=|\frac{\p\Psi}{\p\xi_2}(\xi,\la)-
\frac{\p\Psi}{\p\xi_2}(\xi,0)|=|\sum_{l=1}\la_l\p_{\xi_2}\varepsilon_{l+2}|\le O(r_{\Psi(\xi,\la)})$,
we have finally that
\be
\label{VIII.11}
|d\Psi\ti{J}_{(\xi,\la)}e_1- d\Psi_{(\xi,\la)}\cdot e_2|\le O(r_{\Psi(\xi,\la)})\quad .
\ee
Finally, we have for $1<l\le p$
\be
\label{VIII.12}
\begin{array}{l}
\ds d\Psi\ti{J}_{(\xi,\la)}e_{2l}=J_{\Psi(\xi,\la)}\cdot d\Psi_{(\xi,\la)} e_{2l}= -d\Psi_{(\xi,\la)} e_{2l-1}
\\[5mm]
\ds\quad+[d\Psi_{(\xi,\la)}-d\Psi_{(\xi,0)}]e_{2l-1} + J_{\Psi(\xi,0)}\cdot 
[d\Psi_{(\xi,\la)}-d\Psi_{(\xi,0)}e_{2l}]\\[5mm]
\ds\quad+[J_{\Psi(\xi,\la)}-J_{\Psi(\xi,0)}]\cdot d\Psi_{(\xi,\la)}\quad .
\end{array}
\ee
Using the estimates we used in the above lines, (\ref{VIII.12}) becomes for $1<l\le p$
\be
\label{VIII.13}
|d\Psi\ti{J}_{(\xi,\la)}e_{2l}+d\Psi_{(\xi,\la)} e_{2l-1}|\le O(r_{\Psi(\xi,\la)})\quad .
\ee
Thus combining (\ref{VIII.11}) and (\ref{VIII.13}), we obtain that 
\be
\label{VIII.14}
\forall (\xi,\la)\quad\quad\mbox{ s. t. }|\la|\le r_{\Psi(\xi,\la)}\quad |\ti{J}_{(\xi,\la)}-J_0|\le
 K\ r_{\Psi(\xi,\la)}\quad.
\ee
\section{The unique continuation argument.}
\reset

In this part we show that, assuming ${\mathcal P}_{Q-1}$ (recall that the definition is given by (\ref{zz}),
a point $x_0$ in ${\mathcal C}_{Q}\setminus{\mathcal C}_{Q-1}$  is isolated in 
${\mathcal C}_{Q}\setminus{\mathcal C}_{Q-1}$ unless all points in ${\mathcal C}_\ast$ in a neighborhood
from $x_0$ ar of multiplicity $Q$. This fact has been already proved in chapter IV in the case where $C_{0,x_0}$
was not $Q$ times a same flat holomorphic disk (the easy case). Here we assume that we are in the difficult case
$C_{0,x_0}=Q[D_0]$ and we adopt the coordinate system about $x_0$ constructed in section VIII.2.
We denote by $\Pi$ the map that assigns the first complex coordinate $\xi=\xi_1+i\xi_2$.
Assuming there exists a sequence of points $x_n\in{\mathcal C}_{Q}\setminus{\mathcal C}_{Q-1}$ different from
$x_0$ and converging to $x_0$, the goal of this chapter is to show that $C$ in a neighborhood is a $Q$ times
a same graph. The strategy is inspired by \cite{Ta} chapter 1 : in our coordinates, the points in
${\mathcal C}_{Q}\setminus{\mathcal C}_{Q-1}$ are contained in the disk $\la_i=0$ for $i=1,\cdots,2p-2$
and we shall use a unique continuation argument technique based on the proof of some Carleman estimate to show
that our assumption imposes that the whole cycle in the neighborhood of $x_0$ is included in that disk.
Let $(\xi_n,0)$ be the coordinates of $x_n\rightarrow x_0$. We can allways extract a subsequence such that
$|\xi_{n+1}|\le |\xi_n|^2$. We then introduce the function $g_N(\xi):=\prod_{j=1}^N(\xi-\xi_n)$. Because of 
the speed of convergence of our sequence $\xi_n$ to zero it is not difficult to check that there exists
 a constant $K$ independent of $N$ such that for any $\xi\in B^2_{\rho_{x_0}}$ the following holds
\be
\label{VIII.15}
\frac{K^{-1}}{|\xi|^{N-N_\xi}}\frac{1}{|\xi-\xi_{N_\xi}|}\prod_{j=1}^{N_\xi-1}\frac{1}{|\xi_j|}
\le |g_N^{-1}|(\xi)\le 
\frac{K}{|\xi|^{N-N_\xi}}\frac{1}{|\xi-\xi_{N_\xi}|}\prod_{j=1}^{N_\xi-1}\frac{1}{|\xi_j|}
\ee
where $N_\xi$ is the index less than $N$ such that $|\xi-\xi_{N_\xi}|$ is minimal among the $|\xi-\xi_n|$.
It is also straightforward to check that
\be
\label{VIII.16}
\begin{array}{l}
\ds|\nabla g^{-1}_N|(\xi)\le 
\frac{K(N-N_\xi)}{|\xi|^{N-N_\xi+1}}\frac{1}{|\xi-\xi_{N_\xi}|}\prod_{j=1}^{N_\xi-1}\frac{1}{|\xi_j|}\\[5mm]
\ds +\frac{K}{|\xi|^{N-N_\xi}}\frac{1}{|\xi-\xi_{N_\xi}|^2}\prod_{j=1}^{N_\xi-1}\frac{1}{|\xi_j|}
+\frac{K}{|\xi|^{N-N_\xi}}\frac{1}{|\xi-\xi_{N_\xi}|}\sum_{l=1}^{N_\xi-1}
\frac{1}{|\xi_l|}\prod_{j=1}^{N_\xi-1}\frac{1}{|\xi_j|}
\end{array}
\ee
and we have a corresponding estimate for $|\nabla^kg^{-1}_N|(\xi)$ for arbitrary $k$ in general.
Let $\xi\in \Pi({\mathcal C}_{Q-1})$. $\xi$ belong to some $B^2_{\rho_i}(\xi_i)$ of the covering constructed
in chapter VI. To every such a $i$ we assign $k_i$ an index such that $|\xi_{k_i}|+\rho_{k_i}\le|\xi_i|-\rho_i$
and such that $|\xi_{k_i}-\xi_i|\le K\rho_i$ and such that $\rho_i$ and $\rho_{k_i}$ are comparable :
\be
\label{VIII.17}
K^{-1}\rho_i\le\rho_{k_i}\le K\rho_i\quad.
\ee
This is always possible due to the Whitney-Besicovitch nature of our covering, moreover for every $k$ there exists
a uniformly bounded number of $i$ such that $k_i=k$. Observe also, because of the relative Lipshitz estimate
(\ref{IV.1b}) with constant $\ep$ and because of the ``splitting stage'' of $C_{\xi_i,\rho_i}$ characterised by
(\ref{V.aa3}) we have that for any $\delta>0$ one may choose $\ep$ small enough compare to $\al$ defined in 
chapter V such that for any $\xi\in \Pi({\mathcal C}_{Q-1})$
\be
\label{VIII.18}
dist\lf(\xi,\Pi({\mathcal C}_{Q}\setminus{\mathcal C}_{Q-1})\rg)\ge \delta^{-1} \rho_i\quad,
\ee
where $\xi\in B^2_{\rho_i}(\xi_i)$. Combining (\ref{VIII.17}) and (\ref{VIII.18}) we get that
\be
\label{VIII.19}
\begin{array}{l}
\ds\forall i\in I\quad\forall \xi\in B_{\rho_i}^2(\xi_i)\quad\forall\zeta\in B^2_{\rho_{k_i}}(\xi_{k_i})\\[5mm] 
\ds\quad  \frac{1}{2}dist\lf(\xi,\Pi({\mathcal C}_{Q}\setminus{\mathcal C}_{Q-1})\rg)
\le dist\lf(\zeta,\Pi({\mathcal C}_{Q}\setminus{\mathcal C}_{Q-1})\rg) 
\le 2 dist\lf(\xi,\Pi({\mathcal C}_{Q}\setminus{\mathcal C}_{Q-1})\rg)
\end{array}
\ee
From (\ref{VIII.15}), (\ref{VIII.16}) and (\ref{VIII.19}) we get that for any $N\in {\N}$,
for any $\xi\in B^2_{\rho_i}(\xi_i)$
and for any $\zeta\in B^2_{\rho_{k_i}}(\xi_{k_i})$ 
\be
\label{VIII.20}
|g^{-1}_N|(\xi)\le K |g^{-1}_N|(\zeta)\quad .
\ee
and
\be
\label{VIII.21}
\rho_i\frac{|\nabla g_N|}{|g_N|^2}(\xi)+\rho_i^2\frac{|\nabla^2 g_N|}{|g_N|^2}(\xi)\le K\frac{1}{|g_N|}(\zeta)
\quad.
\ee

Let $\chi_{\rho_{x_0}}$ be a cut-off function identically equal to 1 in $B^2_{\rho_{x_0}/2}(0)$
and equal to 0 outside $B^2_{\rho_{x_0}}$. In $B^2_{\rho_{x_0}}(z_0)\times{\R}^{2p-2}$ we introduce the cycle
$C^{g_N}$ which is given by
\be
\label{VIII.22}
\forall \xi\in B^2_{\rho_{x_0}}\quad\lf<C^{g_N},\Pi,\xi\rg>:={g_N^{-1}(\xi)}_\ast\lf<C,\Pi,\xi\rg>
\ee
In other words if $\Psi\ :\ \Sigma\rightarrow B^2_{\rho_{x_0}}(0)\times {\R}^{2p-2}$ is a parametrisation
of a piece of $C$, a parametrisation of the corresponding piece in $C^{g_N}$ is given by 
$(\Psi_\xi,g^{-1}_N\circ\Pi\circ\Psi\ \Psi_\la)$ where $(\Psi_\xi,\Psi_\la)$ are the coordinates of $\Psi$.
Since $C^{g_N}$ is a cycle in $B^2_{\rho_{x_0}}(z_0)\times{\R}^{2p-2}$, we have, denoting $\Om:=\sum_{l=1}^{p-1}d\la_{2l-1}\wedge d\la_{2l})$
\be
\label{VIII.23}
C^{g_N}(\chi_{\rho_{x_0}}\circ\Pi\ \Om)=
-C^{g_N}(d\chi_{\rho_{x_0}}\circ\Pi \wedge\sum_{l=1}^{p-1}\la_{2l-1} d\la_{2l})\quad .
\ee
We split $C^{g_N}(\chi_{\rho_{x_0}}\circ\Pi\ \Om)=
C^{g_N}\res B^2_{\rho_{x_0}/2}\times{\R}^{2p-2}(\Om)+C^{g_N}
\res(B^2_{\rho_{x_0}}\setminus B^2_{\rho_{x_0}/2})
\times{\R}^{2p-2}(\chi_{\rho_{x_0}}\circ\Pi\ \Om)
$ and we have
\be
\label{VIII.24}
C^{g_N}\res B^2_{\rho_{x_0}/2}\times{\R}^{2p-2}(\Om)=
\sum_{i\in I}C^{g_N}\res B^2_{\rho_{x_0}/2}\times{\R}^{2p-2}(\vphi_i\circ\Pi\sum_{l=1}^{p-1}\Om) \quad,
\ee
where we recall that the partition of unity was constructed in (\ref{V.ax14}) adapted to the covering
$B^2_{\rho_i}(\xi_i)$. 
Let $\Psi_i\ :\ \Sigma_i\rightarrow B^2_{2\rho_i}(\xi_i)\times B^{2p-2}_{2\rho_i}(0)$ be a smooth parametrisation
of $C\res B^2_{2\rho_i}(\xi_i)\times B^{2p-2}_{\rho_{x_0}}(0)$ and denote by $\eta_i$ the map from $\Sigma_i$
into ${\R}^{2p}$  given by \cite{Ri3} - Proposition A.3 - such that $\Psi_i+\eta_i$ is $J_0-$holomorphic.
Since $J$ in $B^{2p}_{2\rho_i}((\xi_i,0))$  is closed to $J_0$ at a distance comparable to $\rho_i$ - see
(\ref{VIII.14}) - we have
\be
\label{VIII.25}
\|\nabla\eta_i\|_{L^2(\Sigma_{\frac{3}{2},i})}\le \rho_i^2\quad
\ee
where we recall that $\Sigma_{\frac{3}{2},i}=\Sigma_i\cap\Psi_i^{-1}(B^2_\frac{3\rho_i}{2}(\xi_i)\times{\R}^2)$.
Using now lemma II.2 of \cite{Ri3}, which do not require $J$ to be $C^1$ in these coordinates but just the metric 
$g$ to be close to the flat one, one has
\be
\label{VIII.25b}
\|\eta_i\|_{L^2(\Sigma_{\frac{4}{3},i})}\le \rho_i^3\quad
\ee
Using the parametrisation $\Psi_i=(\Psi_{i,\xi},\Psi_{i,\la})$, we have
\be
\label{VIII.26}
\begin{array}{l}
\ds C^{g_N}\res B^2_{\rho_{x_0}/2}\times{\R}^{2p-2}(\vphi_i\circ\Pi\ \Om)\\[5mm]
\ds =\int_{\Sigma_i}\vphi(\Psi_{i,\xi}) \sum_{l=1}^{p-1}d\lf[\frac{\Psi_{i,\la}^{2l-1}}{g(\Psi_{i,\xi})}\rg]
\wedge d\lf[\frac{\Psi_{i,\la}^{2l}}{g(\Psi_{i,\xi})}\rg]
\end{array}
\ee
We compare this quantity with 
\be
\label{VIII.27}
\begin{array}{l}
\ds C^{g_N}_i\res (\vphi_i\circ\Pi\ \Om)\\[5mm]
\ds:=
\int_{\Sigma_i}\vphi(\Psi_{i,\xi}) \sum_{l=1}^{p-1}d\lf[
\frac{\Psi_{i,\la}^{2l-1}+\eta_{i,\la}^{2l-1}}{g(\Psi_{i,\xi}+\eta_{i,\xi})}\rg]
\wedge d\lf[\frac{\Psi_{i,\la}^{2l}+\eta_{i,\la}^{2l}}{g(\Psi_{i,\xi}+\eta_{i,\xi})}\rg]\quad .
\end{array}
\ee
One has
\be
\label{VIII.28}
\begin{array}{l}
\ds |(C^{g_N}_i-C^{g_N})(\vphi_i\circ\Pi\sum_{l=1}^{p-1}d\la_{2l-1}\wedge d\la_{2l})|\le \\[5mm]
\ds K\int_{\Sigma_i}
\lf|\nabla\lf[\frac{\Psi_{i,\la}+\eta_{i,\la}}{g(\Psi_{i,\xi}+\eta_{i,\xi})}-
\frac{\Psi_{i,\la}}{g(\Psi_{i,\xi})}\rg]\rg|\ 
\lf|\nabla\lf[\frac{\Psi_{i,\la}+\eta_{i,\la}}{g(\Psi_{i,\xi}+\eta_{i,\xi})}\rg]\rg|\\[5mm]
\ds\le \delta K\int_{\Sigma_i}\vphi_{k_i}\lf|\nabla\lf[\frac{\Psi_{i,\la}+\eta_{i,\la}}{g(\Psi_{i,\xi}+\eta_{i,\xi})}\rg]\rg|^2
\\[5mm]
\ds+\frac{K}{\delta}\int_{\Sigma_i}\lf|\nabla\lf[\frac{\Psi_{i,\la}+\eta_{i,\la}}{g(\Psi_{i,\xi}+\eta_{i,\xi})}-
\frac{\Psi_{i,\la}}{g(\Psi_{i,\xi})}\rg]\rg|^2
\end{array}
\ee
We have
\be
\label{VIII.29}
\begin{array}{l}
\ds\int_{\Sigma_i}\lf|\nabla\lf[\frac{\Psi_{i,\la}+\eta_{i,\la}}{g(\Psi_{i,\xi}+\eta_{i,\xi})}-
\frac{\Psi_{i,\la}}{g(\Psi_{i,\xi})}\rg]\rg|^2\\[5mm]
\ds\le\int_{\Sigma_i}\lf|\nabla\lf[\Psi_{i}\lf(\frac{1}{g(\Psi_{i,\xi}+\eta_{i,\xi})}-\frac{1}{g(\Psi_{i,\xi})}\rg)
\rg]\rg|^2\\[5mm]
\ds+\int_{\Sigma_i}{|\nabla\eta_i|}sup_{\xi\in B^2_{\rho_i}(\xi_i)}\frac{1}{|g(\xi)|^2}+
\int_{\Sigma_i}|\eta_i|^2\ sup_{\xi\in B^2_{\rho_i}(\xi_i)}\frac{|\nabla g|^2}{|g|^4}(\xi)\quad.
\end{array}
\ee
Let $f(\xi)$ be the flat norm of the slice of $C^{g_N}$ by $\Pi^{-1}(\xi)$,
\[
f(\xi)={\mathcal F}(<C^{g_N},\Pi,\xi>)\quad.
\]
Using  (\ref{V.b3}),(observe that the difference of the densities for the metrics $\om(\cdot,J\cdot)$
and $\om_0(\cdot,J_0\cdot)$ is as small  as we want for $\rho_{x_0}$ chosen small enough,
using also (\ref{VIII.20}) and (\ref{VIII.21}), we have
\be
\label{VIII.30}
\begin{array}{l} 
\ds\int_{\Sigma_i}{|\nabla\eta_i|^2}sup_{\xi\in B^2_{\rho_i}(\xi_i)}\frac{1}{|g(\xi)|^2}+
\int_{\Sigma_i}|\eta|^2\ sup_{\xi\in B^2_{\rho_i}(\xi_i)}\frac{|\nabla g_N|^2}{|g_N|^4}(\xi)\\[5mm]
\ds\le K\int_{B^2_{\rho_{k_i}}(\xi_{k_i})}|f|^2
\end{array}
\ee
where we have used the fact that $\int_{\Sigma_{k_i}}\vphi_{k_i}{|\Psi_{k_i,\la}|^2}\ge K \rho_i^4$. This lower bound is a crucial point
in our paper it comes from the fact that $C$ restricted to $B_{\rho_{k_i}}(x_{k_i})$ is splitted : we have that 
$\rho_{k_i}^{-2}M(C\res B_{\rho_{k_i}}(x_{k_i}))$ is less that $\pi Q-K_0\al$ (see (\ref{VI.1c})) where $K_0$ and $\al$ only depend
on $p$, $Q$, $J$ and $\om$. If $\Psi_{k_i,\la}$ would  have been too close to $0$ in $L^2$ norm, since the intersection number
between $(\Psi_{k_i})_\ast[\Sigma_{k_i}]$ and the $2p-2$-planes $\Pi^{-1}(\xi)$ for $\xi\in B_{\rho_{k_i}}^2(\xi_{k_i})$ is $Q$,
 $\rho_{k_i}^{-2}M(C\res B_{\rho_{k_i}}(x_{k_i}))$ would have been too large would have contradicts the upper bound (\ref{VI.1c}).
The first term in the right-hand-side of (\ref{VIII.29}) 
can be bounded as follows
\be
\label{VIII.31}
\begin{array}{l}
\ds\int_{\Sigma_i}\lf|\nabla\lf[\Psi_{i}\lf(\frac{1}{g_N(\Psi_{i,\xi}+\eta_{i,\xi})}-\frac{1}{g_N(\Psi_{i,\xi})}\rg)
\rg]\rg|^2\\[5mm]
\ds\le \int_{\Sigma_i}\lf|\frac{1}{g_N(\Psi_{i,\xi}+\eta_{i,\xi})}-\frac{1}{g_N(\Psi_{i,\xi})}\rg|^2
+\int_{\Sigma_i}\rho_i^2
\lf|\nabla\lf(\frac{1}{g_N(\Psi_{i,\xi}+\eta_{i,\xi})}-\frac{1}{g_N(\Psi_{i,\xi})}\rg)\rg|^2
\\[5mm]
\ds\le
K\int_{\Sigma_i}\rho_i^4\lf|\sup_{\xi\in B^2_{\rho_i}(\xi_i)}\frac{|\nabla g_N|}{|g_N|^2}\rg|^2(\xi)\\[5mm]
\ds+K\int_{\Sigma_i}\rho_i^2|\nabla\eta_i|^2\ 
\lf|\sup_{\xi\in B^2_{\rho_i}(\xi_i)}\frac{|\nabla g_N|}{|g_N|^2}\rg|^2(\xi)\\[5mm]
\ds+K\int_{\Sigma_i}\rho_i^2|\eta_i|^2\sup_{\xi\in B^2_{\rho_i}(\xi_i)}\lf|\frac{|\nabla g_N|^2}{|g_N|^3}\rg|^2(\xi)\\[5mm]
\ds+K\int_{\Sigma_i}\rho_i^2|\eta_i|^2\sup_{\xi\in B^2_{\rho_i}(\xi_i)}\lf|\frac{|\nabla^2 g_N|}{|g_N|^2}\rg|^2(\xi)
\end{array}
\ee
Using (\ref{VIII.20}) and (\ref{VIII.21}) like above and combining (\ref{VIII.28})...(\ref{VIII.31}), we have
finally for every $\delta>0$
\be
\label{VIII.32}
\begin{array}{l}
\ds |(C^{g_N}_i-C^{g_N})(\vphi_i\circ\Pi\ \sum_{l=1}^{p-1}d\la_{2l-1}\wedge d\la_{2l})|\le \\[5mm]
\ds\delta K\int_{\Sigma_i}\vphi_{i}
\lf|\nabla\lf[\frac{\Psi_{i,\la}+\eta_{i,\la}}{g(\Psi_{i,\xi}+\eta_{i,\xi})}\rg]\rg|^2
+\frac{K}{\delta}K\int_{B^2_{\rho_{k_i}}(\xi_{k_i})}|f|^2
\end{array}
\ee
Since $\Psi_i+\eta_i/g\circ\Pi\circ\Psi_i+\eta_i$ is an holomorphic map into ${\C}^p$, the $\la-$coordinates of it,
$\Psi_{i,\la}+\eta_{i,\la}/g\circ\Pi\circ\Psi_i+\eta_i$ is also a holomorphic map but into ${\C}^{p-1}$ and one has
\be
\label{VIII.33}
\begin{array}{l}
\ds\frac{\Psi_{i,\la}+\eta_{i,\la}}{g(\Psi_{i,\xi}+\eta_{i,\xi})}^\ast(\sum_{l=1}^{p-1}d\la_{2l-1}\wedge d\la_{2l})=\frac{1}{2}
(\Psi_{i,\la}+{\eta_{i,\la}})^\ast\lf(\sum_{l=1}^{p-1}d\La_l\wedge d\ov{\La}_l\rg)\\[5mm]
\ds=\lf|\nabla\frac{\Psi_{i,\la}+\eta_{i,\la}}{g(\Psi_{i,\xi}+\eta_{i,\xi})}\rg|^2 \ d\zeta\wedge d\ov{\zeta}
\end{array}
\ee
where $\zeta$ denotes local complex coordinates on $\Sigma_i$, and $\La_l$ is the complex coordinate
$\La_l=\la_{2l-1}+i\la_{2l}$. Therefore, combining (\ref{VIII.32})
and (\ref{VIII.33}) we have for $\delta$ chosen such that $\delta K<\frac{1}{2}$ 
\be
\label{VIII.34}
\begin{array}{l}
\ds C^{g_N} (\vphi_i\circ\Pi\sum_{l=1}^{p-1}d\la_{2l-1}\wedge d\la_{2l})\ge\\[5mm]
\ds\frac{1}{2} \int_{\Sigma_i}\vphi_i\ \lf|\nabla\frac{\Psi_{i,\la}+\eta_{i,\la}}{g(\Psi_{i,\xi}+\eta_{i,\xi})}\rg|^2 
-\frac{K}{\delta}\int_{B^2_{\rho_{k_i}}(\xi_{k_i})}|f|^2
\end{array}
\ee
Let $\{a_i^l(\xi)\}_{l=1\cdots Q}$ be the holomorphic $Q-$valued graph realised by $(\Psi_{i,\xi}+\eta_{i,\xi},
g_N^{-1}(\Psi_{i,\xi}+\eta_{i,\xi})\Psi_{i,\la}+\eta_{i,\la})$, we have that 
\be
\label{VIII.35}
\int_{\Sigma_i}\vphi_i\ \lf|\nabla\frac{\Psi_{i,\la}+\eta_{i,\la}}{g(\Psi_{i,\xi}+\eta_{i,\xi})}\rg|^2 =
\int_{B^2_{\rho_i}}\vphi_i\sum_{l=1}^{p-1}|\nabla a_i^l|^2(\xi)\ d\xi\wedge d\ov{\xi}
\ee
Clearly this quantity is larger than $\int_{B^2_{\rho_i}}\vphi_i|\nabla \ti{f}|^2$ where $\ti{f}(\xi)$ is the
Flat norm of the slice of $C^{g_N}_i$ by $\Pi^{-1}(\xi)$. Replacing $\ti{f}$ by $f$ itself, the slice
of $C^{g_N}$ by $\Pi^{-1}(\xi)$,
in the integral $\int_{B^2_{\rho_i}}\vphi_i|\nabla \ti{f}|^2$, induces error terms which can be controled
by $\int_{B^2_{\rho_{k_i}}(\xi_{k_i})}|f|^2$ like in the computation of the error between $C^{g_N}$ and 
$C^{g_N}_i$ above. Therefore we have
\be
\label{VIII.34z}
\ds C^{g_N} (\vphi_i\circ\Pi\ \Om)\ge
\frac{1}{2} \int_{\Sigma_i}\vphi_i\ |\nabla f|^2 
-\frac{K}{\delta}\int_{B^2_{\rho_{k_i}}(\xi_{k_i})}|f|^2
\ee
Because of the relative lipschitz estimate, $f$ extends as a $W^{1,2}$ function on all of $B^2_{\rho_{x_0}}(0)$.
Standard Poincar\'e estimates yields
\be
\label{VIII.35z}
\int_{B^2_{\rho_{x_0}}}|\chi_{\rho_{x_0}}f|^2\le K\rho_{x_0}^2\int_{B^2_{\rho_{x_0}}}|\nabla(\chi_{\rho_{x_0}}f)|^2
\quad .
\ee
Taking $\rho_{x_0}$ small enough we can ensure that $K\rho_{x_0}^2>O(\delta)$ and combining (\ref{VIII.23}), 
(\ref{VIII.34z}) and (\ref{VIII.35z}) we get that
\be
\label{VIII.36}
\int_{B^2_{\rho_{x_0}/2}}|f|^2\le \int_{B^2_{\rho_{x_0}}\setminus B^2_{\rho_{x_0}/2}}|f|^2+|\nabla f|^2+
C^{g_N}(d\chi_{\rho_{x_0}}\circ\Pi \wedge\sum_{l=1}^{p-1}\la_{2l-1} d\la_{2l})
\ee
By taking the sequence $\xi_n$ such that the largest $|\xi_1|$ 
satisfies $|\xi_1|\le \frac{\rho_{x_0}}{4}$, if $f$ is not identically zero near the origin we would have that
$ \lf(\frac{\rho_{x_0}}{4}\rg)^{-N}\int_{B^2_{\rho_{x_0}/2}}|f|^2$ would tend to infinity, whereas,
it is not difficult to check that the right-hand-side of (\ref{VIII.36}) which involves quantities
supported in $B^2_{\rho_{x_0}}\setminus B^2_{\rho_{x_0}/2}$ is bounded by $K_{\rho_{x_0}}N^2
\lf(\frac{\rho_{x_0}}{2}\rg)^{N}$. The multiplication of it by $\lf(\frac{\rho_{x_0}}{4}\rg)^{-N}$ tends
clearly to zero as $N$ tend 
to infinity. We have then obtained a contradiction and we have proved that any point inside 
${\mathcal C}_{Q}\setminus{\mathcal C}_{Q-1}$ is surrounded in ${\mathcal C}_\ast$ by points which are
all in ${\mathcal C}_{Q}\setminus{\mathcal C}_{Q-1}$ or by points which are all in
 ${\mathcal C}_{Q-1}$. It remains to show that a point in ${\mathcal C}_{Q}\setminus{\mathcal C}_{Q-1}$
is not an accumulation point of $\cup_{q\le Q-1}Sing^q$. This is the purpose of the next chapter.

\section{Points in ${\mathcal C}_{Q}\setminus{\mathcal C}_{Q-1}$ are not accumulation points of 
$\cup_{q\le Q-1}Sing^q$.}
\reset

In this chapter we prove, Assuming ${\mathcal P}_{Q-1}$, that points in ${\mathcal C}_{Q}\setminus{\mathcal C}_{Q-1}$ are not accumulation points
of $\cup_{q\le Q-1}Sing^q$ and combining this fact with the result in the previous chapter we will have proved ${\mathcal P}_Q$.

Let then $x_0\in {\mathcal C}_{Q}\setminus{\mathcal C}_{Q-1}$, and assume then that $x_0$ is an accumulation point of ${\mathcal C}_{Q-1}$,
which means, using the monotonicity formula, lemma~\ref{lm-IV.1a} together with the result obtained in  the previous chater, that there exists a 
radius $\rho$ such that 
${\mathcal C}_\ast\cap B_\rho(x_0)\subset{\mathcal C}_Q$ and that $({\mathcal C}_{Q}\setminus{\mathcal C}_{Q-1})\cap
B_\rho(x_0)=\{x_0\}$. 
From the assumed hypothesis ${\mathcal P}_{Q-1}$, we have then that there exists a Riemann surface $\Sigma$ and a smooth $J-$holomorphic
map $\Psi$ such that $C\res B_{r_0}(x_0)=\Psi_\ast[\Sigma]$. The goal is to show that $\Sigma$ has a finite topology and that
it is a closed Riemann Surface. The idea is to perturb $\Psi$ by finding $\eta\in L^\infty(\Sigma)$ such that $\Psi+\eta$ is
$J_0-$holomorphic and $(\Psi+\eta)_\ast[\Sigma]$ is a cycle. 

For any $r<r_0$, we denote $\Sigma_r$ the finite Riemann surface obtained by taking $\Sigma\cap\Psi^{-1}(B_\rho(x_0)\setminus B_r(x_0))$
and we shall denote $\Gamma_r$ the part of the boundary of $\Sigma_r$ which is disjoint from $\p \Sigma\subset (|\Psi-x_0|)^{-1}(r_0)$.
On $\Sigma_{r}$ we consider $\eta_r$ the map which is given by Proposition A.3 in \cite{Ri3}. It satisfies in particular, using
the complex coordinates induced by $J_0$,

\be
\label{IX.1}
\begin{array}{l}
\ds \ov{\p}(\Psi+\eta_r)=0\quad\quad\mbox{ in }\Sigma_r\\[5mm]
\ds\forall r\le r_0\quad\quad \int_{\Sigma_r}|\nabla\eta_r|^2\le\int_{\Sigma_r}|J(\Psi)-J_0|^2|\nabla\Psi|^2\le K r_0^4\quad,
\end{array}
\ee
where we have used that, for the induced metric by $\Psi$ on $\Sigma$, ($\Psi$ is an isometry) we have
$\int_{\Sigma}|\nabla\Psi|^2=M(C\res B_{r_0}(x_0))\le K r_0^2$.  
Using local $\xi_1$ $\xi_2$ coordinates in $\Sigma_{r}$, we have
for all $k=1\cdots 2p$
\[
\frac{\p\Psi^k_i}{\p \xi_1}=-\sum_{l=1}^{2p}J^k_l(\Psi_i)\frac{\p\Psi^l}{\p \xi_2}\quad\quad\mbox{and}
\quad\quad\frac{\p\Psi^k_i}{\p \xi_2}=\sum_{l=1}^{2p}J^k_l(\Psi_i)\frac{\p\Psi^l}{\p \xi_1}
\]
Taking respectively the $\xi_1$ derivative and the $\xi_2$ derivative of these two equations we obtain
\be
\label{IX.2}
\forall k=1\cdots 2p\quad\Delta_{\Sigma_r}\Psi^k_i=\ast\lf(\sum_{l=1}^{2p}d(J^k_l(\Psi_i))\wedge d\Psi_i^l
\rg)\quad.
\ee
From (\ref{IX.1}) we deduce that $\Delta_{\Sigma_r}(\Psi+\eta_r)=0$ therefore this yields
\be
\label{IX.3}
\forall k=1\cdots 2p\quad\Delta_{\Sigma_{r}}\eta^k_r=-\ast\lf(\sum_{l=1}^{2p}d(J^k_l(\Psi_i))\wedge d\Psi_i^l
\rg)\quad.
\ee
Let $\delta_r^k$ given by 
\be
\label{IX.4}
\lf\{
\begin{array}{l}
\Delta_{\Sigma_r}\delta^k_r=\ast\lf(\sum_{l=1}^{2p}d(J^k_l(\Psi_i))\wedge d\Psi_i^l\rg)\quad\quad\mbox{ in }\Sigma_{r}\\[5mm]
\delta_r^k=0\quad\quad\mbox{ on }\p\Sigma_r
\end{array}
\rg.
\ee
From \cite{Ge} and \cite{To} there exists a universal constant $K$  such that
\be
\label{IX.5}
\|\delta_r\|_{L^\infty(\Sigma_r)}+\|\nabla\delta\|_{L^2(\Sigma_r)}\le K\|J\|_{C^1}\ \int_{\Sigma_r}|\nabla\Psi|^2\le K r_0^2\quad.
\ee
Because of the above estimates, taking some sequence $r_n\rightarrow 0$, one can always extract a subsequence $r_{n'}$ such that
$\eta_{r_{n'}}$ and $\delta_{r_{n'}}$ converge to limits $\eta_0$ and $\delta_0$ that satisfy in particular
\be
\label{IX.6}
\begin{array}{l}
\ds \ov{\p}(\Psi+\eta_0)=0\quad\quad\mbox{ in }\Sigma\\[5mm]
\ds\Delta_{\Sigma}(\eta_0+\delta_0)=0\quad\quad\mbox{ in }\Sigma\\[5mm]
\ds\|\nabla\delta\|_{L^2(\Sigma)}+\|\delta\|_{L^\infty(\Sigma)}+\|\nabla\delta\|_{L^2(\Sigma)}\le K r_0^2
\end{array}
\ee
For any $k=1,\cdots,2p$ we consider the harmonic function $u^k:=\eta^k+\delta^k$. Using the coarea formula we have for any $r<r_0$
\be
\label{IX.7}
\int_0^r ds\int_{\Gamma_s}|\nabla u^k|=\int_{\Sigma\setminus\Sigma_r}|\nabla u^k|\ |\nabla|\Psi||\le 
r\lf(\int_{\Sigma\setminus\Sigma_r}|\nabla u^k|^2\rg)^\frac{1}{2}\quad.
\ee
Therefore using a mean formula, for any $\ep>0$ there exists $s>0$ such that
\be
\label{IX.8}
\int_{\Gamma_s}|\nabla u^k|\le\ep\quad .
\ee
We have
\be
\label{IX.8a}
0=\int_{\Sigma_s}\Delta_\Sigma u^k=\int_{\Gamma_s}\frac{\p u^k}{\p\nu}+\int_{\p\Sigma}\frac{\p u^k}{\p\nu}
\ee
By choosing $\ep$ smaller and smaller and taking the corresponding $s$ given by (\ref{IX.8}),
one gets
\be
\label{IX.8b}
\int_{\p\Sigma}\frac{\p u^k}{\p\nu}\quad .,
\ee
Let $m<M$ two values such that $sup_{\p\Sigma}u^k< m$ and consider the truncation $T_m^Mu^k$ equal tu $m$ if $u^k\le m$
equal to $M$ if $u^k\ge M$ and equal to $u^k$ otherwise. We have
\be
\label{IX.9}
0=\int_{\Sigma_s} T_m^Mu^k\ \Delta_{\Sigma}u^k=-\int_{\Sigma_s}|\nabla T_m^Mu^k|^2+\int_{\Gamma_s}T_m^Mu^k\frac{\p u^k}{\p\nu}+m\int_{\p\Sigma}\frac{\p u^k}{\p\nu}
\ee
Therefore
\be
\label{IX.10}
\int_{\Sigma_s}|\nabla T_m^Mu^k|^2\le M\int_{\Gamma_s}|\nabla u^k|
\ee
and by choosing againg $s$ tending to zero according to (\ref{IX.8}), one gets that $T_m^Mu^k$ is identically equal to $m$ and we deduce that $u^k\le m$
in $\Sigma$. Similarly one gets that $u^k$ is bounded from below and the we have proved that $\|u\|_{L^\infty(\Sigma)}<+\infty$. Combining this
fact with (\ref{IX.6}) we have that 
\be
\label{IX.11}
\|\eta_0\|_{L^\infty(\Sigma)}<+\infty\quad.
\ee
Being more carefuhl above by taking eventually $\Sigma_r$ instead of $\Sigma$ for some $r\in [r_0/2,r_0]$, and
using \cite{Ri3} we could have shown that
$\|\eta_0\|_{L^\infty(\Sigma)}<K\ r_0^2$. We claim now that $\p(\phi+\eta)_\ast[\Sigma]=(\phi+\eta)_\ast[\p\Sigma]$, that is :
for any smooth 1-form $\phi$ equal to zero in a neighborhood of $(\phi+\eta)(\p\Sigma)$, one has
\be
\label{IX.12}
\int_\Sigma(\Psi+\eta_0)^\ast d\Psi=0\quad.
\ee
We have
\be
\label{IX.13}
\lf|\int_{\Sigma_s}(\Psi+\eta_0)^\ast d\Psi\rg|=\lf|\int_{\Gamma_s}(\Psi+\eta_0)^\ast\Psi\rg|\le K_\phi \int_{\Gamma_s}|\nabla\Psi+\eta_0|\quad.
\ee
Argueing like for proving (\ref{IX.8}), for any $\ep$ we can find $s$ such that $\int_{\Gamma_s}|\nabla\Psi+\eta_0|\le\ep$ and we then deduce
(\ref{IX.12}). Thus in $B^{2p}_{r_0}(x_0)\setminus \Psi+\eta_0(\p\Sigma)$, $(\phi+\eta)_\ast[\Sigma]$ is a integre rectifiable holomorphic
cycle. Using the results of Harvey-Shiffman and King (\cite{HS} and \cite{Ki}) we have that there exists a compact Riemann surface
with boundary $\Sigma'$ and an holomorphic map $\Psi'$ such that $(\phi+\eta)_\ast[\Sigma]=\Psi'_\ast[\Sigma']$. $(\phi+\eta)(\Sigma)$ is therefore
an holomophic curve - with boundary - in ${\C}^{2p}$. We claim that $\Psi+\eta_0$ is a  holomorphic simple covering of $\Sigma'$. Indeed let $\om_\Sigma$ be the pull-back by $\Psi$ of the symplectic form $\om$ in ${\R}^{2p}$
we have $\int_\Sigma\om_\Sigma=\int_\Sigma\Psi^\ast\om=\int_{\Sigma}|\nabla\Psi|^2\ge \pi Q r_0^2$, because of the monotonicity formula ($x_0\in {\mathcal C}_Q\setminus{\mathcal C}_{Q-1}$). Let $\om_{\Sigma'}$ be the restriction of $\om_0$ to $\Sigma'$. We have $\int_\Sigma(\Psi+\eta_0)^\ast\om_{\Sigma'}=\int_\Sigma(\Psi+\eta_0)^\ast\om_0
=\int_{\Sigma}|\nabla (\Psi+\eta_0)|^2$. Because of (\ref{IX.1}), we have that the holomorphic covering 
$\Psi+\eta_0$ from $\Sigma$ onto $\Sigma'$ satisfies
\be
\label{IX.14}
\lf|\int_{\Sigma}\om_{\Sigma}-\int_{\Sigma}(\Psi+\eta_0)\om_{\Sigma'}\rg|=
o_{r_0}\lf(\int_{\Sigma}\om_{\Sigma}\rg)\quad.
\ee
Therfore, for $r_0$ small enough this covering has to be a simple one and $\Sigma$ is a compact Riemann surface.
$\Psi$ is now a $J-$holomorphic map from a compact Riemann surface $\Sigma$ into $(B^{2p},J)$, it 
is then smooth and $C\res B^{2p}_{r_0}$ is a $J-$holomorphic curve.

\appendix
\section{Appendix}
\reset
\begin{Lma}
\label{lm-V.3}
Let $U$ be an open subset of ${\R}^2$, let $0<\la<1$ and let $(B^2_{r_i}(z_i))_{i\in I}$ a covering of $U$
 which is locally finite : there exists $n\in {\N}$ such that
\be
\label{V.ax1}
\forall z\in U\quad\quad Card\lf\{i\in I\quad\mbox{ s. t. }\quad z\in B^2_{r_{z_i}}(z_i)\rg\}\le N\quad ,
\ee
moreover one assumes that 
\be
\label{V.ax2}
\forall i, j\in I\quad\quad B_{r_i}(z_i)\cap B_{r_j}(z_j)\ne\emptyset\quad\Longrightarrow\quad r_i\ge\la r_j\quad .
\ee
Then there exists $\delta$ and $P\in{\N}$ depending on $\la$ only such that
\be
\label{V.ax3}
\forall i\in I\quad\quad \mbox{Card}\lf\{j\in I\quad
\mbox{ s. t. }B_{r_j}(z_j)\cap B_{(1+\delta) r_i}(z_i)\ne\emptyset\rg\}\le P\quad .
\ee
\end{Lma}
{\bf Proof of lemma~\ref{lm-V.3}. } 
We argue by contradiction. Assume there exists $\delta_n\rightarrow 0$, a sequences of coverings of $U$, 
$(B^2_{r_{n,i}}(z_{n,i})$ for $i\in I$ satisfying
(\ref{V.ax1}) and (\ref{V.ax2}) and a sequence of indices $i_n$ such that
\be
\label{V.ax4}
\mbox{Card}\lf\{j\in I\quad
\mbox{ s. t. }B_{r_{j,n}}(z_{j,n})\cap B_{(1+\delta_n) r_{i_n,n}}(z_{i_n,n})\ne\emptyset\rg\}\rightarrow +\infty
\quad\mbox{ as }n\rightarrow +\infty
\ee
After an eventual rescaling of the whole covering and a translation we can assume that $r_{i_n,n}=1$ and 
$z_{i_n,n}=0$. After extraction of a subsequence, we can ensure that there exists $A\in \p B_1(0)$ such that for
any $r>0$
\be
\label{V.ax5}
\mbox{Card}\lf\{j\in I\quad
\mbox{ s. t. }B_{r_{j,n}}(z_{j,n})\cap B_r(A)\ne\emptyset\rg\}\rightarrow +\infty
\quad\mbox{ as }n\rightarrow +\infty\quad .
\ee
For a given $r$ and $n$ we take the longuest sequence of distinct balls
of our covering $B_{r_{j_p,n}}(z_{j_p,n})$ for $p=0\cdots P_n$ satisfying 
\begin{itemize}
\item[i)]
\[
B_{r_{j_0,n}}(z_{j_0,n})=B_1(0)
\]
\item[ii)]
\[
\forall p\le P_{n-1}\quad
B_{r_{j_p,n}}(z_{j_p,n})\cap B_{r_{j_{p+1},n}}(z_{j_{p+1},n})\ne\emptyset\quad,
\]
\item[iii)]
\[
\forall p\le  P\quad\quad B_{r_{j_p,n}}(z_{j_p,n})\cap B_r(A)\ne \emptyset\quad .
\]
\end{itemize}
It is clear that for a given $r$ $P_n\rightarrow +\infty$, indeed if it would not be the case, i.e. 
$P_n\le P_\ast<+\infty$
this would imply that the minimal radius for the balls of the covering intersecting $B_r(A)$ is $\la^{P_\ast}$ and 
combining this fact with (\ref{V.ax5}) would contradict (\ref{V.ax1}). Therefore
 we can then find $r_m\rightarrow +\infty$ and ${n_m}\rightarrow +\infty$ as $m\rightarrow +\infty$ and sequences
$(B_{r_{j_p,n_m}})$ for $1\le p\le Q_m$ and for $m=0\cdots +\infty$ such that
\begin{itemize}
\item[i)]
\[
B_{r_{j_0,n_m}}(z_{j_0,n_m})=B_1(0)
\]
\item[ii)]
\be
\label{V.ax6}
\forall p\le Q_{m-1}\quad
B_{r_{j_p,n_m}}(z_{j_p,n_m})\cap B_{r_{j_{p+1},n_m}}(z_{j_{p+1},n_m})\ne\emptyset\quad,
\ee
\item[iii)]
\be
\label{V.ax7}
\forall p\le  Q_m\quad\quad B_{r_{j_p,n_m}}(z_{j_p,n_m})\cap B_{r_m}(A)\ne \emptyset\quad .
\ee
\item[iv)]
\[
Q_m\rightarrow +\infty\quad .
\]
\end{itemize}
Since $\la\le r_{j_1,n_m}\le\la^{-1}$ and since the distance $|z_{j_1,n_m}|$ is bounded, we can extract
from $n_m$ a subsequence that we still denote $n_m$ such that $B_{r_{j_1,n_m}}(z_{j_1,n_m})$ converges to
a limiting ball $B_{r_{1,\infty}}(z_{1,\infty})$ with $\la\le r_{1,\infty}\le\la^{-1}$, $z_{1,\infty}\le 2$
and $A\in \ov{B_{r_{1,\infty}}(z_{1,\infty})}$. This procedure can be iterated and using a diagonal
argument we can assume that 
\[
\forall p\in {\N}\quad r_{j_p,n_m}\rightarrow r_{p,\infty}\quad z_{j_p,n_m}\rightarrow z_{p,\infty}
\]
such that
\[
\forall p\in{\N}\quad\la^{p}\le r_{p,\infty}\le \la^{-p}\quad |z_{p,\infty}|\le 2\quad
 \mbox{ and }A\in\ov{B_{r_{p,\infty}}(z_{p,\infty})}
\]
Moreover because of (\ref{V.ax1}) we have that  
\be
\label{V.ax8}
\forall z\in {\R}^2\quad\quad
 Card\lf\{p\in {\N}\quad\mbox{ s. t. }\quad z\in B^2_{r_{p,\infty}}(z_{p,\infty})\rg\}\le
 N\quad .
\ee
Because of this later fact, since $A\in\ov{B_{r_{p,\infty}}(z_{p,\infty})}$ for all $p$ it is clear that
$r_{p,\infty}\rightarrow +\infty$ as $p\rightarrow +\infty$. Because of (\ref{V.ax8}) again, the 
number of open balls $B_{r_{p,\infty}}(z_{p,\infty})$ containing $A$ is bounded by $N$ and we can therefore
forget them while considering the sequence and assume that 
\[
\forall p\quad\quad A\in \p B^2_{r_{p,\infty}}(z_{p,\infty})\quad .
\]  
Let $\vec{t}_p(A)\in S^1$ be the unit exterior normal to $\p B^2_{r_{p,\infty}}(z_{p,\infty})$ at $A$. 
Let $\vec{t}_\infty$ be an accumulation unit vector of the sequence $\vec{t}_p(A)$. Given a direction
$\vec{t}$ and an open disk containing $A$ in it's boundary and whose exterior normal at $A$ is given by $\vec{t}$
any other open disk containing $A$ in it's boundary and whose exterior unit at $A$ is not $-\vec{t}$ as a non
empty intersection with that disk. Taking $B^2_{r_{p_0,\infty}}(z_{p_0,\infty})$ such that 
$\vec{t}_{p_0}(A)\ne -\vec{t}_\infty$, there exists then infinitely many disks $B^2_{r_p,\infty}(z_{p,\infty})$
having a non empty intersection with $B^2_{r_{p_0,\infty}}(z_{p_0,\infty})$ but then, because of (\ref{V.ax2}) that
passes to the limit, all these infinitely many disks have a radius which is bounded from below by a positive number
and this contardicts the fact that $r_{p,\infty}\rightarrow +\infty$ implied by (\ref{V.ax8}).
Thus lemma~\ref{lm-V.3} is proved.\cqfd


\begin{thebibliography}{99}
\bibitem[All]{All} W.K. Allard, ``On the first variation of a varifold'', Ann. Math., 95 (1972), 417-491.
\bibitem[Alm]{Alm} F. Almgren ``Almgren's Big Regularity Paper''
World Sci. Mono. Series in Math. vol. 1, Ed. V.Scheffer, J. Taylor, World Scientific, 2000.
\bibitem[Ch]{Ch} S.X.-D. Chang ``Two dimensional area minimizing integral currents are classical
minimal surfaces'' J. A.M.S., 1, no 4, (1988), 699-778. 
\bibitem[Do]{Do} S.K. Donaldson ``Infinite determinants, stable bundles and curvature'' Duke Math. J. 54 (1987), 231-247.
\bibitem[DNF]{DNF} B.Dubrovine, A.Fomenko and S.Novikov ``Modern Geometry - Methods and applications'' Part I, Springer, GTM 93, (1992).
\bibitem[FK]{FK} H.M. Farkas and I. Kra ``Riemann Surfaces'',  GTM 71, Springer Verlag, 1991.
\bibitem[Fe]{Fe} H. Federer ``Geometric Measure Theory'', Springer-Verlag, 1969.
\bibitem[Ge]{Ge} Y. Ge ``Estimations of the best constant involving the $L^2$ norm in Wente's inequality
and compact $H-$Surfaces in Euclidian space. C.O.C.V., 3, (1998), 263-300.
\bibitem[Gi]{Gi} M. Giaquinta ``Multiple Integrals in the calculus of variations and nonlinear elliptic
systems'' Annals of Math. Studies, 105, Princeton University Press, (1983).
\bibitem[HS]{HS} F.R.Harvey and B.Shiffman ``A characterization
of holomorphic chains'' Ann. of Math., 99, (1974), 553-587.
\bibitem[Ki]{Ki} J.R.King ``The currents defined by analytic
  varieties'', Acta. Math. 127, (1971), 185-220.
\bibitem[Li]{Li} F.H.Lin "Gradient estimates and blow-up analysis for
  stationary harmonic maps" Ann. Math., 149, (1999), 785-829.
\bibitem[PR]{PR} D.Pumberger and T.Rivi\`ere ``The regularity of non calibrated $J-$holomorphic Rectifiable
 2-cycles'' in preparation. 
\bibitem[Re]{Re} E.R. Reifenberg ``An epiperimetric inequality related to the analyticity
of minimal surfaces'' Ann. of Math., {\bf 80} (1964), 1-14.
\bibitem[Ri1]{Ri1} T. Rivi\`ere ``Bubbling and regularity issues in geometric non-linear analysis''
Proceedings of the I.C.M. Beijing, High. Edu. Press, 2002, vol. III, 197-208 .
\bibitem[Ri2]{Ri2} T. Rivi\`ere ``A lower Epiperimetric Inequality for Minimizing
Surfaces'' preprint (2003).
\bibitem[Ri3]{Ri3} T. Rivi\`ere ``Approximating $J-$holomorphic curves by holomorphic ones'' preprint (2003).
\bibitem[RT]{RT} T. Rivi\`ere and G. Tian ``The Singular set of $J-$holomorphic maps into projective algebraic varieties'' to appear in J. Reine Angewan. Math. (2003).
\bibitem[ST]{ST} B. Siebert and G. Tian ``Weierstrass polynomials and plane pseudo-holomorphic curves'' Chinese Ann. Math. Ser. B, 23, (2002), 1-10.
\bibitem[Si]{Si} L. Simon ``Lectures on geometric measure theory'' Proc. Center for Math. Ana., vol. 3, 
Australian National University, 1984.
\bibitem[Ta]{Ta} C.Taubes ``Seiberg-Witten and Gromov Invariants for
  symplectic 4-manifolds'' International Press (2000).
\bibitem[Ti]{Ti} G.Tian "Gauge theory and calibrated geometry, I"
  Ann. Math., 151 (2000), 193-268.
\bibitem[To]{To} P. Topping ``The optimal constant in Wente's $L^\infty$ estimate, Comm. Math. Helv. 72, (1997), 316-328.
\bibitem[UY]{UY} K.Uhlenbeck and S.T.Yau ``On the existence of Hermitian-Yang-Mills connections in stable vector bundles'' C.P.A.M. 39, (1986), S257-S293.
\bibitem[Wh]{Wh} B. White ``Tangent cones to two-dimensional area-minimizing integral currents
are unique'' Duke Math. J., {\bf 50}, no 1, (1983), 143-160.
\end{thebibliography}
\end{document}